\documentclass[a4paper,reqno,11pt]{article}
\usepackage[a4paper,margin=3cm, marginparwidth=2cm]{geometry}
\usepackage[T1]{fontenc}
\usepackage{amsfonts,amsthm,amssymb,amsmath, mathdots, bbm,mathabx, bm}
\usepackage{amssymb, amscd}
\usepackage[makeroom]{cancel}
\usepackage[pdftex]{color,graphicx}
\usepackage[all,cmtip]{xy}
\usepackage{enumitem}
\usepackage{dsfont}
\usepackage{authblk}

\usepackage[pagebackref, colorlinks = true, linkcolor = black, urlcolor  = black, citecolor = red]{hyperref}
\numberwithin{equation}{section}

\newcommand{\un}{\mathbbm{1}}

\newtheorem{theorem}{Theorem}[section]

\newtheorem{proposition}[theorem]{Proposition}
\newtheorem{lemma}[theorem]{Lemma}
\newtheorem{corollary}[theorem]{Corollary}

\newcommand{\be}{\begin{equation}}
\newcommand{\ee}{\end{equation}}

\newcommand{\bE}{{\mathbb{E}}}

\newcommand{\cP}{{\mathcal{P}}}

\newcommand{\Tr}{{\mathrm{Tr}}}

\definecolor{armygreen}{rgb}{0.29, 0.33, 0.13}

\renewcommand\leq\leqslant

\renewcommand\geq\geqslant

\title{From  higher order free cumulants to non-separable hypermaps}
\author{Luca Lionni\footnote{lionni@thphys.uni-heidelberg.de, ORCID:0000-0002-8329-2146.}}
\date{}
\affil{Heidelberg University, Institut für Theoretische Physik, \\Philosophenweg 19, 69120 Heidelberg, Germany.}

\begin{document}

\maketitle
\vspace{-15pt}
\begin{abstract}
Higher order free moments and cumulants, introduced by Collins, Mingo, \'Sniady and Speicher in 2006, describe the fluctuations of unitarily invariant random matrices in the limit of infinite size. The functional relations between their generating functions  were only found last year by Borot, Garcia-Failde, Charbonnier, Leid and Shadrin and a combinatorial derivation is still missing. 
We simplify these relations and show how their combinatorial derivation reduces to the computation of generating functions of planar non-separable hypermaps with prescribed vertex valencies and weighted hyper-edges. The functional relations obtained by Borot et al.~involve some remarkable simplifications, which can be formulated as identities satisfied by these generating functions. The case of third order free cumulants, whose combinatorial understanding was already out of reach, is derived explicitly.  
\end{abstract}

\tableofcontents

\section{Introduction}

First order free probability \cite{Voicu, Spei94, SpeicherReview, NicaSpeicher, Collins03} provides an efficient understanding of expectations of traces of $N\times N$ random matrices $B$ at leading order in the limit of infinite size $N$:
$$
\varphi (b^\lambda) = \lim_{N\rightarrow \infty} \frac {1} N\,  \mathbb E \mathrm{Tr}(B ^\lambda), 
$$
in generalizing the concepts of independence, cumulants, convolution, and so on, to freeness, free cumulants, free convolution... A central role is  played by the free cumulants $\kappa_\lambda[b]$, whose knowledge is equivalent to the asymptotic moments $\varphi (b^\lambda)$ through the free moment-cumulant formulas,  which involve summations over non-crossing partitions and are thus  invertible in the associated lattice. 
Free cumulants are the useful tools for computing the asymptotic spectrum of a sum or product of independent random matrices, for instance. In practice however, such computations are rendered efficient by going to generating functions  of moments and cumulants, and using the functional analog of these free moment-cumulant formulas, allowing for the use of  complex analysis.  

Higher order free probability \cite{Second-order-1, Second-order-2, Second-order-3, CMSS} deals with the fluctuations of random matrices, that is, with the connected correlations of their traces at leading order in the limit $N\rightarrow \infty$: 
$$
\varphi (b^{\lambda_1}, \ldots, b^{\lambda_p}) = \lim_{N\rightarrow \infty} N^{p-2}  k_p\bigl( \mathrm{Tr}(B ^{\lambda_1}), \ldots, \mathrm{Tr}(B ^{\lambda_p})\bigr), 
$$
where $k_p$ is the $p$-th classical cumulant (connected correlation). The role of free cumulants is then played by their higher order counterparts, $\kappa_{\lambda_1, \ldots, \lambda_p}[b]$, and the notion of freeness generalizes to higher order freeness. These quantities are related to the $\varphi (b^{\lambda_1}, \ldots, b^{\lambda_p}) $ by higher-order free moment-cumulant formulas, also invertible, and which involve summations over pairs consisting of a permutation of $n=\sum_i \lambda_i$ elements, and a partition of its blocks - called partitioned permutations, constrained by a leading order condition on their blocks which generalizes the non-crossing condition. Like for the first order, concrete computations would be rendered more efficient within a functional reformulation of these relations, however until recently, only the relation for second order had been derived in \cite{CMSS} in 2006. A year ago, this was settled by \cite{Analytic-higher}, by allowing for all-genus contributions (not sticking to the leading order in $N$), and  going to the Fock space. The formulas giving the generating series of the $\varphi (b^{\lambda_1}, \ldots, b^{\lambda_p})$ for instance involve a sum over bipartite trees with black and white vertices, whose black vertices are weighted by generating functions of higher-order free cumulants and whose white vertices carry differential operators.  A combinatorial understanding of these functional relations is still missing. 

In this paper, we instead work on  a straightforward combinatorial derivation of the functional relations between the generating functions of higher order cumulants and moments. The approach adopted is a direct generalization of that of \cite{CMSS} for the second order case. The first step (Sec.~\ref{sec:tree-structure}) is a reformulation of the higher order moment-cumulant formulas at the coefficient level, splitting the condition on partitioned permutation in independent conditions, to obtain sums over bipartite combinatorial maps (with black and white  vertices), and tree-like partitions of their black vertices  (Sec.~\ref{sec:mom-cumulants-two-instead-four}). This explains the summations over trees of the functional relations of \cite{Analytic-higher}, as clarified in Sec.~\ref{sec:bij-with-trees}.

The relations derived in Sec.~\ref{sec:bij-with-trees} however involve numbers of planar bipartite maps with prescribed vertex valencies. The next step (Sec.~\ref{sec:funct-rel}) is to factorize the bipartite maps into non-crossing partitions, one per white vertex. For such factorized planar bipartite maps, the functional relations can be derived and take the same form as the functional relations of \cite{Analytic-higher} (Sec.~\ref{sec:factorized}). The key step to put the coefficient formulas in such a factorized form is a decomposition similar to that of planar bipartite maps in their non-separable components \cite{Tutte, Brown, Brown-Tutte, CombEnum}, but where the non-separability is required only  for the white vertices of the map, thus justifying the name non-separable hypermaps. This is done in Sec.~\ref{sec:factorization}. Through this procedure however, the weights associated to the black vertices of the bipartite trees involve additional corrections with respect to  the simpler weights found in \cite{Analytic-higher}, implying some important simplifications among these corrections. 

In Sec.~\ref{sec:corrections}, we show that all the corrections can be computed by acting directly on the generating functions of planar non-separable hypermaps with vertices given by a fixed permutation, and whose hyperedges of valency $q$ are weighted by (first order) free cumulants $\kappa_q[b]$. 
With this, the understanding of these corrections and their simplifications relies solely on deriving these generating functions of non-separable hypermaps, a purely combinatorial task (it does not matter anymore that the weights of the hyper-edges are free cumulants). 
The simplifications occurring among the corrections remain to be  fully understood. They must amount to recursive identities satisfied by the generating functions of planar non-separable hypermaps. The case of third-order free cumulants - already out of reach (see the discussion page 24 of  \cite{Analytic-higher}) - is derived explicitly by direct enumeration. A more efficient approach is needed to derive the identities satisfied by the generating functions of non-separable hypermaps (e.g.~bijective \cite{Jacq-Schaeff, DLDRP}), necessary to understand the simplifications occurring at higher order. 
\section*{Acknowledgements}

My work is supported by the European Research Council (ERC) under the European Union’s Horizon 2020 research and innovation program (grant agreement No818066) and by Deutsche Forschungsgemeinschaft (DFG, German Research Foundation) under Germany's Excellence Strategy  EXC-2181/1 - 390900948 (the Heidelberg STRUCTURES Cluster of Excellence). I warmly thank Matteo Maria Maglio for his help with Mathematica. 

\section{Context}

\subsection{Notations}

Let $S_n$ be the group of permutations of $n$ elements and $S_n^*$ the set of permutations different from the identity. For $\sigma\in S_n$,   $\#(\sigma)$ denotes the number of disjoint cycles of $\sigma$  and $\lvert \sigma \rvert$  the minimal number of transpositions required to obtain $\sigma$, satisfying, $ 
\#(\sigma) + \lvert \sigma \rvert = n \;.
$

We denote by $\pi,\pi',\pi_1,\pi_2$ and so on partitions of a set, and $\cP(n)$ the set of all such partitions. The notation $\#(\pi)$ is used for the number of blocks of $\pi$, while  $B\in \pi$ denotes the blocks, and $|B|$ the cardinal of the block $B$. $\le$ signifies the refinement partial 
order: $\pi'\le \pi$ if all the blocks of $\pi'$ are subsets of the blocks of $\pi$. Furthermore,  $\vee$ denotes the joining of partitions: $\pi\vee\pi'$ is the finest partition which is coarser than both $\pi$ and $\pi'$.
$1_n$ and $0_n$ respectively denote the one-block and the $n$ blocks partitions of $\{1,\dots, n\}$.

The partition induced by the cycles of the permutation $\nu$ is denoted by $\Pi(\nu)$, hence $ \#(\Pi(\nu)) = \#(\nu)$. $d_p(\nu)$ denotes the number of cycles of $\nu$ with $p$ elements.

The notation $\lambda\vdash n$  means that $\lambda$ is  a partition of the integer $n$, that is, a multiplet of integers $\lambda_1\ge \ldots \ge\lambda_p$ such that $\sum_{i=1}^p\lambda_i =n$.  The $\lambda_i$ are called the parts of $\lambda$, $p$ is the number of parts of $\lambda$, denoted by $\#(\lambda)$, and $d_i(\lambda)$ is the number of parts of $\lambda$ equal to $i$, so that $n = \sum_{i=1}^n i d_i(\lambda) $ and $p=\#(\lambda) = \sum_{i=1}^n  d_i(\lambda) $.  To a sequence $\{d_i\}_{1\le i \le k}$, we may therefore always associate the partition $\Lambda_{\{d_i\}_{1\le i \le k}}$ of $\sum_{i=1}^k i d_i $, which has $d_i$ parts equal to $i$. 
If $\pi \in \mathcal{P}(n)$, $\Lambda(\pi)$  is the partition of $n$ given by the sizes of the blocks of $\pi$, and if 
$\sigma\in S_n$, $\Lambda(\sigma) = \Lambda(\Pi(\sigma) )$,
 so that for instance $\#(\Lambda(\sigma))=\#(\sigma)$, and the conjugacy class of a permutation $\sigma$ is $C_{\Lambda(\sigma)}$.

If $\lambda\vdash n$, $\gamma_\lambda  =\gamma_{\lambda_1, \ldots, \lambda_p} $ is the permutation
 \vspace{-1ex}
 \be
 \label{eq:gamma-cycles}
\gamma_\lambda   =  \gamma_{\lambda_1, \ldots, \lambda_p}  = (1 2 \cdots \lambda_1)(\lambda_1+1 \cdots \lambda_1+\lambda_2) \ldots (\sum_{i=1}^{p-1}\lambda_i \cdots \sum_{i=1}^{p}\lambda_i ). 
 \ee

\subsection{Classical cumulants}
\label{sub:classical cumulants}
 Consider some random variables $x_1, \ldots, x_p$, then the classical cumulants are defined as:
\be
\label{eq:classical-cumulant}
k_n(x_1, \ldots, x_n) = \sum_{\pi\in \mathcal{P}(n)}  (-1)^{|\pi|-1} (|\pi|-1)!  \prod_{G\in \pi} \bE\Bigl[\prod_{i\in G}x_i\Bigr] \,.
\ee
If all $x_i$ are equal to $x$, we use the notation $k_n(x)$. 
By M\"obius inversion in the lattice of partitions: 
\be
\label{eq:classical-moment-cumulant}
\bE[x_1\cdots x_n]   = \sum_{\pi\in \mathcal{P}(n)} \prod_{G\in \pi}  k_{\lvert G \rvert}(\{x_i\}, i\in G) \,.
\ee

More generally, we can define for $\pi'\in \mathcal{P}(n)$:
\be
\label{eq:classical-cumulant-2}
k_{\pi'}(x_1, \ldots, x_n) =  \prod_{G\in \pi'}  k_{\lvert G \rvert}(\{x_i\}, i\in G)  =  \sum_{\pi, \pi \le \pi'}  \lambda_{\pi, \pi'} \prod_{G\in \pi} \Bigl[\prod_{i\in G}x_i\Bigr] \,,
\ee
and by M\"obius inversion in the lattice of partitions: 
\be
\label{eq:classical-cumulant-20}
\prod_{G\in \pi'} \bE\Bigl[\prod_{i\in G}x_i\Bigr] = \sum_{\pi, \pi \le \pi'} \prod_{G\in \pi}  k_{\lvert G \rvert}(\{x_i\}, i\in G) \,.
\ee

\subsection{Bipartite trees}
\label{sub:graphs}

A \emph{bipartite graph} is given by two sets $B$ and $W$, whose elements represent the black and white vertices, and a set $E$ of pairs $(w,b)$ where $w\in W$ and $b\in B$, whose elements represent the edges (a pair may appear more than once, corresponding to multiple edges linking two vertices). A bipartite graph is a forest if $L= \lvert E\rvert - \lvert B\rvert - \lvert W\rvert + \lvert K\rvert =0$, where $\lvert K\rvert$ is its number of connected components, and it is a \emph{tree} if in addition $\lvert K\rvert=1$. 

\begin{figure}[h!]
\centering
\includegraphics[scale=0.9]{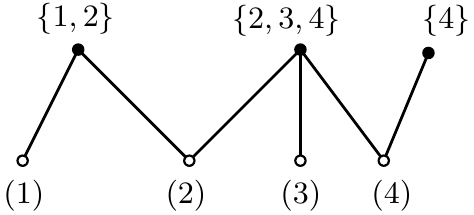}
\caption{An exemple of tree in $\mathcal T_4$, with $\mathcal{I}(T)=\left\{\{1,2\}, \{2,3,4\}, \{4\}\right\}$.
}
\label{fig:ex-arbre}
\end{figure}

 A \emph{labeled bipartite tree} $T$ (Fig.~\ref{fig:ex-arbre}) is a bipartite tree whose white vertices are labeled from 1 to $p=\lvert W\rvert $. It can be encoded by a set $\mathcal{I}(T)$ whose elements are subsets $I\subset \{1, \ldots, p\}$ each representing a black vertex. If $b\in W$, the corresponding set is  
$I_b=\{i\in \{1, \ldots, p\} \mid (b, i)\in E\}$. The condition that $T$ is connected is equivalent to $\bigvee_{I\in \mathcal{I}(T)} \Pi(I) = 1_p$, where $\Pi(I)$ is the partition of $ \{1, \ldots, p\}$ obtained by supplementing $I$ by $1_{p-\lvert I \rvert}$. The condition that $T$ is a tree is translated to 
\be
\label{eq:excess-lab-tree}
\sum_{I\in \mathcal{I}(T)} (\lvert I \rvert - 1) - p + 1 =0.
\ee

Note that this implies that if $I, J \in I\in \mathcal{I}(T)$, $\lvert I \cap J \rvert \le 1$ (i.e.~there are no multiple edges). For  $i\in \{1, \ldots, p\}$, we denote by $\mathcal I_i(T)$ the subset of  $\mathcal{I}(T)$ of sets containing $i$, and 
 $\mathrm{deg}_T(i) = \lvert \mathcal I_i(T) \rvert$ the valency of the white vertex number $i$. It will be useful to introduce the notations $\mathcal I'_i(T)$ and $\mathcal{I}'(T)$ for the subsets of $\mathcal I_i(T)$ and $\mathcal{I}(T)$ which do no contain elements $I$ with only one element. 
 The number of edges incident to the white vertex $i$ and whose other extremity is a black leaf is denoted by $\mathrm{deg}^{{\mathsf{L}}}_{T}(i)$. 
 The number of edges of $T$ can be reformulated as $\lvert E (T)\rvert = \sum_{I\in \mathcal{I}(T)} \lvert I \rvert = \sum_{i=1}^p \mathrm{deg}_T(i)$.
Note that the black leaves incident to a given white vertex $i$ of a labeled bipartite tree are not distinguishable, as they all correspond to the same $I=\{i\}$.  
  
 The set of labeled bipartite trees with $p$  white vertices is 
 \be
 \label{eq:Tp}
 \mathcal{T}_p = \biggl\{\mathcal{I} = \{I_1, \ldots I_k\}, k\ge 1, I_a\subset\{1, \ldots, p\} \mid \bigvee_{I\in \mathcal{I}} \Pi(I) = 1_p \, ; \,  \sum_{I\in \mathcal{I}} (\lvert I \rvert - 1) - p + 1 =0\biggr\}, 
 \ee
 and the subset of  $\mathcal{T}_p$ of labeled trees for which all black vertices have valencies two or higher ($\forall I \in \mathcal{I}(T),\, \vert I\rvert \ge 2$) is denoted by $ \mathcal{G}_p$. By convention, $ \mathcal{G}_1$ has a single element, with one white vertex, no edge and no black vertex.

\subsection{Bipartite maps}
\label{sub:bip-map}

 Two permutations $(\sigma_1, \sigma_2)\in S_n^2$ define a \emph{bipartite map}. They bijectively encode a graph with $n$ labeled edges, embedded on a collection of $\#(\Pi(\sigma_1)\vee\Pi(\sigma_2))$ surfaces, up to orientation preserving homeomorphisms. The genus $g = g(\sigma_1, \sigma_2)$ of that surface is given by the Euler characteristics
\be
\label{eq:genus}
\#(\sigma_1) + \#(\sigma_2)  + \#(\sigma_1\sigma_2) - n = 2\#(\Pi(\sigma_1)\vee\Pi(\sigma_2)) - 2g(\sigma_1, \sigma_2). 
\ee
The graph is bipartite and the vertices of each kind - black or white -  are respectively encoded by the disjoint cycles of $\sigma_1$ and $\sigma_2$, and the edges only link vertices of two different types. Each cycle of $\sigma_i$ is a cyclic sequence of numbers which corresponds to the labeled  edges encountered cyclically when turning around one of the type $i$ vertices, counterclockwise. The number of vertices of each kind are respectively given by $\#(\sigma_1)$ and $\#(\sigma_2)$.  Two consecutive edges around a vertex form a corner. 
The connected components of the complement of the graph on the surface - called faces -  are therefore bijectively encoded by the disjoint cycles of $\sigma_1\sigma_2$, or equivalently of $\sigma_2\sigma_1$, so that in particular, $\#(\sigma_1\sigma_2)=\#(\sigma_2\sigma_1)$, and $g(\sigma_1, \sigma_2)=g(\sigma_2, \sigma_1)$. A bipartite map  $(\sigma_1, \sigma_2)$ with $g( \sigma_1, \sigma_2)=0$ is \emph{planar}.

\begin{figure}[h!]
\centering
\includegraphics[scale=1.1]{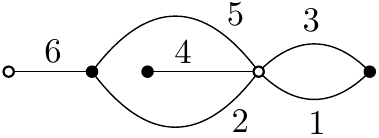}
\caption{An exemple of (planar) bipartite map, with black vertices given by $\sigma_1=(13)(256)(4)$ and white vertices encoded  by $\sigma_2=(13542)(6)$.}
\label{fig:ex-bip-map}
\end{figure}

 The Tutte dual is an operation which consists in adding a grey vertex in each face, and for each new vertex, an edge linking that  vertex to each corner around a white vertex in the face it belongs to. This operation preserves the number of connected components and the genus, and the grey vertices are bijectively encoded by the disjoint cycles of $\sigma_1^{-1}\sigma_2^{-1}$ or $\sigma_2^{-1}\sigma_1^{-1}$.  Said otherwise, 
\be
\label{eq:id-pis}
\#(\Pi(\sigma_1)\vee\Pi(\sigma_2)) = \#(\Pi(\sigma_i)\vee\Pi(\sigma_1\sigma_2)) = \#(\Pi(\sigma_i)\vee\Pi(\sigma_2\sigma_1)), 
\ee
and 
 \be
 \label{eq:id-gs}
 g(\sigma_1, \sigma_2) = g(\sigma_i, \sigma_1^{-1}\sigma_2^{-1}) = g(\sigma_i, \sigma_2^{-1}\sigma_1^{-1}) .
 \ee
 Finally, reversing the convention for the ordering of the edges from counterclockwise to clockwise obviously does not change the genus, so that  $g(\sigma_1, \sigma_2)=g(\sigma_1^{-1}, \sigma_2^{-1})$. 

\ 

\begin{figure}[h!]
\centering
\raisebox{2ex}{\includegraphics[scale=0.9]{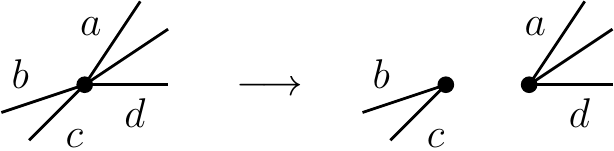}}\hspace{3cm} \includegraphics[scale=1.1]{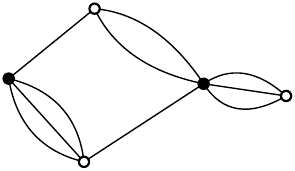}
\caption{Left: splitting the vertex along two corners. Right: a non-separable hypermap.}
\label{fig:Vert-split}
\end{figure}

Consider a bipartite map $(\sigma, \tau)$, $\sigma, \tau \in S_n$, and  one of its vertices (a cycle of either permutation), and two distinct corners around that vertex (two pairs of non-necessarily consecutive integers $a, b$ and $c, d$ in the cycle - $a,d$ on one hand and $b, c$ on the other are not necessarily distinct). The operation which consists in duplicating this vertex in two as shown in Fig.~\ref{fig:Vert-split} is called \emph{splitting the vertex (along these two corners)}. It corresponds to composing from the left the corresponding permutation (e.g.~$\sigma$) by the transposition $(d,b)$ (the new map is  $(\sigma\circ (d,b), \tau)$).  A \emph{cut-vertex} can be split in a way that raises the number of connected components of the map.

 A map for which there are no cut-vertices is usually called non-separable or 2-connected \cite{Tutte, Brown, Brown-Tutte, CombEnum, Jacq-Schaeff, DLDRP}. A bipartite map which has no white cut-vertices will be called a  \emph{non-separable hypermap}.\footnote{Hypermaps are a generalization of maps where the vertices are linked by hyper-edges that may connect any number of vertices. They are nothing more than bipartite maps \cite{Walsh}, but still a different role is played by the vertices (here the white vertices) and the hyper-edges (here the black vertices). This justifies calling non-separable hypermap  a bipartite map which is non-separable on the white vertices.}

\begin{figure}[h!]
\centering
\includegraphics[scale=1.1]{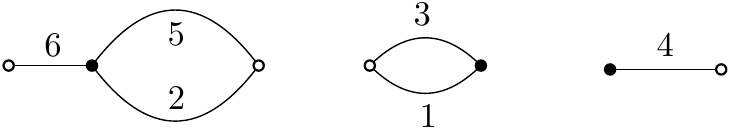}
\caption{Decomposition of the bipartite map of Fig.~\ref{fig:ex-bip-map} in non-separable hypermaps.}
\label{fig:bip-map-NS-dec}
\end{figure}

It is known that any map has a unique decomposition in non-separable components (see e.g.~\cite{CombEnum}). This is also true for a non-separable hypermap, with a simple modification of the procedure: recursively split the white vertices along two corners whenever it raises the number of connected components until there are no such pairs of corners. The operation is commutative and the resulting terminal bipartite map uniquely defined. The connected components of the terminal bipartite map provide a unique decomposition of the original map in non-separable connected hypermaps. See the example in Fig.~\ref{fig:bip-map-NS-dec}.

We stress that by Tutte duality (but linking the newly added grey vertices to the black vertices), non-separable hypermaps are dual to \emph{simple maps} \cite{Relating-ordinary-fully, Simple-maps-topoRec, Analytic-higher}, yet not fully simple, and  such that all faces are boundary faces. 
Splitting a white vertex then amounts to splitting a grey vertex, so that the decomposition in non-separable hypermaps corresponds to a decomposition is simple maps\footnote{It is different from the operation used in \cite{Relating-ordinary-fully} leading to fully simple maps, which duplicates the vertex but keeping the number of connected components fixed.}.

\subsection{Constellations and Hurwitz numbers}

We let for $\bar \pi\in \mathcal P (n)$ and  $\nu\in S_n$ such that $\bar \pi \ge \Pi(\nu)$, and $l, k\in \mathbb{N}$:
\be 
\label{eq:defM}
M(\bar \pi, \nu; l, k) =  \mathrm{Card}\Bigl\{\rho_1,\ldots , \rho_{k}\,\in\, S_{n}^*\ \mid\
  \Pi(\hat \rho) = \bar \pi\ ;\,  \ \rho_1\cdots\; \rho_{k}=\nu\ ;\,   \ \sum_{i=1}^{k} \lVert\rho_{i}\rVert =l \Bigr\}\;,
\ee
where $\Pi(\hat \rho) = \bigvee_{i=1}^k \Pi(\rho_i)$. The systems enumerated by $M$ are called $k$-constellations, and correspond to a specific kind of bipartite maps of genus $g$.  The relation between $l$ and $g$ is then given by: 
\be 
\label{eq:RH}
l= 2g + n-2 + \#(\nu).
\ee

If $\alpha \vdash n$ is a partition of $n$, we let $C_\alpha$ be the conjugacy class  associated to $\alpha \vdash n$, whose cardinal is where $\lvert C_\alpha \rvert = n! / ( \prod_{p\ge 1}  p^{d_p(\alpha)} d_p(\alpha) ! )  $.
For $\alpha \vdash n$ and $g\ge 0$, the genus-$g$ monotone Hurwitz number $\vec H_g(\alpha)$ is defined combinatorially by 
$$
\vec H_g(\alpha) = \lvert C_\alpha \rvert \times \mathrm{Card}\biggl\{ \tau_1, \ldots, \tau_l \mid \rho_\alpha = \tau_1 \cdots \tau_l \ \   ;\ \bigvee_{1\le i \le l}\Pi(\tau_i) = 1_n \biggr\}, 
$$
where $\rho_\alpha $ is any permutation in $C_\alpha$, and the $\tau_i$ are transpositions with weakly monotone maxima, that is, they can be written as $\tau_i=(p_i q_i)$ where $\forall i, \ p_i<q_i$ and $q_1\le \ldots \le q_l$.

Divided by $n!$, $\vec H_g(\alpha)$ counts connected $n$-sheeted weighted  ramified coverings of the 2-sphere by a surface of genus $g$, with $l+1$ ordered branch points, $l$ of which have simple ramification (they have $n-1$ preimages), the last point having ramification profile $\alpha$, where $l$ is given by \eqref{eq:RH}. This relation is then called the Riemann-Hurwitz formula. The condition that the transpositions have weakly monotone maxima restricts the admissible coverings.

If $l\ge 0$ and $\nu \in S_n$,  we define $ \gamma_l(\nu)$ as 
\be
 \gamma_l(\nu)  = \sum_{k\ge 0} (-1)^k M(1_n, \nu; l, k). 
\ee
If $\nu\in C_\alpha$,   $ \gamma_l(\nu)$ can be related to $\vec H_g(\alpha)$  by  \cite{Col-Matsu, Matsu-Nov, CGL}:
\be 
\label{eq:Simple-Hurw-gamma}
\gamma_{\#(\nu) + n + 2g - 2}(\nu) = \vec H_g(\alpha)\; \frac{(-1)^{\#(\nu) + n}}{ \lvert C_\alpha \rvert }\; \; .
\ee 
 For $g=0$, 
from \eqref{eq:RH} $l= n-2 + \#(\nu)$, and $\gamma(\nu) := \gamma_{\#(\nu) + n - 2}(\nu)$, which is related to the genus~0 monotone Hurwitz number by
\be
\label{eq:gamma-hurwitz}
\gamma(\nu) = \vec H_0(\alpha)\; \frac{(-1)^{\#(\nu) + n}}{ \lvert C_\alpha \rvert }\; \; .
\ee
This  is computed explicitly as \cite{Collins03}:
\be
\label{eq:explicit-gamma}
\gamma(\nu) = (-1)^{\#(\nu) + n}  \frac{(2n + \#(\nu) - 3)!}{(2n)!} \prod_{p\ge 1}\left[ \frac{(2p)!}{p! (p-1)!} \right] ^{d_p(\nu)},
\ee
where $d_p(\nu)$ is the number of cycles of $\nu$ with $p$ elements.

\

More generally, we define for $l\ge 0$ and $\bar \pi\in\mathcal{P}(n)$ satisfying $\Pi(\nu)\le \bar \pi < 1_n$:
\be
\label{eq:def-gamma-l-pi}
 \gamma_l(\bar \pi, \nu)  = \sum_{k\ge 0} (-1)^k M(\bar \pi, \nu; l, k).
\ee

\begin{lemma}[\cite{Collins03}, \cite{CGL} Sec.~4.2]
\label{lem:summing-M}
 It can be shown that:
\be
\label{eq:summing-M}
 \gamma_l(\bar \pi, \nu) = \sum_{ \substack{ \{ l_{G} \}_{G \in \bar \pi}\ge 0     \\ \sum_{G} l_{G}  = l } }  \;   
  \prod_{G\in \bar \pi} \gamma_{l_{G}}({\nu}_{\lvert _G} ) . 
\ee
\end{lemma}
Note that this implies that $(-1)^{\#(\nu) + n} \gamma_l(\bar \pi, \nu) \ge 0$. For each $G$, \eqref{eq:RH} writes
\be
\label{eq:lG}
   l_{G} =2g_{G} +  |G|  -2 +  \#({\nu}_{\lvert _G}) , 
\ee
where $g_{G} $ is the genus of the constellations in the alternating sum, or of the coverings counted by the corresponding Hurwitz number.
The proof is not succinct and we don't repeat it here.

\begin{lemma}[\cite{CGL}, Sec.~4.2] 
\label{lem:L}
Let $\bar \pi, \tilde \pi\in \mathcal P (n)$ and $\nu\in S_n$. It is always true that if $\bar \pi\ge \Pi(\nu)$ and $\tilde \pi\ge  \Pi(\nu)$:
\be 
\label{eq:defL}
L[\tilde \pi, \bar \pi; \Pi(\nu)] = \#(\nu) - \#(\tilde \pi)-\#(\bar \pi)+ \#(\tilde \pi\vee \bar \pi)\ge 0.
\ee
Furthermore, if $\pi, \tilde \pi$ such that  $\pi \ge \tilde \pi\ge \Pi(\nu)$, there exists a $\bar \pi\ge  \Pi(\nu)$ such that $\bar\pi\vee \tilde \pi=\pi$ and $L[\tilde \pi, \bar \pi; \Pi(\nu)]=0$.
\end{lemma}
\proof We construct an incidence graph $\mathcal{G}[\tilde \pi, \bar \pi; \Pi(\nu)]$ as follows (see Sec.~\ref{sub:graphs}). To each block of $\tilde \pi$ (resp.~$\bar \pi$ ) we associate a black (resp.~white) vertex. Each block of $\Pi(\nu)$ is contained both in a block of $\tilde \pi$ (corresponding to a black vertex) and a block of $\bar \pi$ (corresponding to a white vertex): for every block of $\Pi(\nu)$, we draw an edge between the corresponding black and white vertices (so that the graph is bipartite).
The graph $\mathcal{G}[\tilde \pi, \bar \pi; \Pi(\nu)]$ has $\#(\tilde \pi\vee \bar \pi)$ connected components, so that $L[\tilde \pi, \bar \pi; \Pi(\nu)]$ is its first Betti number (excess). The inequality is saturated if $\mathcal{G}[\tilde \pi, \bar \pi; \Pi(\nu)]$ is a forest. 
Showing the second assertion is equivalent to showing that given a collection of black vertices, each with a certain number of incident edges such that the other extremity of each edge is not connected to any vertex, it is always possible to group these free extremities to form new white vertices, so that the resulting graph is a tree. This of course is  always true. \qed 

\

The following is a slight generalization of the results of \cite{Collins03} and \cite{CGL}, Sec.~4.2.

\begin{lemma}
\label{lem:min-l}
We fix $\pi, \tilde \pi$ and $\nu$ such that $\pi\ge \tilde \pi\ge \Pi(\nu)$. 
The minimal value of $l$ such that \begin{enumerate}[label=$-$]
\item $\gamma_l(\bar \pi, \nu)$ is non-vanishing 
\item $ \bar \pi$  satisfies $\pi \ge  \bar \pi\ge \Pi(\nu)$ and $\tilde \pi \vee \bar \pi =\pi$, 
\end{enumerate}
is given by 
\be
\label{eq:min-bound-l}
\ell(\nu, \pi, \tilde\pi) =n - \#( \nu) + 2 \bigl(\#( \tilde \pi )- \#(\pi )\bigr).
\ee
For this value to be reached, $\bar \pi$ must be such that 
$$
\#(\bar \pi) = \#(\nu) - \#(\tilde \pi) + \#(\pi), 
$$ 
which is always possible from the second point of Lemma~\ref{lem:L}.
Furthermore,  
\be
\label{eq:summing-M-min}
\gamma_{\ell(\nu, \pi, \tilde\pi)}(\bar \pi, \nu) = 
  \prod_{G\in \bar \pi} \gamma({\nu}_{\lvert _G} ) . 
\ee 
\end{lemma}

\proof Summing the relations \eqref{eq:lG} for the blocks of $\bar \pi$ in Lemma~\ref{lem:summing-M}, we get:
\be
\label{eq:summ-lG}
   l = 2\sum_{G\in \bar \pi} g_G + n - 2 \#(\bar \pi)  + \#(\nu). 
\ee
If $\nu$ is fixed, the minimal value of $l$ is reached for $\bar \pi=0_n$, but this is excluded if the condition $\tilde \pi \vee \bar \pi =\pi$ is imposed for  $\tilde \pi < \pi$. Instead, we may rewrite \eqref{eq:summ-lG}  as: 
\be
l =\#(\nu) + 2 \bigl(\#( \tilde \pi ) - \#(\pi )\bigr) + 2\sum_{G\in \bar \pi} g_G + 2L[\tilde \pi, \bar \pi; \Pi(\nu)],
\ee
so that $l\ge  \ell(\nu, \pi, \tilde\pi)$, with equality if and only if $L[\tilde \pi, \bar \pi, \Pi(\nu)]=0$ and for all $G\in \bar \pi$, $g_G=0$. 
 The second condition fixes all $l_G$ \eqref{eq:lG} to their minimal value $ |G|  -2 +  \#({\nu}_{\lvert _G})$, so that \eqref{eq:summing-M} simplifies to  \eqref{eq:summing-M-min}. \qed

\subsection{Higher order free cumulants}

\subsubsection{Finite $N$}

Higher order free cumulants were introduced in \cite{CMSS}.  As in Notation 4.1 p.~20 of this reference, we define for $\pi\ge \Pi(\tau)$, $\pi\in \mathcal{P}(n)$, $\tau\in S_n$, and $\vec B = (B_1, \ldots, B_n)$: 
\be
\label{eq:Phi-higher-order-0}
\varphi (\pi, \tau) [\vec B]=  \prod_{ G\in \pi} k_{\#(\tau \lvert_{_G})}\biggl(\Bigl\{\Tr\Bigl(\vec \prod _{j\in \tau_i}B_j\Bigr)\Bigr\}_{\tau_i \textrm{ cycle of }\tau \lvert_{_G}}\biggr),
\ee
where the classical cumulants $k_n$ have been introduced in Sec.~\ref{sub:classical cumulants}. For $\vec B = (B, \ldots, B)$ (in which case we remove the arrow), this relation simplifies to:
\be
\label{eq:Phi-higher-order}
\varphi (\pi, \tau) [B]=  \prod_{ G\in \pi} k_{\#(\tau \lvert_{_G})}\biggl(\Bigl\{\Tr\bigl(B^{\lvert V \rvert }\bigr)\Bigr\}_{V\in \Pi(\tau_{\lvert_{_G}})}\biggr).
\ee

\

The following are then defined  for $\sigma\in S_n$ (eq (24) and (14) in \cite{CMSS}): 
\be
\label{eq:Cum-higher-order-fin-N-gen}
\kappa (\pi, \sigma) [\vec B]=  \sum_{\tau\in S_n} \sum_{\substack{{\pi'\in \mathcal{P}(n)}\\{\pi\ge \pi' \ge \Pi(\tau)}}} \varphi(\pi', \tau) [\vec B] \sum_{\substack{{\pi''\in \mathcal{P}(n)}\\{\pi\ge \pi''  \ge \Pi(\sigma)\vee\pi'}}} \lambda_{\pi'', \pi} \prod_{G\in \pi'' } W^{(N)}\Bigl(\left(\sigma\tau^{-1}\right)_{\lvert_{_G}}\Bigr),
\ee
where the $W^{(N)}$ are the Weingarten functions, which can be shown (see \cite{ColSni}) to have the  expansion:
\begin{equation}
\label{eq:expansion-wg-1}
W^{(N)}(\nu) = N^{-n}\sum_{k\ge 0}\;\sum_{\substack{{\rho_1,\ldots , \rho_k\,\in\, S_n^*,}\\{ 
 \rho_1\cdots \rho_k =\nu}}}\;(-1)^k\,N^{-\sum_{i=1}^k\lVert\rho_i\rVert } \; ,
\end{equation}
which can be written using \eqref{eq:def-gamma-l-pi} as:
\begin{equation}
\label{eq:expansion-wg-2}
W^{(N)}(\nu) = \sum_{\substack{\pi\in\mathcal{P}(n) \\ \pi \ge \Pi(\nu)}}\sum_{l\ge 0} \, N^{-n - l }\; \gamma_l(\pi, \nu).
\end{equation}
Note that \eqref{eq:Cum-higher-order-fin-N-gen} is well-defined, because since $\pi' \ge \Pi(\tau) = \Pi(\tau^{-1})$, from \eqref{eq:id-pis}: $ \Pi(\sigma)\vee\pi' = \Pi(\sigma)\vee\Pi(\tau^{-1})\vee\pi' = \Pi(\sigma)\vee\Pi(\sigma \tau^{-1})\vee\pi' \ge \Pi(\sigma \tau^{-1})$, so that $\pi'' \ge \Pi(\sigma \tau^{-1})$.

The inverse relation of \eqref{eq:Cum-higher-order-fin-N-gen} is given in \cite{CMSS}:
\be
\label{eq:Mom-higher-order-fin-N-gen}
\varphi (\pi, \sigma) [\vec B]=  \sum_{\substack{{\tau\in S_n, \,\pi'\in \mathcal{P}(n)}\\[+0.2ex]{\pi\ge \pi' \ge \Pi(\tau)}\\{\pi'\vee \Pi(\sigma\tau^{-1}) = \pi}}} \kappa(\pi', \tau) [\vec B] \, \cdot \, N^{\#(\sigma\tau^{-1})}.
\ee

\subsubsection{Asymptotics}
From \cite{CMSS}, $\varphi (\pi, \sigma) [\vec B]$ is of order $N^{2\#(\pi) - \#(\sigma) }$, so that we define
 \be
\label{eq:rescaled-higher-moments}
 \varphi( \pi, \sigma) [\vec b] = \lim_{N\rightarrow \infty } N^{\#(\sigma) - 2\#(\pi ) }  \varphi(\pi, \sigma) [\vec B], 
\ee
 and $\kappa (\pi, \sigma) [\vec B]$ is of order $N^{2\#(\pi )-n -\#(\sigma)}$. Defining the \emph{higher order free cumulants} by
\be
\label{eq:rescaled-higher-cumulants}
 \kappa( \pi, \sigma) [\vec b] = \lim_{N\rightarrow \infty } N^{\#(\sigma) + n- 2\#(\pi ) }  \kappa(\pi, \sigma) [\vec B], 
\ee
the asymptotics of \eqref{eq:Mom-higher-order-fin-N-gen} are derived in \cite{CMSS}. We describe them below.
\begin{proposition}[\cite{CMSS}] 
\label{}For $\pi\in \mathcal{P}(n)$ and $\sigma\in S_n$ satisfying $\pi \ge \Pi(\sigma)$,
\be
\label{eq:Mom-higher-order-asympt}
\varphi (\pi, \sigma) [\vec b]=  \sum_{\substack{{\tau\in S_n, \,\pi'\in \mathcal{P}(n)}\\[+0.2ex]{\pi\ge \pi' \ge \Pi(\tau)}\\{\pi'\vee \Pi(\sigma\tau^{-1}) = \pi}\\{\#(\tau) - 2\#(\pi') - \#(\sigma\tau^{-1}) + n = \#(\sigma) - 2\#(\pi)  }}} \kappa(\pi', \tau) [\vec b].
\ee
\end{proposition}
This is rewritten through the use of a dot-product (see Def.~4.9 and Notation 4.6 of \cite{CMSS})
\be
\label{eq:dot-product}
(\pi_1, \sigma_1)\cdot(\pi_2, \sigma_2) = (\pi, \sigma) \Leftrightarrow \left \{\begin{array}{lll}\pi = \pi_1\vee \pi_2\\ \sigma = \sigma_1\sigma_2 \\ \#(\sigma_1) - 2\#(\pi_1) + \#(\sigma_2) - 2\#(\pi_2) + n = \#(\sigma) - 2\#(\pi)\end{array}\right. , 
\ee
so that:
\be
\label{eq:Mom-higher-order-asympt-dot}
\varphi (\pi, \sigma) [\vec b]=  \sum_{\substack{{\tau\in S_n, \,\pi'\in \mathcal{P}(n)}\\[+0.2ex]{\pi' \ge \Pi(\tau)}\\[+0.2ex]{(\pi', \tau)\cdot(\Pi(\sigma\tau^{-1}), \;\sigma\tau^{-1}) = (\pi, \sigma)}}} \kappa(\pi', \tau) [\vec b].
\ee

\ 

 A function $f : \pi \ge \sigma \rightarrow f(\pi, \sigma)$  is said to be \emph{multiplicative} if $f(\pi, \sigma) = \prod_{G\in \pi} f(1_{{\lvert G\rvert}}, \sigma_{\mid{_{\lvert G\rvert}}})$, and $f(1_n, \sigma)$ depends only on $\Lambda(\sigma)$.  
 in the sense that for all $\nu\in S_n$ and denoting by $\vec b_\nu=(b_{\nu(1)}, \ldots, b_{\nu(n)})$ whereas $\vec b=(b_{1}, \ldots, b_{n})$:
 \be
 \varphi (1_n, \nu^{-1} \sigma \nu) [\vec b_\nu]=   f(1_n, \sigma ) [\vec b],
 \ee
and similarly for $\kappa$.  Both $\varphi (\pi, \sigma) [\vec b]$ and $\kappa(\pi, \sigma) [\vec b]$ are multiplicative functions. Given $\pi, \tau$,  $\pi \ge \Pi(\tau)$, and two multiplicative functions $f$ and $g$,  their \emph{convolution} is defined as:
 \be
 \label{eq:def-convolution}
 (f\star g)(\pi, \sigma) = \sum_{\substack{{\sigma_1, \sigma_2\in S_n,\ \pi_1, \pi_2 \in \mathcal{P}(n)}\\[+0.2ex]{\pi_1 \ge \Pi(\sigma_1),\  \pi_2 \ge \Pi(\sigma_2)}\\[+0.2ex]{(\pi_1, \sigma_1)\cdot(\pi_2, \sigma_2)=(\pi, \sigma)}}}  f(\pi_1, \sigma_1) g(\pi_2, \sigma_2). 
 \ee
The relations \eqref{eq:Mom-higher-order-asympt} and  \eqref{eq:Mom-higher-order-asympt-dot} can in turn be expressed as
\be
\label{eq:phi-in-conv}
\varphi (\pi, \sigma) [\vec b] =  (f\star \zeta)(\pi, \sigma), 
\ee
where $f : \pi \ge \sigma \rightarrow \kappa(\pi, \sigma)[\vec b]$ and $\zeta: \pi \ge \sigma \rightarrow \delta_{\pi, \Pi(\sigma)}$.  From the multiplicativity and the rules of this convolution, \eqref{eq:phi-in-conv} and thereby \eqref{eq:Mom-higher-order-asympt} can be inverted, yielding (see Def 7.4 of \cite{CMSS}):
\begin{theorem}[\cite{CMSS}] For $\pi\in \mathcal{P}(n)$ and $\sigma\in S_n$ satisfying $\pi \ge \Pi(\sigma)$,
\be
\label{eq:cum-mom-CMSS-prod}
\kappa (\pi, \sigma) [\vec b]=\sum_{\substack{{\sigma_1, \sigma_2\in S_n,\ \pi_1, \pi_2 \in \mathcal{P}(n)}\\[+0.2ex]{\pi_1 \ge \Pi(\sigma_1),\  \pi_2 \ge \Pi(\sigma_2)}\\[+0.2ex]{(\pi_1, \sigma_1)\cdot(\pi_2, \sigma_2)=(\pi, \sigma)}}} \varphi(\pi_1, \sigma_1)[\vec b] \, \,\mu(\pi_2, \sigma_2), 
\ee
where $\mu$ is given as a series by 
\be
\label{eq:moebius-convolution}
\mu(\pi_2, \sigma_2) = \delta(\pi_2, 0_n) \delta(\sigma_2, \mathrm{id}) + \sum_{k\ge 0} \, (-1)^k\sum_{\substack{{\rho_1, \ldots, \rho_k \in S_n^*}\\{(\pi_2, \sigma_2) = (\Pi(\rho_1), \rho_1)\cdot \ldots \cdot(\Pi(\rho_k), \rho_k) }}} 1,
\ee
where the $\delta$ function are Kronecker functions.
\end{theorem}

More explicitly:
\be
\label{eq:cum-mom-CMSS-expl}
\kappa (\pi, \sigma) [\vec b]=\sum_{\substack{{\sigma_1, \sigma_2\in S_n,\ \pi_1, \pi_2 \in \mathcal{P}(n)}\\[+0.1ex]{\pi_1 \ge \Pi(\sigma_1),\  \pi_2 \ge \Pi(\sigma_2)}\\[+0.5ex]{\pi_1\vee \pi_2 = \pi, \ \sigma_1\sigma_2=\sigma}\\[+0.4ex]{\#(\sigma_1) - 2\#(\pi_1 )  + \#(\sigma_2)-  2\#(\pi_2 ) + n = \#(\sigma) -  2\#(\pi)}}} \varphi(\pi_1, \sigma_1)[\vec b] \ \, \,\mu(\pi_2, \sigma_2). 
\ee

Whereas this result is obtained in \cite{CMSS} using the properties of the convolution to invert \eqref{eq:phi-in-conv}, we will reformulate  \eqref{eq:cum-mom-CMSS-prod} in Sec.~\ref{sec:mom-cumulants-two-instead-four} as a tree-like sum involving only one partition and one permutation instead of two each, and show that this simpler expression can be derived directly by taking the asymptotics of \eqref{eq:Cum-higher-order-fin-N-gen} using the properties of the Weingarten functions (Sec.~\ref{sub:first-ref-cum-direct}).

\

Assuming that $\vec b = (b, \ldots, b)$, the higher order moments and cumulants only depend on the partition $\lambda=\Lambda(\sigma) = (\lambda_1, \ldots, \lambda_p)\vdash n $ of the permutation $\sigma \in S_n$ (see \eqref{eq:Phi-higher-order}), that is, for any $\sigma\in C_\lambda$, where $\#(\lambda)=p$:
$$\varphi_{p} ( \lambda) [b] := \varphi (1_n, \sigma) [b] = \lim_{N\rightarrow \infty} N^{p - 2}k_{p}\Bigl(\bigl\{\Tr(B^{\lambda_i})\bigr\}_{1 \le i \le p}\Bigr).$$
The following notation is also common:
\be
\varphi_{ p}(b^{\lambda_1}, \ldots, b^{\lambda_{p}}) =  \varphi_{p} ( \lambda) [b].
\ee
In the same way, for any $\sigma\in C_\lambda$, where $\#(\lambda)=p$:
\be
\label{eq:kappa-dep-only-part} 
 \kappa
 (\lambda)[b] := \kappa (1_n, \sigma) [b] = \kappa (1_n,\gamma_{\lambda_1, \ldots, \lambda_{p}}) [b] ,
\ee
where  $\gamma$ was defined in  \eqref{eq:gamma-cycles}. We use the following notation for the higher order free cumulants:
\be
\kappa_{\lambda_1, \ldots, \lambda_p} [b] =  \kappa_{\{\lambda_i\}} [b] =  \kappa
(\lambda)[b].
\ee

The free cumulants (of first order) are obtained for $\lambda$  the one-part partition ($p=1$) and are denoted by 
\be 
 \kappa_n[b] : = \kappa (1_n, \gamma_n) [b].
 \ee

\subsubsection{Resolution} 

The following reformulation of \eqref{eq:Mom-higher-order-asympt-dot} has been derived in \cite{Voicu} for $p=\#(\lambda)=1$, in \cite{CMSS} for $p=2$, and in  \cite{Analytic-higher} for $p\ge3$ (Thm.~3.12). These combinatorial relations correspond - coefficientwise - to the analytic functional relations between the generating series of higher order moments and cumulants (we do not repeat the functional relations here in full generality). We recall that $\mathcal{G}_p$ is the set of  bipartite trees with $p$ labeled white vertices and no black leaves, and define the following generating functions:
\be
\label{eq:def-gen-func-order2}
M_2(X_1, X_2) = \sum_{a_1, a_2\ge 1}  \varphi_2(b^{a_1}, b^{a_2}) X_1^{a_1}X_2^{a_2}, \qquad \mathrm{and} \qquad C_2(Y_1, Y_2) =\sum_{a_1, a_2\ge 1}  \kappa_{a_1, a_2}[b]\  Y_1^{a_1}Y_2^{a_2}.
\ee 
 
\begin{theorem}[\cite{Voicu, CMSS, Analytic-higher}] 
\label{thm:analytic-higher}
For $p\ge 1$ and $\lambda \vdash n$, $n\ge 1$,
\begin{align*}
& \tilde \varphi_p(b^{\lambda_1}, \ldots, b^{\lambda_{p}}) = 
 \sum_{T\in \mathcal{G}_{p}} \sum_{i=1}^{p}  \sum_{l_i\ge  1}
 \sum_{\substack{
 {\{j^I_i\ge 1\}_{I\in \mathcal{I}_i(T)}}
 \\{\{m_i^{(s)}\ge 1\}_{1\le s\le  l_i }}
 \\[+0.5ex]{\sum_{I\in \mathcal I_i(T)} j^I_i + \sum_{s=1}^{l_i} m_i^{(s)} = \lambda_i}
 }}
 \left(\prod'_{I\in \mathcal{I}(T)} \kappa_{\{j_i^I\}_{i\in I}}[b]\right) \\[+0.5ex]&\hspace{7cm}\times\left(\prod_{i=1}^{p}  \frac{\lambda_i !}{l_i!(\lambda_i + 1 - l_i -  \mathrm{deg}_T(i))! } \prod_{s=1}^{l_i} 
\kappa_{m_i^{(s)}}[b]\right) ,
\end{align*}
where  $\mathcal{I}_i(T)=\{I\in \mathcal{I}(T), i\in I\}$, $\mathrm{deg}_T(i)=\lvert \mathcal{I}_i(T) \rvert $ is the valency of the $i$th white vertex of $T$, and $\prod'$ means that if $\lvert I \rvert =2$, $\kappa_{j_1, j_2}[b]$ should be replaced with 
\be
\kappa_{j_1, j_2}[b] \quad \rightarrow \quad \alpha_{j_1, j_2} [b] = \langle Y_1^{j_1}\rangle\langle Y_2^{j_2}\rangle  \Biggl\{ C_2(Y_1, Y_2) +  \frac{Y_1 Y_2}{(Y_1 - Y_2)^2} \Biggr\},
\ee
and $ \tilde \varphi_p$ is $\varphi_p$ for $p\neq 2$, and 
\be 
\tilde \varphi_2(b^{\lambda_1}, b^{\lambda_{2}}) = \langle X_1^{\lambda_1}\rangle\langle X_2^{\lambda_2}\rangle  \Biggl\{ M_2(X_1, X_2) + \frac{X_1 X_2}{(X_1 - X_2)^2} \Biggr\}.
\ee
 \end{theorem}
 
\paragraph{Immediate reformulations.}
Each sum over the variables $m_i^{(s)}$ can be traded for a sum over $k_1^{(i)}, \ldots, k_{\lambda_i}^{(i)}$  for the cost of a multinomial $(\sum_{a} k_a^{(i)})!/\prod_a k_a^{i}!$, with $l_i=\sum_{a=1}^{\lambda_i} k_a^{(i)}$:
\begin{align}
\label{eq:com-form-with-k}
& \varphi_p(b^{\lambda_1}, \ldots, b^{\lambda_{p}}) = 
 \sum_{T\in \mathcal{G}_{p}} \sum_{i=1}^{p}  \sum_{l_i\ge 1}
\sum_{\substack{{\{j^I_i\ge 1\}_{I\in \mathcal{I}_i(T)}}\\
 {k_1^{(i)}, \ldots, k^{(i)}_{\lambda_i}\ge 0, \, \, \sum_{a=1}^{\lambda_i}  k_a^{(i)} =l_i }
 \\[+0.5ex]{\sum_{I\in \mathcal I_i(T)}j^I_i + \sum_{a=1}^{\lambda_i} a k_a^{(i)} = \lambda_i}
 }}
 \left(\prod'_{I\in \mathcal{I}(T)} \kappa_{\{j_i^I\}_{i\in I}}[b]\right) \\&\hspace{7cm}\times\left(\prod_{i=1}^{p}  \frac{\lambda_i !}{ (\lambda_i+1  - l_i - \mathrm{deg}_T(i) )! } \prod_{a=1}^{\lambda_i} \frac 1 {k_a^{(i)}!}\,  
 \kappa_{a}[b]^{k_a^{(i)}} \right)\nonumber.
\end{align}
We may sum over $T'\in \mathcal{T}_p$ instead of $T\in \mathcal{G}_p$, by adding $l_i$ black leaves to the vertex $i$.  For  $T'$, $\mathrm{deg}_T(i)$ is the number of edges not incident to leaves in $T'$,  $l_i = \mathrm{deg}^\mathsf{L}_{T'}(i)$, so that  $l_i +  \mathrm{deg}_T(i) = \mathrm{deg}_{T'}(i) $. 
Some weights have been assigned to the edges of the tree: 
\begin{itemize}
\item the edges incident to black vertices of valency larger than one are identified by a pair $(i, I)$, $1\le i \le p$, $I\in \mathcal{I}'(T)$, and have been assigned a weight $j_i^I$,
\item the black leaves incident to the white vertex $i\in \{1, \ldots, p\}$ are undistinguishable, so to each  white vertex $i$ is associated an unordered list of weights, that is, a  list $\{k_a^{(i)}\}_{1\le a \le \lambda_i}$, where $k_a^{(i)}$ is the number of black leaves associated to the white vertex $i$ that have a weight equal to $a$.
\end{itemize}
The formula of the theorem may therefore be re-stated as follows:
\begin{align} \varphi_p(b^{\lambda_1}, \ldots, b^{\lambda_{p}}) =  \sum_{\substack{{T\in \mathcal{T}_{p}}\\{\forall i, \deg_T^{\mathsf L} (i) \ge 1}}} \prod_{i=1}^{p}  \frac{\lambda_i !}{ (\lambda_i + 1  - \mathrm{deg}_T(i) )! } \sum_{i=1}^{p}   \sum_{\substack{ {\{j^I_i\ge 1\}_{I\in \mathcal{I}_i(T)}} \\[+0.5ex]{\sum_{I\in \mathcal I_i(T)} j^I_i  = \lambda_i} }} \prod'_{I\in \mathcal{I}(T)} \kappa_{\{j_i^I\}_{i\in I}}[b] \prod_{i=1}^p\prod_{a=1}^{\lambda_i} \frac 1 {k_a^{(i)}!},\end{align}
where as discussed above, the weights $j_i^{\{i\}}$ are an unordered list of $\mathrm{deg}_T^{\mathsf L}(i)$ positive integers, and $k_a^{(i)}$ is the number of these weights equal to $a$.


\paragraph{First order.} 
 For $p=1$, we use the formulation \eqref{eq:com-form-with-k} of the theorem. $\lambda$ has one block $\lambda_1=n$. The convention is that a tree with one white vertex but no leaves is just a single vertex (so that $\mathcal{I}(T)$ is empty). Thus,  \eqref{eq:com-form-with-k} simplifies to:
\be \varphi_1(b^{n}) = \sum_{l = 1}^n \sum_{\substack{{\{k_i\ge 0\}_{1\le i \le n}}\\[+0.2ex]{\sum_{i=1}^{n} k_i =l  }\\[+0.2ex]{\sum_{i=1}^{n} i k_i =n  }}} \frac{n!}{(n+1-l)!} \prod_{i=1}^n \frac{\kappa_i[b]^{k_i}}{k_i!} = \sum_{\substack{{\{k_i\ge 0\}_{1\le i \le n}}\\[+0.2ex]{\sum_{i=1}^{n} i k_i =n  }}} \frac{n!}{\prod_{i} k_i ! (n+1-\sum_{i=1}^n k_i)!} \prod_{i=1}^n \kappa_i[b]^{k_i}, \ee 
which as detailed in Corollary 2.6  of  \cite{Voicu} or  Corollary 1 of \cite{Spei94} is a combinatorial reformulation of  the relation  between generating functions 
\be
\label{eq:first-order-funct}
C_1(XM_1(X)) = M_1(X),
\ee
where
\be
M_1(X) = 1 +  \sum_{a\ge 1}  \varphi_1(b^a) X^a, \qquad \mathrm{and} \qquad C_1(Y) = 1 +  \sum_{a\ge 1}  \kappa_a [b] Y^a.
\ee 

\

\paragraph{Second order.}  For $p=2$, there is a single tree in $\mathcal{G}_2$, with two white vertices and one black vertex  of valency two ($\mathcal I$ consists of a single set $\{1,2\}$). Theorem~\ref{thm:analytic-higher} simplifies to: 
\begin{align}
 \tilde \varphi_2(b^{\lambda_1}, b^{\lambda_{2}}) = 
 \sum_{\substack{{j_1, j_2\ge 1}}}
 \alpha_{j_1, j_2}[b]\ 
 \prod_{i=1}^{2} \Biggl(\sum_{l_i\ge 1} \binom{\lambda_i}{l_i} \sum_{\substack{{\{m_i^{(s)}\ge 1\}_{1\le s\le  l_i }}
 \\[+0.5ex]{ \sum_{s=1}^{l_i} m_i^{(s)} = \lambda_i - j_i}
 }}
\prod_{s=1}^{l_i} \kappa_{m_i^{(s)}}[b]\Biggr).
\end{align}
To carry out the sums on the right hand side, we notice that since $\lambda>j$ (as $s\ge 1$):
\be
\label{eq:sum-the-l}
\sum_{l\ge 1} \binom{\lambda}{l} \sum_{\substack{{\{m^{(s)}\ge 1\}_{1\le s\le  l }}
 \\[+0.5ex]{ \sum_{s=1}^{l} m^{(s)} = \lambda - j}
 }}
\prod_{s=1}^{l} \kappa_{m^{(s)}}[b]
 =\sum_{l= 1}^{\lambda - j} \binom{\lambda}{l} \langle z^{\lambda - j}\rangle\biggl(\sum_{a=1}^{\lambda - j}\kappa_a[b]z^a\biggr)^{l}
 = \langle z^{\lambda - j}\rangle \biggl(1+\sum_{a=1}^{\lambda - j}\kappa_a[b]z^a\biggr)^{\lambda} , 
\ee
where $\langle z^n \rangle (\sum_{a\ge 0} c_a z^a) = c_n$.
By Lagrange inversion (see \cite{LagrangeGessel}, Eq.~(2.2.1)):
\be
\label{eq:Lagrange}
\langle z^{\lambda - j}\rangle \biggl(1+\sum_{a\ge1}\kappa_a[b]z^a\biggr)^{\lambda-K} =  \frac {\lambda - K} {j-K}\  \langle t^{\lambda - K} \rangle \biggl(t + \sum_{a\ge 1}\varphi_1(b^a)t^{a+1} \biggr)^{j-K},
\ee 
so that (here $K=0$):
\begin{align}
\tilde  \varphi_2(b^{\lambda_1}, b^{\lambda_{2}}) = 
 \sum_{\substack{{j_1, j_2\ge 1}}}
 \alpha_{j_1, j_2}[b]\ 
 \prod_{i=1}^{2} \Biggl(\frac {\lambda_i} {j_i} \sum_{\substack{{n_i^{(1)}, \ldots, n_{i}^{(j_i)} \ge 0}\\{\sum_{s=1}^{j_i} n_i^{(s)} =\lambda_i - j_i}}} \prod_{s=1}^{j_i} \varphi_1(b^{n_i^{(s)}})\Biggr), 
\end{align}
with the convention that $\varphi_1(b^{0})=1$. Noticing that
\be
\label{eq:getting-ni}
j \sum_{\substack{{\{n^{(s)}\ge 0\}_{1\le s\le  j }}
 \\[+0.5ex]{ \sum_{s=1}^{j} n^{(s)} = P}
 }}
\hspace{-0.3cm}n^{(1)}\prod_{s=1}^{j} \varphi_1\Bigl(b^{n^{(s)}}\Bigr) 
= \sum_{s=1}^j \sum_{\substack{{\{n^{(s)}\ge 0\}_{1\le s\le  j }}
 \\[+0.5ex]{ \sum_{s=1}^{j} n^{(s)} = P}
 }}
\hspace{-0.3cm}n^{(s)}\prod_{s=1}^{j} \varphi_1\Bigl(b^{n^{(s)}}\Bigr) 
=
P \sum_{\substack{{\{n^{(s)}\ge 0\}_{1\le s\le  j }}
 \\[+0.5ex]{ \sum_{s=1}^{j} n^{(s)} = \lambda - j}
 }}
\prod_{s=1}^{j} \varphi_1(b^{n^{(s)}}), \ee
where $P=\lambda - j$, we obtain the following:
\begin{align}
\label{eq:comb-form-p2}
 \tilde \varphi_2(b^{\lambda_1}, b^{\lambda_{2}}) = \hspace{-0.1cm}
 \sum_{\substack{{j_1, j_2\ge 1}}}\hspace{-0.1cm}
 \alpha_{j_1, j_2}[b]\ 
 \prod_{i=1}^{2} \Biggl(\sum_{\substack{{\{n_i^{(s)}\ge 0\}_{1\le s\le  j_i }}
 \\[+0.5ex]{ \sum_{s=1}^{j_i} n_i^{(s)} = \lambda_i - j_i}
 }}
\hspace{-0.3cm}(1+n_i^{(1)})\prod_{s=1}^{j_i} \varphi_1\Bigl(b^{n_i^{(s)}}\Bigr)\Biggr),
\end{align}
which is to be compared with the key combinatorial reformulation of \eqref{eq:Mom-higher-order-asympt} derived in  Prop.~6.2 of \cite{CMSS}. Going to generating functions, we indeed recover from \eqref{eq:comb-form-p2} the  second order functional relations
\be
\label{eq:second-order-functional-rel}
M_2(X_1, X_2)  + \frac{X_1 X_2}{(X_1 - X_2)^2} = \frac{\mathrm{d} \ln Y_1}{\mathrm{d} \ln X_1} \frac{\mathrm{d} \ln Y_2}{\mathrm{d} \ln X_2}\left ( C_2(Y_1, Y_2)  + \frac{Y_1 Y_2}{(Y_1 - Y_2)^2}\right),
\ee
where $M_1$ has been defined in \eqref{eq:first-order-funct}, $M_2$ and $C_2$ in \eqref{eq:def-gen-func-order2}, and 
\be
 Y_i = X_i M_1(X_i), \qquad \mathrm{and} \qquad \frac{\mathrm{d} \ln Y_i}{\mathrm{d} \ln X_i} = \frac {\frac {\mathrm{d}} {\mathrm{d}X_i} Y_i}{M_1(X_i)}.
 \ee

\

For $n\ge 3$, the formulas are derived 
through Fock space methods, and deriving a combinatorial proof is still open (see the discussion page 24 of  \cite{Analytic-higher}).

\subsubsection{Simpler formulation of the generating functions formulas}
\label{subsub:simplificaion-funct-relations-gen}

In view of the computations carried above for $p=2$, we may reformulate Thm.~\ref{thm:analytic-higher} as follows.

\begin{align}
\label{eq:comb-good-for-proof}
& \varphi_p(b^{\lambda_1}, \ldots, b^{\lambda_{p}}) = 
 \sum_{T\in \mathcal{G}_{p}} 
 \ \sum_ {\{j^I_i\ge 1\}_{\substack{{1\le i \le p}\\{I\in \mathcal{I}_i(T)}}}}
 \left(\prod'_{I\in \mathcal{I}(T)} \kappa_{\{j_i^I\}_{i\in I}}[b]\right) \\&\hspace{4cm}\times \prod_{i=1}^{p}  \frac{\lambda_i !}{(\lambda_i+1  -   \mathrm{deg}_T(i))! }  
\langle z^{\lambda_i - \sum_{I\in \mathcal{I}_i(T)} j_i^I}\rangle \biggl(1 + \sum_{a\ge 1} \kappa_a[b] z^a\biggr) ^{\lambda_i + 1 - \mathrm{deg}_T(i) }.
\nonumber
\end{align}
By Lagrange inversion again (\eqref{eq:Lagrange} for $K=\mathrm{deg}_T(i) - 1$), and then using \eqref{eq:getting-ni} with $J_i - \mathrm{deg}_T(i) + 1$ instead of $j$, where $J_i = \sum_{I\in \mathcal{I}_i(T)} j_i^I$ and with $P=\lambda_i - J_i$, we obtain the generalization for arbitrary $p\ge 2$ of the expression  analogous to \eqref{eq:comb-form-p2} for $p=2$:
\begin{align}
\label{eq:comb-to-gen}
& \varphi_p(b^{\lambda_1}, \ldots, b^{\lambda_{p}}) = 
 \sum_{T\in \mathcal{G}_{p}} 
 \ \sum_ {\{j^I_i\ge 1\}_{\substack{{1\le i \le p}\\{I\in \mathcal{I}_i(T)}}}}
 \left(\prod'_{I\in \mathcal{I}(T)} \kappa_{\{j_i^I\}_{i\in I}}[b]\right) \\&\hspace{3cm}\times\prod_{i=1}^{p}  \frac{\lambda_i !}{(\lambda_i + 1 -   \mathrm{deg}_T(i))! }  
  \sum_{\substack{{\{n_i^{(s)}\ge0\}_{1\le s\le  J_i  - \mathrm{deg}_T(i) + 1 }}
 \\[+0.5ex]{ \sum_{s=1}^{J_i} n_i^{(s)} = \lambda_i - J_i}
 }}
\hspace{-0.3cm}(1+n_i^{(1)})\prod_{s=1}^{J_i- \mathrm{deg}_T(i) + 1 } \varphi_1\Bigl(b^{n_i^{(s)}}\Bigr)   ,
\nonumber
\end{align}
with the convention that $\varphi_1(b^{0})=1$. Noticing that:
\begin{align}
 \sum_{\substack{{\{n^{(s)}\ge0\}_{1\le s\le  J-K}}
 \\[+0.5ex]{ \sum_{s=1}^{J} n^{(s)} = \lambda - J}
 }}
\hspace{-0.3cm}(1+n^{(1)})\prod_{s=1}^{J-K} \varphi_1\Bigl(b^{n^{(s)}}\Bigr)  &= \langle X^{\lambda-K}\rangle\ X \frac {\mathrm d} {\mathrm d X} \bigl(X M_1(X)\bigr) \bigl(X M_1(X)\bigr)^{J-K -1} \\&= \langle X^{\lambda}\rangle\ \frac{1}{M_1(X)^{K+1}} \frac {\mathrm d} {\mathrm d X} \bigl(X M_1(X)\bigr) \bigl(X M_1(X)\bigr)^{J}, 
\end{align}
and that 
\be 
\label{eq:identity1}
\frac {\lambda ! }{(\lambda - K)!}\langle X^\lambda\rangle f(X)= \langle X^{\lambda-K}\rangle\frac {\mathrm d^{K}} {\mathrm d X^{K}} f(X) = \langle X^{\lambda}\rangle X^K  \frac {\mathrm d^{K}} {\mathrm d X^{K}} f(X),
\ee 
we may write: 
\begin{align}
& \varphi_p(b^{\lambda_1}, \ldots, b^{\lambda_{p}}) = \biggl( \prod_{i=1}^p \langle X_i^{\lambda_i}\rangle  \biggr)
 \sum_{T\in \mathcal{G}_{p}}  \ \sum_ {\{j^I_i\ge 1\}_{\substack{{1\le i \le p}\\{I\in \mathcal{I}_i(T)}}}}
 \left(\prod'_{I\in \mathcal{I}(T)} \kappa_{\{j_i^I\}_{i\in I}}[b]\right) \\&\hspace{4cm}\times
 \left( \prod_{i=1}^p 
 X_i^{\mathrm{deg}_T(i) -1} \frac {\mathrm d^{\mathrm{deg}_T(i) -1}} {\mathrm d X_i^{\mathrm{deg}_T(i) -1}}
  \left\{ \frac{\frac {\mathrm d} {\mathrm d X_i} (X_iM_1(X_i))}{M_1(X_i)^{\mathrm{deg}_T(i)}}
 \bigl(X_iM_1(X_i)\bigr)^{J_i}  \right\} \right) .
\nonumber
\end{align}
Introducing the higher order generating functions 
\be
\begin{split}
&M_p(X_1, \ldots, X_p) = \sum_{a_1, \ldots, a_p\ge 1}  \varphi_n(b^{a_1}, \ldots, b^{a_p}) X_1^{a_1}\ldots X_p^{a_p}, \\ 
&C_p(Y_1, \ldots, Y_p) =\sum_{a_1,\ldots, a_p\ge 1 }  \kappa_{a_1, \ldots, a_p}[b]\  Y_1^{a_1}\ldots Y_p^{a_p}, 
\end{split}
\ee
we get the following equivalent reformulation of Thm.~\ref{thm:analytic-higher}:

\begin{theorem} 
\label{thm:theorem-functional-reformulated}
For $p\ge 3$, the generating functions $M_p$ and $C_p$ satisfy the relation 
\be
M_p(X_1, \ldots, X_p)  = \sum_{T\in \mathcal{G}_p} \biggl( \prod_{i=1}^p 
  X_i^{\mathrm{d}_T(i) -1} \frac {\mathrm d^{\mathrm{d}_T(i) -1}} {\mathrm d X_i^{\mathrm{d}_T(i) -1}}\biggr)\Biggl\{\prod_{i=1}^p
   {\frac {\mathrm{d} Y_i} {\mathrm d X_i}} \prod'_{I\in \mathcal{I}(T)} \ \frac{C_{\lvert I \rvert}\bigl(Y_I\bigr)}{\prod_{i\in I}M_1(X_i)}\Biggr\},   
\ee
where $\mathrm{d}_T(i) = \mathrm{deg}_T(i)$,  $Y_I=\{Y_i\}_{i\in I}$, and
$Y_i = X_i M_1(X_i)$,  
and $\prod'$ means that any occurence of $C_2(Y_1, Y_2)$ should be replaced by
\be
C_2 (Y_1, Y_2) \quad \rightarrow \quad C_2 (Y_1, Y_2) + \frac {Y_1Y_2}{(Y_1 - Y_2)^2}.
\ee
The formula remains true for $p=2$ if $M_2(X_1, X_2)$ is replaced with $M_2(X_1, X_2) + \frac{X_1X_2}{(X_1 - X_2)^2}$, and for $p=1$ if ``${\frac {\mathrm d^{ -1}} {\mathrm d X_i^{ -1}}{\frac {\mathrm{d} Y_i} {\mathrm d X_i}}  = Y_i}$''.
\end{theorem}

This formula can be seen to be equivalent to that derived by Hock in Corollary 5.4 of \cite{Hock} (if in that formula, $\mathrm{Aut}(T)=1$, which is consistent with \cite{Analytic-higher}) by noticing the following: if $f(X)=\sum_{k\ge r} c_k X^k$, then 
\be
\frac {\lambda ! }{(\lambda - K)!}\langle X^\lambda\rangle f(X)= \langle X^{\lambda} \rangle\frac 1 {X}\biggl (X^2\frac {\mathrm d} {\mathrm d X}\biggr)^K \left\{\frac{1}{X^{K-1}}f(X)\right \} , 
\ee
and using this identity in \eqref{eq:comb-good-for-proof} together with  \eqref{eq:comb-to-gen} instead of using the simpler identity \eqref{eq:identity1}.

\newpage

\section{Looking for the trees}
\label{sec:tree-structure}

In this section, we re-express the moment-cumulant formulas using one permutation and one partition, instead of two each, as in \eqref{eq:Mom-higher-order-asympt-dot} and \eqref{eq:cum-mom-CMSS-prod}. This way of rewriting the formulas shows that they involve a sum over bipartite trees as in Thm.~\ref{thm:analytic-higher}. This last point is clarified in Sec.~\ref{sec:bij-with-trees}.

\subsection{Simplification of the moment-cumulant formulas}
\label{sec:mom-cumulants-two-instead-four}

\begin{theorem}
\label{thm:first-tree-moments} For any $\sigma\in S_n$, $\pi\in \mathcal{P}(n)$ satisfying $\pi \ge \Pi(\sigma)$,  the higher order  moment-cumulant formula can be expressed as:
\be
\label{eq:mom-cum-CMSS-facto-rewrite}
\varphi (\pi, \sigma) [\vec b]=
\sum_{\substack{{\tau\in S_n,}\\[+0.3ex]{g(\sigma, \tau^{-1})=0}}} \hspace{0.2cm}
\sum_{\substack{{\pi\ge \pi' \ge \Pi(\tau)}\\[+0.3ex]{\pi'\vee\Pi(\sigma)= \pi}\\[+0.3ex] {L[\pi',\, \Pi(\sigma)\vee \Pi(\tau)\, ; \,\Pi(\tau)] =0
}}}  \kappa(\pi', \tau)[\vec b] ,
\ee
 where $g$ has been defined in \eqref{eq:genus} and we recall the forest-like condition \eqref{eq:defL}: 
 $$
L[\pi', \Pi(\sigma)\vee \Pi(\tau) ; \Pi(\tau) ] =  \#(\tau) - \#(\pi') - \#(\Pi(\sigma)\vee\Pi(\tau)) + \#(\pi'\vee\Pi(\sigma))
.$$
The higher order cumulant-moment relation can be written:
\be
\label{eq:final-higher-order-cum-mom}
\kappa (\pi, \sigma) [\vec b]=   \sum_{\substack{{\tau\in S_n}\\[+0.2ex]{g(\sigma, \tau^{-1}) =0}}} \,\ \sum_{\substack{{\pi'\in \mathcal{P}(n)}\\[+0.2ex]{\pi \ge \pi' \ge \Pi(\tau)}
\\[+0.3ex]{ 
L[\pi', \, \Pi(\sigma)\vee\Pi(\tau)\, ; \, \Pi(\tau)] = 0}}}  \varphi(\pi', \tau) [\vec b] \ \, \Gamma\bigl[\sigma \tau^{-1}, \pi, \pi'\vee \Pi(\sigma)\bigr], 
\ee
where $\Gamma$ is given by a sum of trees-like partitions whose blocks are weighted by rescaled genus 0 monotone Hurwitz numbers $\gamma$ \eqref{eq:gamma-hurwitz}, \eqref{eq:explicit-gamma}: 
\be
\label{eq:comb-expr-leading-cum-Weing}
\Gamma[\nu, \pi, \tilde\pi] = 
\sum_{\substack{{\pi''\ge \Pi(\nu)}\\{\tilde\pi\vee\pi''= \pi }\\[+0.3ex]{ 
L[\pi'', \,\tilde \pi \, ; \, \Pi(\nu)]=0}}} \ 
\prod_{G\in \pi''} \gamma(\nu_{|_G}) \;.
\ee
\end{theorem}

We thus have simpler higher order moment cumulant formulas, in the sense that they only involve one pair $(\pi, \sigma)$ with $\pi\ge \Pi(\sigma)$ satisfying some conditions, instead of two such pairs in \eqref{eq:cum-mom-CMSS-prod}, satisfying the conditions of the dot-product \eqref{eq:dot-product}. The role of the M\"obius inverse is now played by $\Gamma$ instead of $\mu$ \eqref{eq:moebius-convolution}.

\proof  The key point is that the last condition in the sum in \eqref{eq:Mom-higher-order-asympt}
is in fact  equivalent to three conditions (in the flavor of  Thm.~5.6 of \cite{CMSS}, but yet a different rewritting):
\begin{lemma}
\label{lemma:split}
 Let $\pi_1, \pi_2\in \mathcal{P}(n)$ and $\sigma_1, \sigma_2\in S_n$ be such that $\pi_i \ge \Pi(\sigma_i)$. Then:
\begin{align}
\nonumber
&\#(\sigma_1) - 2\#(\pi_1 )  + \#(\sigma_2)-  2\#(\pi_2 ) + n = \#(\sigma_2\sigma_1) -  2\#(\pi_1\vee \pi_2) \\[+1ex] \nonumber
&\hspace{4.5cm}\Leftrightarrow \left\{
    \begin{array}{lll}
      & \#(\sigma_1) + \#(\sigma_2) + \#(\sigma_2\sigma_1) - n =2 \#(\Pi(\sigma_1) \vee \Pi(\sigma_2\sigma_1)) \\
       &\#(\sigma_1) - \#(\pi_1 )  =  \#(\Pi(\sigma_1) \vee \Pi(\sigma_2\sigma_1) ) - \#(\pi_1 \vee \Pi(\sigma_2\sigma_1))\\
	& \#(\sigma_2) - \#(\pi_2 )  =  \#(\pi_1\vee\Pi(\sigma_2\sigma_1)) -\#(\pi_1\vee \pi_2)
	 \\ \nonumber
	     \end{array}
    \right. \\ \label{eq:facto}
&\hspace{7cm}\Leftrightarrow \left\{
    \begin{array}{lll}
      &g(\sigma_1, \sigma_1^{-1}\sigma_2^{-1} )=0 \\
       & L[\pi_1, \Pi(\sigma_1)\vee \Pi(\sigma_2\sigma_1) ; \Pi(\sigma_1)]=0\\
	& L[\pi_2, \pi_1\vee \Pi(\sigma_2\sigma_1) ; \Pi(\sigma_2)]=0
    \end{array}
    \right. .
    \end{align}
    \end{lemma}
    \noindent{\it Proof of the lemma.}  It is true with these assumptions that
    \begin{align}
&    \#(\sigma_1) - 2\#(\pi_1 )  + \#(\sigma_2)-  2\#(\pi_2 ) + n   -  \#(\sigma_2\sigma_1) +  2\#(\pi_1\vee \pi_2) = \\
&\hspace{2.7cm}2 g(\sigma_1, \sigma_1^{-1}\sigma_2^{-1} ) + 2L[\pi_1, \Pi(\sigma_1)\vee \Pi(\sigma_2\sigma_1) ; \Pi(\sigma_1)]  + 2L[\pi_2, \pi_1\vee \Pi(\sigma_2\sigma_1) ; \Pi(\sigma_2)]. \nonumber
    \end{align}
For the last term, one has to justify that $\pi_1\vee \Pi(\sigma_2\sigma_1) \ge \Pi(\sigma_2)$, which is true because $\pi_1\vee \Pi(\sigma_2\sigma_1) \ge\Pi(\sigma_1)\vee \Pi(\sigma_2\sigma_1) = \Pi(\sigma_1)\vee \Pi(\sigma_2)  \ge \Pi(\sigma_2)$, and that $\pi_1\vee \Pi(\sigma_2\sigma_1) \vee \pi_2 = \pi_1\vee \pi_2$, which is true because $\pi_1 \vee \pi_2 \ge \Pi(\sigma_1)\vee \Pi(\sigma_2) \ge \Pi(\sigma_2\sigma_1)$.
Each term on the right hand side is non-negative, so vanishing of the left hand side is equivalent to vanishing of these three quantities. \qed 

\

From \eqref{eq:id-gs}, $g(\sigma_1, \sigma_1^{-1}\sigma_2^{-1} )=0 \Leftrightarrow g(\sigma_2\sigma_1, \sigma_1^{-1})=0$. Let us apply Lemma~\ref{lemma:split} for $\pi_1=\pi'$, $\pi_2=\Pi(\sigma\tau^{-1})$, $\sigma_1=\tau$, $\sigma_2=\sigma\tau^{-1}$. The last of the three conditions in \eqref{eq:facto} then reads
\be
\label{eq:third-condition-int}
\#(\pi'\vee\Pi(\sigma)) = \#(\pi'\vee\Pi(\sigma\tau^{-1})).
\ee
This is always satisfied for $\pi'\ge \Pi(\tau)$, because then (see  \eqref{eq:id-pis}):
$$
\pi'\vee \Pi(\sigma\tau^{-1}) = \pi'\vee \Pi(\tau)\vee\Pi(\sigma\tau^{-1})  = \pi'\vee\Pi(\sigma)\vee\Pi(\tau) = \pi'\vee\Pi(\sigma).
$$
The condition that $\pi'\vee \Pi(\sigma\tau^{-1}) = \pi$ in the sum of \eqref{eq:Mom-higher-order-asympt}) can therefore be replaced with $\pi'\vee \Pi(\sigma) = \pi$. 
The other two conditions are those in \eqref{eq:mom-cum-CMSS-facto-rewrite}, which concludes the proof of \eqref{eq:mom-cum-CMSS-facto-rewrite}.

\

We now turn to the proof of \eqref{eq:final-higher-order-cum-mom}. From \eqref{eq:id-gs}, $g(\sigma_1, \sigma_1^{-1}\sigma_2^{-1} )=0 \Leftrightarrow g(\sigma_2\sigma_1, \sigma_1^{-1})=0$. Using Lemma~\ref{lemma:split}, we can split the sum in \eqref{eq:cum-mom-CMSS-expl} as:
\be
\label{eq:cum-mom-CMSS-facto}
\kappa (\pi, \sigma) [\vec b]=
\sum_{\substack{{\sigma_1, \sigma_2\in S_n,}\\[+0.3ex]{\sigma_1\sigma_2=\sigma}\\{g(\sigma_2\sigma_1, \sigma_1^{-1})=0}}} \hspace{0.2cm}
\sum_{\substack{{\pi\ge \pi_1 \ge \Pi(\sigma_1)}\\[+0.3ex]{
L[\pi_1, \, \Pi(\sigma_1) \vee \Pi(\sigma_2\sigma_1)  \,  ; \, \Pi(\sigma_1)] =0}}}  \varphi(\pi_1, \sigma_1)[\vec b] \, \,A[\sigma_2, \pi, \pi_1],
\ee
where the condition $\pi\ge \pi_1$ is necessary for $\pi_1\vee \pi_2 = \pi$ to be satisfied, and where:
\be
\label{eq:A}
A[\sigma_2, \pi, \pi_1]= \sum_{\substack{{ \pi_2 \ge \Pi(\sigma_2)}\\[+0.3ex]{\pi_1\vee \pi_2 = \pi, }\\[+0.2ex]{ 
L[\pi_2, \, \pi_1 \vee \Pi(\sigma_2\sigma_1)  \,  ; \, \Pi(\sigma_2)] =0}}} 
\mu(\pi_2, \sigma_2). 
\ee
Recalling that $\Gamma$ has been defined in \eqref{eq:comb-expr-leading-cum-Weing}, the rest of this subsection is devoted to showing the following.

\begin{lemma} 
\label{lem:Aisf}
Let  $\pi_1, \pi_2, \pi\in \mathcal{P}(n)$ and $\sigma_1, \sigma_2\in S_n$ be such that $\pi \ge \pi_i \ge \Pi(\sigma_i)$. Then:
\be
A[\sigma_2, \pi, \pi_1] = \Gamma\bigl[\sigma_2, \pi, \pi_1\vee\Pi(\sigma_2\sigma_1)\bigr].
\ee
\end{lemma}

\

Having proved this, changing in \eqref{eq:cum-mom-CMSS-facto} the names of $\pi_1, \sigma_1$ to $\pi', \tau$, suppressing the summation on $\sigma_2$ as $\sigma_2$ is fixed to  $\sigma_2=\sigma\sigma_1^{-1} = \sigma\tau^{-1}$, we recover 
$$
\kappa (\pi, \sigma) [\vec b]=
\sum_{\substack{{\tau\in S_n,}\\[+0.3ex]{g(\sigma, \tau^{-1})=0}}} \hspace{0.2cm}
\sum_{\substack{{\pi\ge \pi' \ge \Pi(\tau)}\\[+0.3ex]{
L[\pi', \, \Pi(\tau) \vee \Pi(\sigma)  \,  ; \, \Pi(\tau)] =0
}}}  \varphi(\pi', \tau)[\vec b] \, \, \Gamma\bigl[ \sigma\tau^{-1}, \pi, \pi'\vee\Pi(\sigma)\bigr],
$$
that is, \eqref{eq:final-higher-order-cum-mom}.

\

\noindent{\it Proof of  Lemma~\ref{lem:Aisf}.} As can easily be seen inductively, the function $\mu(\pi_2, \sigma_2)$ of \eqref{eq:moebius-convolution} can be rewritten as:
\be
\label{eq:rewrite-mu}
\mu(\pi_2, \sigma_2) =\sum_{k\ge 0} \, \, (-1)^k\sum_{\substack{{\rho_1, \ldots, \rho_k \in S_n^*}\\[+0.2ex]{\Pi(\hat \rho)= \pi_2 }\\[+0.6ex]{\rho_1 \cdots \rho_k= \sigma_2}\\[+0.2ex]{\sum_{i=1}^k \#(\rho_i) = 2\#(\pi_2 ) - \#(\sigma_2) + (k-1)n}}} 1,
\ee
where $\Pi(\hat \rho) = \bigvee_{i=1}^k \Pi(\rho_i)$ and with the convention that for $k=0$, $\pi_2$ is the $n$ blocks partition and $\sigma_2$ is the identity. 
The rightmost sum  is the quantity $M$  defined in \eqref{eq:defM} for $l=\#(\sigma_2) - 2\#(\pi_2 ) +n$, so that:
\be
\label{eq:rewrite-mu-M}
\mu(\pi_2, \sigma_2) 
= \gamma_{\#(\sigma_2) - 2\#(\pi_2 ) +n}(\pi_2, \sigma_2),
\ee

 For $\pi_1, \pi_2$ satisfying $\pi \ge \pi_i \ge \Pi(\sigma_2)$ and $\pi_1\vee \pi_2=\pi$, and $\sigma_1$ such that $\pi \ge \Pi(\sigma_1)$, we may again use Lemma~\ref{lem:min-l} for $\pi'':=\pi_2$, $\nu:=\sigma_2$ and $\tilde \pi := \pi_1\vee \Pi(\sigma_2\sigma_1)$.   Indeed, $\pi\ge \Pi(\sigma_1)$ and $\pi\ge \Pi(\sigma_2)$ so $\pi\ge \Pi(\sigma_1)\vee \Pi(\sigma_2)\ge \Pi(\sigma_2\sigma_1)$, so that  $\pi_1\vee\Pi(\sigma_2\sigma_1)\vee \pi_2=\pi\vee \Pi(\sigma_2\sigma_1) =\pi$. Therefore, the minimal value of $l$ such that $\gamma_l(\pi_2, \sigma_2)$ is non-vanishing  is reached for $l= \ell(\sigma_2, \pi_2, \pi_1\vee\Pi(\sigma_2\sigma_1))$, which here reads:
 \be
\label{eq:lowerboundislambda2}
 \ell(\sigma_2, \pi_2, \pi_1\vee\Pi(\sigma_2\sigma_1)) = n - \#(\sigma_2) + 2\bigl(\#(\pi_1\vee\Pi(\sigma_2\sigma_1) ) - \#(\pi )\bigr). 
 \ee
 
 Coming back to \eqref{eq:A}, we see that \emph{under the condition imposed in the sum}, \eqref{eq:lowerboundislambda2} can be re-written as
  \be
\label{eq:lowerboundislambda3}
 \ell(\sigma_2, \pi_2, \pi_1\vee\Pi(\sigma_2\sigma_1)) = n + \#(\sigma_2) - 2\#(\pi_2 ) , 
 \ee
 which is precisely the value taken by $l$ in \eqref{eq:rewrite-mu-M}.

\

 Consequently, from \eqref{eq:summing-M-min}: 
\be
\label{eq:A-2}
A[\sigma_2, \pi, \pi_1]= \sum_{\substack{{\pi_2 \ge \Pi(\sigma_2)}\\[+0.3ex]{\pi_1\vee \pi_2 = \pi, }\\[+0.2ex]{ \#(\sigma_2) - \#(\pi_2 )  =  \#(\pi_1\vee\Pi(\sigma_2\sigma_1)) -\#(\pi)}}} 
\prod_{G\in \pi_2} \gamma({\sigma_2}_{\lvert _G} ).
\ee
Since from \eqref{eq:id-pis}, for $\pi_1\ge \Pi(\sigma_1)$, $\pi_1\vee\Pi(\sigma_2\sigma_1)\vee \pi_2=\pi_1\vee \pi_2$, $A$ can actually be seen as depending only on $\Pi(\sigma_2)$, $\pi$,  and $\pi_1\vee\Pi(\sigma_2)$, and is given by $ \Gamma\bigl[\sigma_2, \pi, \pi_1\vee\Pi(\sigma_2\sigma_1)\bigr]$, which concludes the proof. \qed

\subsection{Direct derivation of the cumulant-moment formula}
\label{sub:first-ref-cum-direct}

The asymptotic  cumulant-moment formula in Thm.~\ref{thm:first-tree-moments} is obtained from  \eqref{eq:cum-mom-CMSS-prod}  by splitting the conditions defining the convolution (through Lemma~\ref{lemma:split}). The relation \eqref{eq:cum-mom-CMSS-prod} is itself derived in  \cite{CMSS} using the properties of the convolution to invert \eqref{eq:phi-in-conv}. In this subsection, we show that   \eqref{eq:final-higher-order-cum-mom} can instead be derived directly by taking the asymptotics of the expression for finite $N$, that is \eqref{eq:Cum-higher-order-fin-N-gen}, which we rewrite here:
$$
\kappa (\pi, \sigma) [\vec B]=  \sum_{\tau\in S_n} \sum_{\substack{{\pi'\in \mathcal{P}(n)}\\{\pi\ge \pi' \ge \Pi(\tau)}}} \varphi(\pi', \tau) [\vec B]  \sum_{\substack{{\pi''\in \mathcal{P}(n)}\\{\pi\ge \pi''  \ge \Pi(\sigma)\vee\pi'}}} \lambda_{\pi'', \pi} \prod_{G\in \pi'' } W^{(N)}\left((\sigma\tau^{-1})_{\lvert_{_G}}\right).
$$

\begin{lemma} 
\label{lem:explicit-expr-WC}
Let $\pi, \tilde \pi \in \mathcal{P}(n)$ and $\nu\in S_n$ such that $\pi \ge \tilde \pi \ge \Pi(\nu)$, then: 
\be
\label{eq:explicit-expr-WC}
\sum_{\substack{{\pi''\in \mathcal{P}(n)}\\{\pi\ge \pi''  \ge \tilde\pi}}} \lambda_{\pi'', \pi} \prod_{G\in \pi'' } W^{(N)}\bigl(\nu_{\lvert_{_G}}\bigr) = \sum_{l\ge 0}  N^{-n - l}  \sum_{\substack{{\pi''\ge \Pi(\nu) }\\{\tilde \pi\vee\pi'' = \pi } } }\ \gamma_l(\pi'', \nu) .
\ee
\end{lemma}
This can be seen from arguments in \cite{Collins03} and  \cite{CGL}, but we reprove it here.
\proof For $\pi''\ge \tilde \pi$, using \eqref{eq:expansion-wg-2} and Lemma~\ref{lem:summing-M}, we write:
$$
\prod_{G\in \pi'' } W^{(N)}\left(\nu{\lvert_{_G}}\right)  = N^{-n} \prod_{G\in \pi''}\Biggl( \sum_{\substack{\pi_G\in\mathcal{P}(\lvert G \rvert)\\\pi_G\ge \Pi(\nu{\lvert_{_G}})}}\sum_{l_G\ge 0} N^{-l_G} \sum_{\substack{\{l_G^H\ge0\}_{H\in \pi_G}\\ \sum_H l_G^H = l_G}}\prod_{H\in \pi_G} \gamma_{l_G^H}(\nu{\lvert_{_H}})\Biggr), 
$$
where we have used the fact that ${(\nu_{\lvert_{_G}})}_{\lvert{_H}}={\nu}_{\lvert_{_H}}$. We have:
$$
\prod_{G\in \pi'' } W^{(N)}\left(\nu{\lvert_{_G}}\right)  = N^{-n} \sum_{\left\{\substack{\pi_G\in\mathcal{P}(\lvert G \rvert)\\\pi_G\ge \Pi(\nu{\lvert_{_G}})}\right\}_{G\in \pi''}}  \sum_{\{l_G^H\ge0\}_{G, H}}N^{-\sum_{G,H}l_G^H} \prod_{G,H}\gamma_{l_G^H}(\nu{\lvert_{_H}}),
$$
where the sums and products with subscripts $G,H$ are for $G\in\pi''$, $H\in \pi_G$. We may group all the partitions $\pi_G$ into a partition $\pi=\sqcup_G \pi_G$ satisfying $\pi'' \ge \pi \ge \Pi(\nu)$. Reciprocally, any such permutation is subdivided uniquely into  the $\pi_G$s. Therefore:
$$
\prod_{G\in \pi'' } W^{(N)}\left(\nu{\lvert_{_G}}\right)= N^{-n} \sum_{\substack{\pi\in\mathcal{P}(n) \\ \pi''\ge \pi \ge \Pi(\nu)}}  \sum_{\{l^H\ge0\}_{H\in\pi}}N^{-\sum_{H\in\pi}l^H} \prod_{H}\gamma_{l^H}(\nu{\lvert_{_H}}) = \hspace{-0.2cm}\sum_{\substack{\pi\in\mathcal{P}(n) \\ \pi''\ge \pi \ge \Pi(\nu) }}  \sum_{l\ge0}\ N^{-n-l} \gamma_l(\pi, \nu),
$$
where we have used Lemma~\ref{eq:summing-M}. Grouping the terms according to the partition $\bar \pi = \pi \vee\tilde\pi$:  
$$
\prod_{G\in \pi'' } W^{(N)}\left(\nu{\lvert_{_G}}\right)   = \sum_{\substack{\bar \pi\in\mathcal{P}(n) \\ \pi''\ge \bar \pi \ge \tilde\pi }}\  \sum_{\substack{\pi\in\mathcal{P}(n) \\  \pi \vee \tilde \pi = \bar \pi  }}\  \sum_{l\ge0}\ N^{-n-l} \gamma_l(\pi, \nu). 
$$
Inverting this relation in the lattice of partitions, we get the desired result:
$$
\sum_{\substack{\pi\in\mathcal{P}(n) \\  \pi \vee \tilde \pi = \bar \pi  }}\  \sum_{l\ge0}\ N^{-n-l} \gamma_l(\pi, \nu) = \sum_{\substack{{\pi''\in \mathcal{P}(n)}\\{\pi\ge \pi''  \ge \tilde\pi}}} \lambda_{\pi'', \pi} \prod_{G\in \pi'' } W^{(N)}\left(\nu{\lvert_{_G}}\right).
$$
\qed 

\ 

We may now combine Lemmas~\ref{lem:explicit-expr-WC} and~\ref{lem:min-l} 
to obtain the dominant contribution in $N$. For $\pi \ge \tilde \pi \ge \Pi(\nu)$:
\be
\label{eq:expression-sum-lambda-wein}
\sum_{\substack{{\pi''\in \mathcal{P}(n)}\\{\pi\ge \pi''  \ge \tilde\pi}}} \lambda_{\pi'', \pi} \prod_{G\in \pi'' } W^{(N)}\bigl(\nu_{\lvert_{_G}}\bigr) = N^{\#(\nu) - 2 \left(\#(\tilde \pi ) - \#(\pi )\right)-2n} \Gamma[\nu, \pi, \tilde \pi](1+O(\frac1{N^2})), 
\ee
with $ \Gamma[\nu, \pi, \tilde\pi]$ corresponding to the coefficient for $l= \ell(\nu, \pi, \tilde\pi)$, given by \eqref{eq:comb-expr-leading-cum-Weing}. Inserting this in \eqref{eq:Cum-higher-order-fin-N-gen}:
$$
\kappa (\pi, \sigma) [\vec B]=  \hspace{-0.1cm} \sum_{\tau\in S_n}\hspace{-0.1cm} \sum_{\substack{{\pi'\in \mathcal{P}(n)}\\{\pi\ge \pi' \ge \Pi(\tau)}}}  \hspace{-0.3cm}N^{\#( \sigma\tau^{-1}) - 2 \left(\#( \Pi(\sigma)\vee\pi' ) - \#(\pi)\right)-2n} \varphi(\pi', \tau) [\vec B]  \  \Gamma\bigl[\sigma \tau^{-1}, \pi, \pi'\vee \Pi(\sigma)\bigr](1+o(1)).
$$
Recalling that $
 \varphi(\pi', \tau) [\vec  b] = \lim_{N\rightarrow \infty } N^{\#(\tau) - 2\#(\pi' )}  \varphi(\pi', \tau) [\vec  B],
$
we obtain:
\begin{align*}
&\kappa (\pi, \sigma) [\vec  B]=   \sum_{\tau\in S_n} \sum_{\substack{{\pi'\in \mathcal{P}(n)}\\{\pi\ge \pi' \ge \Pi(\tau)}}}  N^{2\#(\pi' ) - \#(\tau) + \#( \sigma\tau^{-1}) - 2 \left(\#( \Pi(\sigma)\vee\pi' ) - \#(\pi)\right)-2n} \varphi(\pi', \tau) [\vec  b]\\
&\hspace{9cm} \times \Gamma\bigl[\sigma \tau^{-1}, \pi, \pi'\vee \Pi(\sigma)\bigr](1+o(1)).
\end{align*}

For $\pi'\ge \Pi(\tau)$, we apply Lemma~\ref{lem:L} again:
 $$
L[\Pi(\sigma)\vee\Pi(\tau), \pi'; \Pi(\tau)] = \#(\tau) -  \#(\Pi(\sigma)\vee\Pi(\tau)) -  \#(\pi' )  + \#(\Pi(\sigma)\vee\pi')  \ge 0.
 $$
Furthermore, the  Euler characteristics of the bipartite map $(\sigma, \tau^{-1})$ reads \eqref{eq:genus}:
 $$
\#(\sigma\tau^{-1}) = n + 2\#(\Pi(\sigma)\vee\Pi(\tau) ) - 2g(\sigma, \tau^{-1}) - \#(\sigma) - \#(\tau),
$$
where $g(\sigma, \tau^{-1})\ge 0$, with equality if and only if the bipartite map $(\sigma, \tau^{-1})$ is planar. 
Therefore:
\begin{align*}
&\kappa (\pi, \sigma) [\vec  B]=    N^{2\#(\pi)-n-\#(\sigma)} \sum_{\tau\in S_n} \sum_{\substack{{\pi'\in \mathcal{P}(n)}\\{\pi\ge \pi' \ge \Pi(\tau)}}}  N^{-2L[\Pi(\sigma)\vee\Pi(\tau), \pi'; \Pi(\tau)] - 2g(\sigma, \tau^{-1}) } \varphi(\pi', \tau) [\vec  b] \\
&\hspace{9cm} \times \Gamma\bigl[\sigma \tau^{-1}, \pi, \pi'\vee \Pi(\sigma)\bigr](1+o(1)).
\end{align*}
We thus recover that $\kappa (\pi, \sigma) [\vec  B]$ is of order $N^{2\#(\pi )-n -\#(\sigma)}$, and that  \eqref{eq:final-higher-order-cum-mom}. \qed

\subsection{Reformulation using trees}
\label{sec:bij-with-trees}

From the proof of Lemma~\ref{lem:L}, it is clear that the expression in Thm.~\ref{thm:first-tree-moments}  involves a sum over bipartite trees. Let us make this more explicit, while restricting to the case where  $\pi=1_n$ (which is still general by multiplicativity), and $\vec b =(b, \ldots, b)$.

Considering $\lambda, \lambda'\vdash n$, we let $ \mathcal{M}(\lambda ; \lambda') $ be the number of planar and connected  bipartite maps with $n$ labeled edges, with white vertices encoded by the permutation $\gamma_\lambda = \gamma_{\lambda_1, \ldots, \lambda_p}\in S_n$, $p=\#(\lambda)$, and such that the degrees of the black vertices induce a partition $\lambda'\vdash n$: 
\be
 \mathcal{M}(\lambda ; \lambda')  = \mathrm{Card}\Bigl\{   \tau\in S_n \, \mid \, g(\gamma_\lambda , \tau)=0, \,   \Pi(\gamma_\lambda)\vee\Pi(\tau)=1_n,\, \Lambda(\tau) = \lambda' \Bigr\}. 
\ee

Consider a partition $\pi\in \mathcal{P}(n)$. To label its blocks \emph{by the minimum}, we label by 1 the block  $B_1$ which has the element 1, by 2 the block $B_2$ which has the smallest element not in $B_1$, by 3 the block $B_3$ which has the smallest element not in $B_1\cup B_2$, and so on.

We recall that $\Lambda_{\{d_a\}_a}$ is the partition of $\sum_a a d_a$ that has $d_a$ parts equal to $a$, and that if $T\in \mathcal{T}_p$,  $\mathcal I'_i(T)$ and $\mathcal{I}'(T)$ are the subsets of $\mathcal I_i(T)$ and $\mathcal{I}(T)$ whose elements $I$ satisfy $\lvert I \rvert >1$.

\begin{theorem}
\label{th:first-tree-formula}
 For $\lambda\vdash n$,
 $\vec b =(b, \ldots, b)$, and $\pi=1_n$, we may rewrite \eqref{thm:first-tree-moments} as:
\begin{align*}
& \varphi_p (b^{\lambda_1}, \ldots, b^{\lambda_p})=\sum_{k= 1}^{p}\ 
\sum_{\substack{{\tilde \pi \ge \Pi(\gamma_{\lambda})}\\{\#(\tilde \pi) = k}}}
\sum_{T\in \mathcal{G}_k}
\Biggl(\sum_{i=1}^k \sum_{l_i\ge 0}\sum_{\substack{{\{j_i^I\ge 1\}_{I\in \mathcal{I}_i(T)}}\\{k_1^{(i)}, \ldots, k_{\lvert G_i \rvert}^{(i)}\ge 0,\ \ \sum_{a}k_a^{(i)} = l_i}\\{\sum_{I\in \mathcal{I}_i(T)} j_i^I + \sum_{a=1} ^ {\lvert G_i \rvert} a k_a^{(i)}  = {\lvert G_i \rvert} }}}\Biggr)\Biggl(\prod_{I\in\mathcal{I}(T)}\  \kappa_{\{j_i^I\}_{i\in I}}[b] \Biggr)
 \\[+2ex]&\hspace{6.5cm}\times
\left(\prod_{i=1}^k  \mathcal{M}\Bigl(\lambda_{\lvert _i} ; \Lambda_{\{k_a^{(i)} + q_a^{(i)}\}_a} \Bigr)\prod_{a=1}^{\lvert G_i \rvert} \frac{\bigl(k_a^{(i)} + q_a^{(i)}\bigr)!}{k_a^{(i)} !}\, \kappa_a[b]^{k_a^{(i)}}    \right),
\end{align*}
where for each $\tilde \pi$ in the sum, we label its blocks by the minimum as $G_1, \ldots, G_k$;  we have denoted by 
$\lambda_{\lvert _i}= \Lambda\bigl((\gamma_\lambda){\lvert_{_{G_i}}}\bigr)$; and $q_a^{(i)}$ is the number of integers in the set $\{j_i^I\}_{I\in \mathcal{I}_i(T)} $ that are equal to $a$. 
\end{theorem}

\proof Our starting point is \eqref{eq:mom-cum-CMSS-facto-rewrite}. We fix $\sigma, \tau\in S_n,\, g(\sigma, \tau^{-1})=0$, and let $k= \#(\Pi(\sigma)\vee \Pi(\tau))$. We label the blocks of $\Pi(\sigma)\vee \Pi(\tau)$ by the minimum as $G_1, \ldots, G_k$ and let 
$$\tau_i = \tau_{\lvert_{G_i}}.$$

Let $\mathcal{T}(\sigma, \tau)$ be the set of labeled bipartite connected trees in $\mathcal{T}_k$ with white vertices labeled from 1 to $k$, such that the vertex number $i$ has valency $\mathrm{deg}_T(i)= \#(\tau_i)$. As explained in Sec.~\ref{sub:graphs}, $T$ is bijectively encoded by its set $\mathcal{I}(T)$ satisfying 
\be
\label{eq:excess-lab-tree-rapl}
\sum_{I\in \mathcal{I}(T)} (\lvert I \rvert - 1) - p + 1 =0.
\ee 
For $1\le i \le k$, define $\mathcal{I}_i'(T)=\{I \in \mathcal{I}_i(T)\, \mid \, \lvert I \rvert >1\}$.  We establish a bijection  between the set
\be
\label{eq:Psigmatau}
\mathcal{P}(\sigma, \tau) = \bigl\{\pi' \ge \Pi(\tau),\, \mid\,  \pi'\vee \Pi(\sigma)= 1_n,\, L[\pi', \Pi(\sigma)\vee \Pi(\tau) ; \Pi(\tau) ]  \bigr\}, 
\ee
and the set of triplets $(T, \{\mathcal L_i\}_{1\le i \le p}, \mathcal L ^{\mathsf L})$, where $T\in \mathcal{T}(\sigma, \tau)$ and $ \{\mathcal L_i\}_{1\le i \le p}$, $\mathcal L ^{\mathsf L}$ are \emph{``bijective'' attributions of the blocks of $\tau$ to the edges of $T$}, that is, functions of the following kind:
\begin{align}
\label{eq:attributions}
\\[-0.6cm]
\nonumber
\begin{array}{ll}
\ \bullet\   \mathcal L_i : \mathcal{I}_i'(T) \rightarrow \Pi(\tau_i), \ \forall i\in \{1, \ldots, k\},  \\[-0.5ex] \\
\ \bullet \  \mathcal L ^{\mathsf L} : \{1, \ldots, k\} \rightarrow \mathcal{P}(\Pi(\tau_i)), \textrm{ where an element of }\mathcal{P}(\Pi(\tau_i))\textrm{ is a subset of the blocks of }\Pi(\tau_i), 
 \end{array}
\end{align}
such that a block $V\in \Pi(\tau_i)$ is found exactly one time,  either as an element of the set $\mathcal L ^{\mathsf L}(i)$, either as $\mathcal L_i(I)$ for a unique  $I\in \mathcal{I}_i'(T)$:
$$
\forall V\in \Pi(\tau_i), \ V\notin \mathcal L ^{\mathsf L}(i) \ \Rightarrow\  \exists ! I\in \mathcal{I}_i'(T),\ V=\mathcal L_i(I). 
$$
Note that the black leaves incident to the same white vertex  are not distinguishable and are therefore treated as a group: to the set of leaves incident to the white vertex number $i$ is attributed a group of as many blocks of $\Pi(\tau_i)$.

\begin{figure}[!h]
\centering
\includegraphics[scale=0.65]{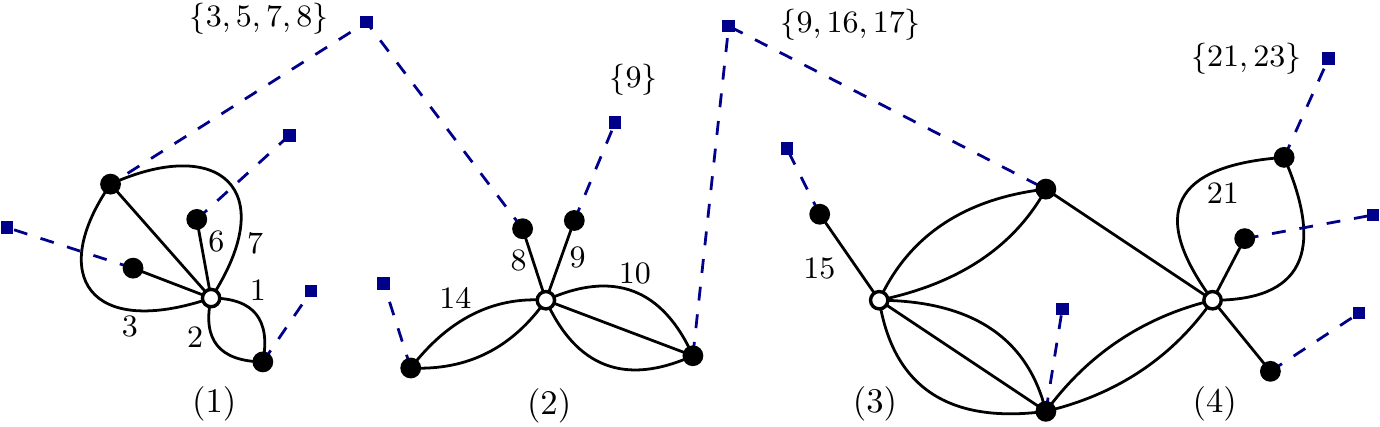}\hspace{0.6cm}\raisebox{6ex}{$\rightarrow$}\hspace{0.6cm}\raisebox{1.5ex}{\includegraphics[scale=0.45]{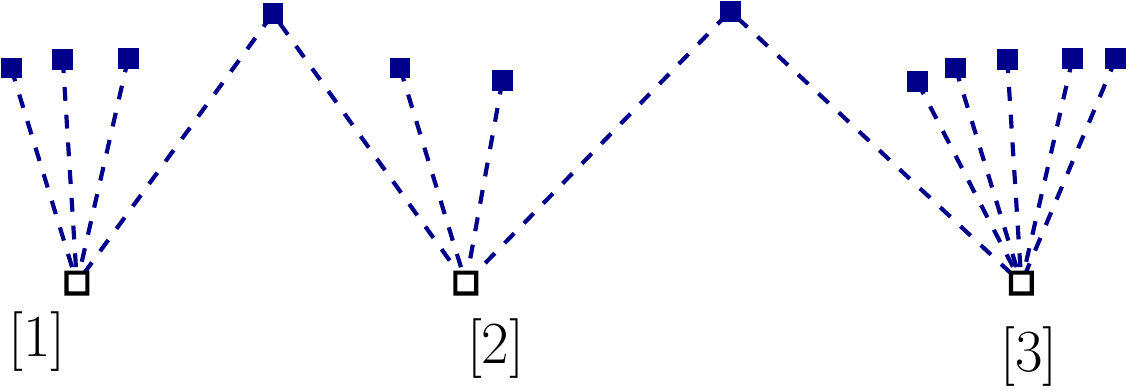}}
\caption{Left: an element of $\mathcal{P}(\gamma_{\lambda}, \tau)$ for $\lambda= \{7, 7, 6, 7\}$ and $\tau$ defining a planar bipartite map $(\gamma_{\lambda}, \tau)$. The blue vertices represent the blocks of $\pi'$. Right: the corresponding tree in    $\mathcal{T}(\gamma_{\lambda}, \tau)$.
}
\label{fig:Bij-arbres}
\end{figure}

\noindent{\it From trees to partitions.} Consider $(T, \{\mathcal L_i\}_{1\le i \le p}, \mathcal L ^{\mathsf L})$ as above, and define
\be
\pi' = \bigl\{  \cup_{i\in I} \mathcal L _i(I) \ \mid\  I\in \mathcal{I}_i'(T) \bigr\} \cup \left(\cup_{1=i}^k \mathcal L ^{\mathsf L}(i) \right),
\ee
that is, to each black vertex $I\in \mathcal{I}(T)$ we associate a block of $\pi'$:  if $\lvert I \rvert >1$ then the corresponding 
block of $\pi'$ is the union of the blocks of $\Pi(\tau)$ attributed to the edges incident to that black vertex $I$. As for the black leaves, this is done in a non-distinguishable manner: all the blocks of $\Pi(\tau)$ associated to black leaves are also blocks of $\pi'$. From the ``bijectivity'' of the block attributions: 
\be
\Pi(\tau) = \bigl\{   \mathcal L _i(I) \ \mid\ I\in \mathcal{I}'(T), i\in I \bigr\} \cup \left(\cup_{1=i}^k \mathcal L ^{\mathsf L}(i) \right), 
\ee
and  $\pi' \ge \Pi(\tau)$. If $i\in I \cap I'$, where $I, I'\in \mathcal{I}'(T)$ then $G_i$ has elements both in $ \cup_{i\in I} \mathcal L _i(I)$ and $ \cup_{i\in I'} \mathcal L _i(I)$, and obvious modifications for black leaves, so the fact that $T$ is connected translates as $\pi' \vee (\Pi(\sigma \vee \Pi(\tau)))=1_n$. 

The fact that  $ \#(\pi) -  \#(\tau) - \#(\pi' )  +  \#(\Pi(\tau) \vee \Pi(\sigma) ) = 0 $ follows from \eqref{eq:excess-lab-tree-rapl}, noticing that $\#(\tau) = \sum_{I\in \mathcal{I}(T)} \lvert I \rvert$ and $\#(\pi') = \sum_{i} (\lvert \mathcal{I}'(T)\rvert + \mathrm{deg}_T^{\mathsf L}(i))$. 

\

\noindent{\it From partitions  to trees.} Starting from $\pi' \in \mathcal{P}(\sigma, \tau)$, we first define $T$ through $\mathcal{I}(T)$. To each block of $\pi'$, we associate 
\be
\label{eq:IB}
I_B = \{i\in \{1, \ldots, k\}\, \mid \, B\cap G_i \neq \emptyset\}, 
\ee
and define $\mathcal{I}(T) = \{I_B\, \mid \,  B\in \pi'\}$. We verify that it defines an element of $\mathcal{T}(\sigma, \tau)$. It is seen to be connected because $\pi' \vee (\Pi(\sigma \vee \Pi(\tau)))=1_n$ and noticing that : 
$$
i \in I_B \cap I_{B'}\  \Leftrightarrow\  \exists\, a,b\in G_i,\ a\in B, b\in B'.
$$
If $L[\pi', \Pi(\sigma)\vee \Pi(\tau) ; \Pi(\tau) ] =0 $, a block of $\pi'$ and a block $G_i$ of $\Pi(\sigma)\vee \Pi(\tau)$ may intersect only if their intersection is precisely a block of $\Pi(\tau)$: 
\be 
\label{eq:unicity-of-intersection}
\mathrm{If}\, L[\pi', \Pi(\sigma)\vee \Pi(\tau) ; \Pi(\tau) ] =0,  \ \forall B \in \pi', \, [B\cap G_i\neq \emptyset \  \Rightarrow \ \exists  V\in \Pi(\tau), B\cap G_i = V]. 
\ee 
From this we see that  
$\mathrm{deg}_T(i)= \lvert \{B\in \pi'\, \mid \,  B\cap G_i \neq \emptyset\}\rvert$ is in fact the number of blocks of $\Pi(\tau)$ in $G_i$, that is, 
$\#(\tau_i)$. It also implies that $\lvert I_B\rvert= \#(\tau_{\lvert_B})$, and $\sum_{I\in \mathcal I(T)} \lvert I \rvert = \#(\tau)$. 
so that the condition $L[\pi', \Pi(\sigma)\vee \Pi(\tau) ; \Pi(\tau) ] =0$ translates to \eqref{eq:excess-lab-tree-rapl}. 

It remains to define the block attributions $ \mathcal L _i$ and $\mathcal L ^{\mathsf L}$: 
\begin{itemize}
\item $ \mathcal L_i : I_B \textrm{ s.t. } \lvert I_B\rvert >1 \rightarrow B\cap G_i, \ \forall i\in \{1, \ldots, k\}, $ (which, has we have seen, corresponds to a unique block of $\Pi(\tau_i)$),
\item $\mathcal L ^{\mathsf L} : i \rightarrow \{V\in \pi' \cap \Pi(\tau),\, V\in G_i\}$.
\end{itemize}
One must verify that the attributions thus defined are ``bijective'': let $V\in \Pi(\tau_i)$. Since $\pi'\ge \Pi(\tau_i)$, $\exists ! B\in \pi',\, V\subset B\cap G_i$, so from \eqref{eq:unicity-of-intersection}, $V =  B\cap G_i$. We distinguish two cases: if $\lvert I_B\rvert >1$, from the definition we indeed have $V=\mathcal L_i ( I_B)$. Otherwise $\lvert I_B\rvert =1$, so $B$ has empty intersections with all others $G_j, j\neq i$, so that $B\cap G_i =B$. Since $V =  B\cap G_i$, $V=B$, so that $V\in \pi'$. Since $V\in \Pi(\tau)$ and $V\subset G_i$, $V\in \mathcal L ^{\mathsf L}(i)$. 

\

Having established this bijection, we may reformulate \eqref{eq:mom-cum-CMSS-facto-rewrite} as follows for $\vec b =(b, \ldots, b)$, $\sigma=\gamma_\lambda, \lambda \vdash n$ and $\pi=1_n$: 
\be
 \varphi_p (b^{\lambda_1}, \ldots, b^{\lambda_p}) =
\sum_{\substack{{\tau\in S_n,}\\[+0.3ex]{g(\gamma_\lambda, \tau^{-1})=0}}} 
\sum_{T \in \mathcal{T}(\gamma_\lambda, \tau)} 
\sum_{\substack{{\{\mathcal L_i\}_{1\le i \le p},\ \mathcal L ^{\mathsf L}}\\{\textrm{bij. attributions}}}} 
\prod_{I \in \mathcal{I}'(T)}  \kappa_{\{\lvert \mathcal L_i(I)\rvert\}_{i\in I}}[b] \prod_{i=1}^k \prod_{V\in \mathcal L ^{\mathsf L}(i)}  \kappa_{\lvert V\rvert}[b], 
\ee
where $k=\#(\Pi(\gamma_\lambda)\vee \Pi(\tau))$. We group the terms in the sum by the partition $\tilde \pi := \Pi(\gamma_\lambda)\vee \Pi(\tau)$:
\begin{align*}
& \varphi_p (b^{\lambda_1}, \ldots, b^{\lambda_p})  =\sum_{k=1}^p 
\sum_{\substack{{\tilde \pi \ge \Pi(\gamma_{\lambda}),}\\[+0.2ex]{\#(\tilde \pi)=k}}} 
 \sum_{\substack{{\tau\in S_n,}\\[+0.2ex]{g(\gamma_\lambda, \tau^{-1})=0}\\[+0.3ex]{\Pi(\gamma_\lambda)\vee \Pi(\tau) =\tilde \pi}}} 
\sum_{\substack{{T \in \mathcal{T}_k}\\{\mathrm{deg}_T(i)= \#(\tau_i)}}} 
\sum_{\substack{{\{\mathcal L_i\}_{i},\ \mathcal L ^{\mathsf L}}\\{\textrm{bij. att.}}}} \\&\hspace{7cm}\times
\prod_{I \in \mathcal{I}'(T)}  \kappa_{\{\lvert \mathcal L_i(I)\rvert\}_{i\in I}}[b] \prod_{i=1}^k \prod_{V\in \mathcal L ^{\mathsf L}(i)}  \kappa_{\lvert V\rvert}[b], 
\end{align*}
where we label the blocks of $\tilde \pi$ by the minimum as $G_1, \ldots, G_k$, and $\tau_i = \tau_{\lvert _{G_i}}$. We can now exchange the sums over $T$ and $\tau$: 
\be
\label{eq:useful-sum-with-attributions}
 \varphi_p (b^{\lambda_1}, \ldots, b^{\lambda_p}) =\sum_{k=1}^p 
\sum_{\substack{{\tilde \pi \ge \Pi(\gamma_{\lambda}),}\\[+0.2ex]{\#(\tilde \pi)=k}}} 
\sum_{T \in \mathcal{T}_k} 
 \sum_{\substack{{\tau\in S_n,}\\[+0.2ex]{g(\gamma_\lambda, \tau^{-1})=0}\\[+0.3ex]{\Pi(\gamma_\lambda)\vee \Pi(\tau) =\tilde \pi}\\[+0.3ex]{\#(\tau_i) = \mathrm{deg}_T(i)}}} 
\sum_{\substack{{\{\mathcal L_i\}_{i},\, \mathcal L ^{\mathsf L}}\\{\textrm{bij. att.}}}} \ 
\prod_{I \in \mathcal{I}'(T)}  \kappa_{\{\lvert \mathcal L_i(I)\rvert\}_{i\in I}}[b] \prod_{i=1}^k \prod_{V\in \mathcal L ^{\mathsf L}(i)}  \kappa_{\lvert V\rvert}[b], 
\ee
and then group the terms by the cardinals of the $\mathcal L_i(I)$ for $I,i$, and of the blocks of $\mathcal L ^{\mathsf L}(i)$, in order to factorize the product of cumulants:
\begin{align*}
& \varphi_p (b^{\lambda_1}, \ldots, b^{\lambda_p}) =\sum_{k=1}^p 
\sum_{\substack{{\tilde \pi \ge \Pi(\gamma_{\lambda}),}\\[+0.2ex]{\#(\tilde \pi)=k}}} 
\sum_{T \in \mathcal{T}_k} 
\sum_{i=1}^k
\sum_{\substack{{\{j_i^I\ge 1\}_{I\in \mathcal{I}'_i(T)}}\\{k_1^{(i)}, \ldots, k_{\lvert G_i \rvert}^{(i)}\ge 0,\ \ \sum_{a}k_a^{(i)} =   \mathrm{deg}_T^{\mathsf L}(i) }\\{\sum_{I\in \mathcal{I}_i(T)} j_i^I + \sum_{a=1} ^ {\lvert G_i \rvert} a k_a^{(i)}  = {\lvert G_i \rvert} }}}
\\&\hspace{3cm}
\times \prod_{I \in \mathcal{I}'(T)}  \kappa_{\{j_i^I\}_{i\in I}}[b] \prod_{i=1}^k \prod_{a=1}^{\lvert G_i \rvert}  \kappa_{a}[b]^{k_a^{(i)}}
 \sum_{\substack{{\tau\in S_n,}\\[+0.2ex]{g(\gamma_\lambda, \tau^{-1})=0}\\[+0.3ex]{\Pi(\gamma_\lambda)\vee \Pi(\tau) =\tilde \pi}\\[+0.4ex]{\Lambda(\tau_i) =  \Lambda_{\{k_a^{(i)}+q_a^{(i)}\}_a}}}} 
\sum_{\substack{{\{\mathcal L_i\}_{i},\ \mathcal L ^{\mathsf L}}\\{\textrm{bij. att.}}\\{\forall i \in \mathcal{I}_i'(T),\, \lvert \mathcal L_i(I)\rvert = j_i^I}\\{\forall i,\, \mathcal L ^{\mathsf L}(i)\textrm{ has }k_a^{(i)}\textrm{ elts. of size } a}}} 1, 
\end{align*}
where we recall that $\mathrm{deg}_T^{\mathsf L}(i) $ is the number of black leaves incident to the white vertex $i$, and $q_a^{(i)}$ is the number of $j_i^I$ equal to $a$ when $I$ spans $\mathcal I_i'(T)$. The number of attributions is found to be:
\be
\sum_{\substack{{\{\mathcal L_i\}_{i},\ \mathcal L ^{\mathsf L}}\\{\textrm{bij. att.}}\\{\forall i \in \mathcal{I}_i'(T),\, \lvert \mathcal L_i(I)\rvert = j_i^I}\\{\forall i,\, \mathcal L ^{\mathsf L}(i)\textrm{ has }k_a^{(i)}\textrm{ elts. of size } a}}} 1 \ = \ \prod_{i=1}^k\prod_{a=1}^{\lvert G_i \rvert} \binom{k_a^{(i)}+q_a^{(i)}}{k_a^{(i)}} q_a^{(i)}!\, .
\ee
Indeed, the blocks of $\Pi(\tau_i)$ are fixed and we know from the condition on $\Lambda(\tau_i)$ that for each $a$, there are $k_a^{(i)}+q_a^{(i)}$ edges incident to the white vertex $i$ that are attributed a block of $\Pi(\tau_i)$ of size $a$. We have to count all possible ways to attribute for every $a,$ $k_a^{i}$ of them to the set of black leaves (without providing the information on which black leaf as they can't be distinguished), and among the $q_a{(i)}$ remaining ones, to decide to which non-leaf black vertex they are assigned to (the edges incident to black leaf vertices are distinguishable). There are $q_a{(i)}!$ ways to do so, thus the result above. This quantity factorizes as it is  independent from $\tau$, and it remains to see that denoting by $\gamma^i_\lambda=(\gamma_\lambda){\lvert_{_{G_i}}}$:
\be
\sum_{\substack{{\tau\in S_n,}\\[+0.2ex]{g(\gamma_\lambda, \tau^{-1})=0}\\[+0.3ex]{\Pi(\gamma_\lambda)\vee \Pi(\tau) =\tilde \pi}\\[+0.4ex]{\forall i,\, \Lambda(\tau_i) =  \Lambda_{\{k_a^{(i)}+q_a^{(i)}\}_a}}}}  1 \ =\  \prod_{i=1}^k \sum_{\substack{{\tau_i\in S_n,}\\[+0.2ex]{g(\gamma^i_\lambda, \tau_i^{-1})=0}\\[+0.3ex]{\Pi(\gamma^i_\lambda)\vee \Pi(\tau_i) = 1_{\lvert G_i\rvert}}\\[+0.4ex]{\Lambda(\tau_i) =  \Lambda_{\{k_a^{(i)}+q_a^{(i)}\}_a}}}} 1 \ =\  \prod_{i=1}^k \mathcal M\left(\lambda_{\lvert _i} ; \Lambda_{\{k_a^{(i)}+q_a^{(i)}\}_a}\right), 
\ee
where we recall that $\lambda_{\lvert _i}= \Lambda\bigl((\gamma_\lambda){\lvert_{_{G_i}}}\bigr)$. The theorem follows by replacing the sum over trees  $T\in\mathcal T_k$ by a sum over  trees $T'\in \mathcal G_k$ (with no black leaves)  and a sum over $l_i\ge 0$, where $\mathrm{deg}_T^{\mathsf L}(i) = l_i$, and $\mathcal I'(T) = \mathcal I(T')$.
\qed 

\

As an example, for $p=2$  there are two trees in $\mathcal{G}_2$ and the formula in the theorem reads:
\begin{align*}
& \varphi_2 (b^{\lambda_1}, b^{\lambda_2}) =  \sum_{\substack{{\{k_i\ge 0\}_{1\le i \le n}}\\[+0.2ex]{\sum_{i=1}^{n} i k_i =n  }}}  \prod_{i=1}^n \kappa_i[b]^{k_i}  \mathcal{M}\Bigl({\{\lambda_1, \lambda_2\}} ; \Lambda_{\{k_a^{(i)} \}_a}\Bigr)
\\&  +     \sum_{\substack{{j_1, j_2\ge 1}}}\kappa_{\{j_1, j_2\}}[b] 
 \prod_{i=1}^{2} \Biggl(\sum_{\substack{{k_1^{(i)}, \ldots, k_{\lambda_i}^{(i)} \ge 0}\\{\sum_{a=1}^{\lambda_i} a k_a^{(i)} =\lambda_i - j_i}}} \frac{(k_{j_i}^{(i)} + 1)!}{k_{j_i}^{(i)}!}  \prod_{a=1}^{\lambda_i} \kappa_a[b]^{k_a^{(i)}} \frac{\lambda_i!}{(\lambda_i  - \sum_a k_a^{(i)})!\prod_{a\neq j} k_a! (k_{j_i}+1)!}\Biggr)     ,
\end{align*}
where $\sigma_i = \sigma_{\lvert_{G_i}}$, and we have used the fact that (see e.g.\cite{Spei94})
\be
\label{eq:unicellular-maps-fixed-partition}
\mathcal{M}\bigl(n; \Lambda_{\{k_a \}_a}\bigr) = \frac{n!}{(n +1 - \sum_a k_a )!\prod_a k_a^{(i)}!}.
\ee
This simplifies to 
\begin{align*}
& \varphi_2 (b^{\lambda_1}, b^{\lambda_2}) =  \sum_{\substack{{\{k_i\ge 0\}_{1\le i \le n}}\\[+0.2ex]{\sum_{i=1}^{n} i k_i =n  }}}  \prod_{i=1}^n \kappa_i[b]^{k_i}  \mathcal{M}\Bigl({\{\lambda_1, \lambda_2\}}  ; \Lambda_{\{k_a^{(i)} \}_a}\Bigr)
\\& \hspace{3cm} +     \sum_{\substack{{j_1, j_2\ge 1}}}\kappa_{\{j_1, j_2\}}[b] 
 \prod_{i=1}^{2} \Biggl(\sum_{\substack{{k_1^{(i)}, \ldots, k_{\lambda_i}^{(i)} \ge 0}\\{\sum_{a=1}^{\lambda_i} a k_a^{(i)} =\lambda_i - j_i}}}  \prod_{a=1}^{\lambda_i} \kappa_a[b]^{k_a^{(i)}} \frac{\lambda_i!}{(\lambda_i  - \sum_a k_a^{(i)})!\prod_{a} k_a! }\Biggr) .
\end{align*}

\ 

The bijection between partitions satisfying a forest-condition $L=0$ and labeled bipartite trees in the proof of Thm.~\ref{th:first-tree-formula} will be used in several occasions in Sec.~\ref{sec:funct-rel}, and may take simpler forms, depending on the cases and choices of variables in the sums. Let us formalize it in a simpler form. For $\pi \in \mathcal{P}(n)$, we let  $\mathcal{G}_{\pi}$ be the set of bipartite trees with  white vertices labeled by the blocks of $\pi$ (to avoid labeling these blocks). The black vertices - given by the set $\mathcal{I}(T)$ - are now sets of strictly more than one blocks of $\pi$. We consider some  series of numbers $c(\{i_j\})_{1\le j \le q}$ for $q\ge 1$ and all $i_j \ge 1$. 
\begin{corollary} 
\label{cor:bij-with-trees-form}
Consider $\bar \pi, \tilde \pi \in \mathcal{P}(n)$ satisfying $\tilde \pi \ge \bar \pi$ and let $k=\#(\tilde \pi \vee \bar \pi)$. The following identity holds:
$$
\sum_{\substack{{\pi''\ge \bar \pi}\\{\tilde\pi\vee\pi''= 1_n }\\[+0.3ex]{ L[\pi'', \,\tilde \pi \, ; \, \bar \pi]=0}}} 
\hspace{-0.1cm}\prod_{G\in \pi''} c(\{\lvert B \rvert\}_{B\in \bar \pi}) =  \sum_{T\in \mathcal{G}_k}
  \sum_{G\in \tilde \pi \vee \bar \pi}\hspace{-0.1cm}\sum_{\substack{{\{B_G^I\}_{I\in \mathcal{I}_G(T)}
}\\[+0.2ex]{\mathrm{s.t.}\, B_G^I\in \bar \pi_{\lvert_G}}\\[+0.3ex]{\mathrm{2-by-2\  different}}}} 
\hspace{-0.1cm} \prod_{I\in\mathcal{I}(T)}\  c({\{\lvert B_G^I \rvert \}_{G\in I}}) 
 \prod_{G\in \tilde \pi \vee \bar \pi}  \prod_{\substack{{B \in \bar \pi}\\{B\notin \{B_G^I\}_{I\in \mathcal{I}_G}} }} \hspace{-0.1cm} c(\lvert B \rvert).
  $$
\end{corollary}
In particular, we may apply this  to reformulate  \eqref{eq:comb-expr-leading-cum-Weing}, which we recall to be:
\be
\Gamma[\nu, 1_n, \tilde\pi] = 
\sum_{\substack{{\pi''\ge \Pi(\nu)}\\{\tilde\pi\vee\pi''= 1_n }\\[+0.3ex]{ L[\pi'', \,\tilde \pi \, ; \, \Pi(\nu)]=0}}} \ 
\prod_{G\in \pi''} \gamma(\nu_{|_G}) \;.
\ee
This can thus be expressed as:
$$
\Gamma[\nu, 1_n, \tilde\pi] = \sum_{T\in \mathcal{G}_k}
  \sum_{G\in \tilde \pi \vee \Pi(\nu)}\sum_{\substack{{\{B_G^I\}_{I\in \mathcal{I}_G(T)}
}\\[+0.2ex]{\mathrm{s.t.}\, B_G^I\in \Pi(\nu_{\lvert_G})}\\[+0.3ex]{\mathrm{2-by-2\  different}}}} 
 \prod_{I\in\mathcal{I}(T)}\  \gamma({\{\lvert B_G^I \rvert \}_{G\in I}}) 
  \prod_{\substack{{B \in \Pi(\nu)}\\{B\notin \{B_G^I\}_{I\in \mathcal{I}_G}} }} C_{\lvert B \rvert - 1} ,
$$
where $C_k=\binom{2k}{k}/(k+1)$ is the $k$th Catalan number, and $\gamma$  are rescaled monotone Hurwitz numbers:
$$
\gamma(\{q_G\}_{G\in I})  = (-1)^{\lvert I \rvert + 1} \frac{(2 \sum_{G\in I} q_G + \lvert I \rvert - 3)!}{(2 \sum_{G\in I} q_G)!}\prod_{G\in I }\frac{(2 q_G )!}{ q_G! (q_G-1)!},
$$
where $\lvert I \rvert$ is the valency of the black vertex $I$.

\section{Functional relations}
\label{sec:funct-rel}

\subsection{Functional relations for factorized permutations}
\label{sec:factorized}

Instead of pursuing with \eqref{th:first-tree-formula}, where the number of planar bipartite maps  with prescribed vertex valencies  $\mathcal{M}$ causes difficulties, we rather follow the strategy of \cite{CMSS}, providing a second tree reformulation from factorized permutations.

\

 Recall that a function $f : \pi \ge \sigma \rightarrow f(\pi, \sigma)$  is said to be \emph{multiplicative} if $f(\pi, \sigma) = \prod_{G\in \pi} f(1_{{\lvert G\rvert}}, \sigma_{\mid{_{\lvert G\rvert}}})$, and $f(1_n, \sigma)$ depends only on $\Lambda(\sigma)$.  Also recall that 
 $$
L[\pi', \Pi(\sigma)\vee \Pi(\tau) ; \Pi(\tau) ] =  \#(\pi) -  \#(\tau) - \#(\pi' )  +  \#(\Pi(\tau) \vee \Pi(\sigma) ).$$

If $\lambda \vdash n$, $\gamma_\lambda = \gamma_{\lambda_1, \ldots, \lambda_p}\in S_n$ is given in \eqref{eq:gamma-cycles}, where $p=\#(\lambda)$, and we denote by
\be
\mathrm{NC}(\lambda_1)\times \cdots\times \mathrm{NC}(\lambda_p) = \Bigl\{ \tau \in S_n
\, \mid \, \Pi(\tau)\le \Pi(\gamma_\lambda) ,\, \textrm{and for }  1\le i \le p,\ g\bigl(\tau_i, \gamma_i^{-1}\bigr)=0\Bigr\},
\ee
where $\tau_i$ and $\gamma_i$ are the restrictions of $\tau$ and $\gamma$ to $\{\sum_{j<i} \lambda_j + 1, \ldots, \sum_{j<i} \lambda_j + \lambda_i \}$ (which is the $i$th block of $\Pi(\gamma_\lambda)$ when labeled by the minimum). We also use the notation $\tau = \tau_1 \times \ldots \times \tau_p$,  $\tau_i\in  \mathrm{NC}(\lambda_i)$.
The following is a generalization at arbitrary order of Prop.~6.2 in \cite{CMSS}.

\begin{theorem}
\label{thm:functional-relations-for-any-H}
Let $g : \pi \ge \sigma \rightarrow g(\pi, \sigma)$ be a multiplicative function and 
$p\ge 1$, and 
consider the generating function
$$
H_p(X_1, \ldots, X_p) = \delta_{p, 1} + \sum_{k_1, \ldots, k_p \ge 1} \tilde \alpha _{k_1, \ldots, k_p} X_1^{k_1}\cdots X_p^{k_p}, \qquad \ \   \tilde \alpha _{k_1, \ldots, k_p} = g(1_{k_1 + \cdots + k_p}, \gamma_{k_1, \ldots, k_p}).
$$
Introduce the following numbers: for $p\ge 1$ and $\lambda_1, \ldots, \lambda_p \ge 1$,
\be
\label{eq:relation-NC-factorisee}
\beta_{\lambda_1, \ldots, \lambda_p} =
\sum_{\substack{{\tau\in\mathrm{NC}(\lambda_1)\times \cdots\times \mathrm{NC}(\lambda_p) 
}}} 
\hspace{0.2cm}
\sum_{\substack{{\pi' \ge \Pi(\tau)}\\[+0.3ex]{\pi'\vee\Pi(\gamma_{\lambda})= 1_n}\\[+0.3ex] {L[\pi',\, \Pi(\gamma_{\lambda})\vee \Pi(\tau)\, ; \,\Pi(\tau)] =0
}}} g(\pi', \tau),
\ee
and the corresponding generating function:
$$
B_p(X_1, \ldots, X_p)= \delta_{p, 1} + \sum_{\lambda_1, \ldots, \lambda_p \ge 1} \beta _{\lambda_1, \ldots, \lambda_p} X_1^{\lambda_1}\cdots X_p^{\lambda_p}.
$$
Then we have the following functional relations:
\be
B_p(X_1, \ldots, X_p)  = \sum_{T\in \mathcal{G}_p} \biggl( \prod_{i=1}^p 
  X_i^{\mathrm{d}_T(i) -1} \frac {\mathrm d^{\mathrm{d}_T(i) -1}} {\mathrm d X_i^{\mathrm{d}_T(i) -1}}\biggr)\Biggl\{\prod_{i=1}^p
   {\frac {\mathrm{d} Y_i} {\mathrm d X_i}} \prod_{I\in \mathcal{I}(T)} \ \frac{H_{\lvert I \rvert}\bigl(Y_I\bigr)}{\prod_{i\in I}B_1(X_i)}\Biggr\},   
\ee
where $\mathrm{d}_T(i) = \mathrm{deg}_T(i)$,  $Y_I=\{Y_i\}_{i\in I}$, and
$Y_i = X_i B_1(X_i)$. 
\end{theorem}

At the level of the coefficients, these relations are equivalent to
\begin{align*}
& \beta_{\lambda_1, \ldots, \lambda_p} =
\sum_{T\in \mathcal{G}_p}
\Biggl(\sum_{i=1}^p \sum_{l_i\ge 0}\sum_{\substack{{\{j_i^I\ge 1\}_{I\in \mathcal{I}_i(T)}}\\{k_1^{(i)}, \ldots, k_{\lambda_i}^{(i)}\ge 0,\ \ \sum_{a}k_a^{(i)} = l_i}\\{\sum_{I\in \mathcal{I}_i(T)} j_i^I + \sum_{a=1} ^ {\lambda_i} a k_a^{(i)}  = {\lambda_i} }}}\Biggr)\Biggl(\prod_{I\in\mathcal{I}(T)}\  \tilde \alpha_{\{j_i^I\}_{i\in I}} \Biggr)
 \\&\hspace{7.5cm}\times
\left(\prod_{i=1}^n   \frac{\lambda_i!}{(\lambda_i + 1 -  l_i  - \mathrm{deg}_T(i))!} \prod_{a=1}^{\lambda_i} \frac{1}{k_a^{(i)} !}\  \tilde \alpha_a ^{k_a^{(i)}}    \right). 
\end{align*}

\proof We can adapt the proof of Th.~\ref{th:first-tree-formula}. The difference lies in the fact that here, $\Pi(\tau)\vee \Pi(\gamma_\lambda) = \Pi(\gamma_\lambda)$, so that the bijection reduces to trees in $\mathcal T_p$ with one white vertex per block of $\Pi(\gamma_\lambda)$.  Using in addition the formula \eqref{eq:unicellular-maps-fixed-partition}:
 \begin{align*}
& \beta_{\lambda_1, \ldots, \lambda_p} =
\sum_{T\in \mathcal{G}_p}
\Biggl(\sum_{i=1}^p \sum_{l_i\ge 0}\sum_{\substack{{\{j_i^I\ge 1\}_{I\in \mathcal{I}_i(T)}}\\{k_1^{(i)}, \ldots, k_{\lambda_i}^{(i)}\ge 0,\ \ \sum_{a}k_a^{(i)} = l_i}\\{\sum_{I\in \mathcal{I}_i(T)} j_i^I + \sum_{a=1} ^ {\lambda_i} a k_a^{(i)}  = {\lambda_i} }}}\Biggr)\Biggl(\prod_{I\in\mathcal{I}(T)}\  \tilde \alpha_{\{j_i^I\}_{i\in I}} \Biggr)
 \\&\hspace{3cm}\times
\left(\prod_{i=1}^n   \frac{\lambda_i!}{(\lambda_i + 1 - \sum_a k_a^{(i)}+ q_a^{(i)})\prod_{a} (k_a^{(i)} + q_a^{(i)})! } \prod_{a=1}^{\lambda_i} \frac{\bigl(k_a^{(i)} + q_a^{(i)}\bigr)!}{k_a^{(i)} !}\, \tilde \alpha_a ^{k_a^{(i)}}    \right),
\end{align*}
where $q_a^{(i)}$ is the number of integers in the set $\{j_i^I\}_{I\in \mathcal{I}_i(T)} $ that are equal to $a$. In particular, $\sum_a q_a^{(i)} = \mathrm{deg}_T(i)$. 
The derivation of the functional formulas  is unchanged from Sec.~\ref{subsub:simplificaion-funct-relations-gen}. There are only two minor points to be justified. The first one is that nothing changes with the sums starting at $l_i\ge0$ instead of $l_i \ge 1$. One verifies that the equation \eqref{eq:sum-the-l} is still valid for $l\ge0$ (and it no longer matters whether $\lambda>j$).
The second point to be justified is that we can indeed use Lagrange inversion as in \eqref{eq:Lagrange}, which is true because  the relation \eqref{eq:relation-NC-factorisee} is $\beta_{n} =
\sum_{\tau\in \mathrm{NC}(n)} g(\Pi(\tau), \tau)$, meaning (see the main theorem in \cite{Spei94}) that $H_1(X B_1(X))= B_1(X)$. \qed

\subsection{Factorization}
\label{sec:factorization}

Let us recall here the relations of Thm.~\ref{thm:first-tree-moments} for $\vec b =(b, \ldots, b)$, $\sigma=\gamma_\lambda, \lambda \vdash n$ and $\pi=1_n$: 
\be
 \varphi_p (b^{\lambda_1}, \ldots, b^{\lambda_p}) =
\sum_{\substack{{\tau\in S_n,}\\[+0.3ex]{g(\gamma_\lambda, \tau^{-1})=0}}}  \hspace{0.2cm}
\sum_{\substack{{\pi' \ge \Pi(\tau)}\\[+0.3ex]{\pi'\vee\Pi(\gamma_\lambda)= 1_n}\\[+0.3ex] {L[\pi',\, \Pi(\gamma_\lambda)\vee \Pi(\tau)\, ; \,\Pi(\tau)] =0
}}}  \kappa(\pi', \tau)[ b] .
\ee
Considering  $f : \pi \ge \sigma \rightarrow g(\pi, \sigma)$ a multiplicative function and 
$p\ge 1$, and
some numbers defined for $p\ge 1$ and $\lambda_1, \ldots, \lambda_p \ge 1$ by (see the left of Fig.~\ref{fig:Bij-arbres}):
\be
\label{eq:assumption-facto}
\beta_{\lambda_1, \ldots, \lambda_p} =
\sum_{\substack{{\tau\in S_n,}\\[+0.3ex]{g(\gamma_\lambda, \tau^{-1})=0}}}
\hspace{0.2cm}
\sum_{\substack{{\pi' \ge \Pi(\tau)}\\[+0.3ex]{\pi'\vee\Pi(\gamma_{\lambda})= 1_n}\\[+0.3ex] {L[\pi',\, \Pi(\gamma_{\lambda})\vee \Pi(\tau)\, ; \,\Pi(\tau)] =0
}}} f(\pi', \tau),
\ee
 for instance  $\beta_{\lambda_1, \ldots, \lambda_p} := \varphi_p(b^{\lambda_1}, \ldots, b^{\lambda_p})$ and $f(1_{k_1 + \cdots + k_p}, \gamma_{k_1, \ldots, k_p}) := \kappa_{k_1, \ldots, k_p}[b]$,  we can always opt for a summation over factorized permutations $\mathrm{NC}(\lambda_1)\times \cdots\times \mathrm{NC}(\lambda_p)$ by a generalization of the trick of Sec.~6 of \cite{CMSS}.

If $\tau = \tau_1 \times \ldots \times \tau_p$, where $\tau_i\in  \mathrm{NC}(\lambda_i)$ and $B$ a block of $\pi'\ge \Pi(\tau)$, we let as in  \eqref{eq:IB}:
\be
I_B = \{i\in \{1, \ldots, p\}\, \mid \, B\cap \Pi(\tau_i) \neq \emptyset\}. 
\ee
If $I\subset \{1, \ldots, p\}$, we let $\mathcal{P}^*(I)$ be the set of partitions of the elements in $I$ that have more than $\lvert I \rvert$ blocks.  
For $d\ge 1$ and $\tau\in S_d$, we let $\mathrm{NS} (\tau)$ be the subset 
$$
\mathrm{NS} (\tau)\subset \Bigl\{\nu \in S_{d} \mid  g(\nu, \tau)=0, \Pi(\tau)\vee \Pi(\nu) = 1_{d}\Bigr\}$$
\emph{for which there are no white cut-vertices}, which we have called planar non-separable hypermaps (see Sec.~\ref{sub:bip-map}).
Note that if 
$\tau=\eta$ has just one cycle, $\mathrm{NS} (\eta)=\{\eta^{-1}\}$ has just one element.

\begin{theorem}
\label{thm:facto}
The relations \eqref{eq:assumption-facto} defining  $\beta_{\lambda_1, \ldots, \lambda_p}$ can be expressed equivalently as:
\be
\label{eq:split-perm-1}
\beta_{\lambda_1, \ldots, \lambda_p} =
\sum_{\substack{{\tau =\bigtimes_{i=1}^p \tau_i}\\[+0.3ex]{\tau_i\in \mathrm{NC}(\lambda_i)}}}\hspace{0.2cm}
\sum_{\substack{{\pi' \ge \Pi(\tau)}\\[+0.3ex]{\pi'\vee\Pi(\gamma_{\lambda})= 1_n}\\[+0.3ex] {L[\pi',\, \Pi(\gamma_{\lambda})\, ; \,\Pi(\tau)] =0
}}} g (\pi', \tau),
\ee
where $g(\pi', \tau)$ is multiplicative over the blocks of $\pi'$, and such that if $B\in \pi'$:
\be
\label{eq:split-perm-2}
g (1_{\lvert B\rvert }, \tau_{\lvert_B}) = f\bigl(\Pi( \tau_{\lvert_B}),  \tau_{\lvert_B}\bigr) + \sum_{\bar \pi \in \mathcal{P}^*(I_B)} \bar g _{\bar \pi }  (\tau_{\lvert_B}),
\ee
where   
for $G\in \bar \pi$, we define $\tau_G=\bigtimes_{i\in G} \bigl({\tau_i}_{\lvert_{B} 
}\bigr)$ (which is also  $\tau_G={\gamma_\lambda}_{\lvert_G}$), 
and:
\be
\label{eq:bar-f-gen}
\bar g _{\bar \pi }  (\tau_{\lvert_B}) = \sum_{\substack{{\nu =\bigtimes_{G\in \bar \pi}\nu_G}\\[+0.3ex]{\nu_G\in \mathrm{NS} (\tau_G)}}} \ 
\sum_{\substack{{\pi''  \ge \Pi(\nu)}\\[+0.3ex]{\pi'' \vee\Pi(\tau_{\lvert_B})= 1_{\lvert B \rvert}}\\[+0.3ex] {L[\pi'',\, \Pi(\tau_{\lvert_B})\vee \Pi(\nu)\, ; \,\Pi(\nu)] =0}}} 
f(\pi'', \nu). 
\ee
\end{theorem}
 See the end of the present subsection for the proof. If $\bar\pi$ has a single block, this simplifies to: 
\be
\label{eq:bar-f-one-block}
\tilde f (\tau_{\lvert_B})  := \bar g_{1_{\lvert I_B \rvert}}  (\tau_{\lvert_B}) = \sum_{\nu\in \mathrm{NS}\left(\tau_{\lvert_B}\right)} \ 
 f(\Pi(\nu), \nu). 
\ee

\paragraph{Weight of a vertex as sum over trees.}We now set  $f(1_{k_1 + \cdots + k_p}, \gamma_{k_1, \ldots, k_p})= \kappa_{k_1, \ldots, k_p}[b]$ and \emph{drop the explicit dependency in} $[b]$ from now on. Let $\mathcal{G}_{\bar \pi}$ be the set of bipartite trees with  white (square) vertices labeled by the blocks of $\bar \pi$ (to avoid labeling these blocks with integers), so that the non-white vertices of the tree, given by the set $\mathcal{I}(T)$ - now blue squares to avoid confusion - are sets of strictly more than one blocks of $\bar \pi$. From \eqref{eq:useful-sum-with-attributions} in the proof of Thm.~\ref{th:first-tree-formula} (see also Cor.~\ref{cor:bij-with-trees-form}), we may rewrite \eqref{eq:bar-f-gen} as:
\be
\bar g _{\bar \pi }    =
\sum_{T \in \mathcal{G}_{\bar \pi}} 
\bar \kappa _{\bar \pi, T}
\ee
where recalling that $\tau_G=\bigtimes_{i\in G} \bigl({\tau_i}_{\lvert_{B}}\bigr)$ (which is also  $\tau_G={\gamma_\lambda}_{\lvert_G}$),
\be
\label{eq:bar-f-tree-1}
\bar \kappa _{\bar \pi, T}  (\tau_{\lvert_B})  =
\sum_{\substack{{\nu =\bigtimes_{G\in \bar \pi}\nu_G}\\[+0.3ex]{\nu_G\in \mathrm{NS} (\tau_G)}\\[+0.3ex]{\#(\nu_G) \ge \mathrm{deg}_T(G) }}}
\
\sum_{\substack{{\{\mathcal L_G\}_{G}
}\\{\textrm{inj. att.}}}} \ 
\prod_{K \in \mathcal{I}(T)}  \kappa_{\{\lvert \mathcal L_G(K)\rvert\}_{G\in K}}  \ \prod_{G\in \bar \pi}\prod_{\substack{{V\in \Pi(\nu_G)}\\{V\notin \mathcal L_G( \mathcal{I}_G(T))} }}  \kappa_{\lvert V\rvert} , 
\ee
in which  $\mathcal{I}(T)=\cup_G \mathcal{I}_G(T)$, and  $\mathcal{I}_G(T)$ is the set of blue square vertices of the tree  incident to the white square vertex $G$, and for each white square vertex $G\in \bar \pi$, $\mathcal L_G$ is a choice of \emph{injective} attribution of the edges incident to that vertex to the cycles of $\nu_G$:
$$
\forall G\in \bar \pi, \ \mathcal L_G : \mathcal{I}_G(T) \rightarrow \Pi(\nu_G), 
$$
which can also be written as in  Cor.~\ref{cor:bij-with-trees-form}. The sums over $k$ and $\tilde \pi$ in \eqref{eq:useful-sum-with-attributions} are no longer needed, as $\tilde \pi$ is fixed to $\bar \pi = \Pi(\nu)\vee\Pi(\tau)$, and we have chosen to sum over trees in $\mathcal G_p$ instead of $\mathcal T_p$, so that the attributions are only injective. 

For $T$ the tree with one blue vertex $\lvert \mathcal{I}(T)\rvert=1$, we denote by $\bar \kappa _{\bar \pi }  : = \bar \kappa _{\substack{{\bar \pi, T},\, \mathrm{s.t.} \, {\lvert \mathcal{I}(T)\rvert=1}}}  $. This sum then simplifies to:
\be
\label{eq:bar-f-one-tree}
\bar \kappa _{\bar \pi }  (\tau_{\lvert_B})  = \sum_{\substack{{\nu =\bigtimes_{G\in \bar \pi}\nu_G}\\[+0.3ex]{\nu_G\in \mathrm{NS} (\tau_G)}}} \ \sum_{\substack{{F=\bigcup_{G\in \bar \pi} F_G}\\{F_G\in \Pi(\nu_G)} }} \  \kappa(1_{\lvert F \rvert}, \nu_{\lvert _{F}}) \prod_{\substack{{V \in \Pi(\nu)}\\{V \notin \{F_G\}} }}  \kappa(1_{\lvert V \rvert}, \nu_{\lvert _{V }}). \ee

If $\tau_G=\eta$ consists of a single cycle, then  $\mathrm{NS} (\tau_G)=\{\eta^{-1}\}$, and the condition $\#(\nu_G) \ge \mathrm{deg}_T(G)$ imposes $\mathrm{deg}_T(G)=1$.

\

\begin{figure}[!h!]
\centering
\includegraphics[scale=0.7]{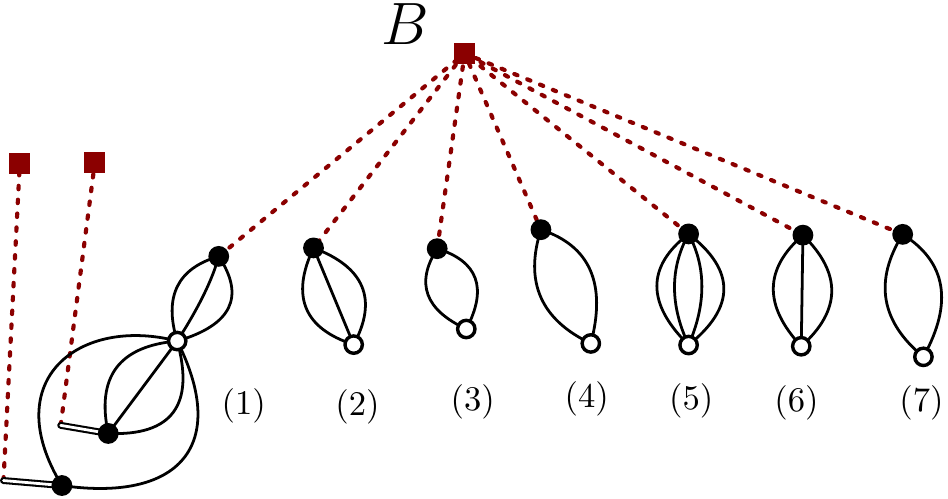}

\includegraphics[scale=0.7 ]{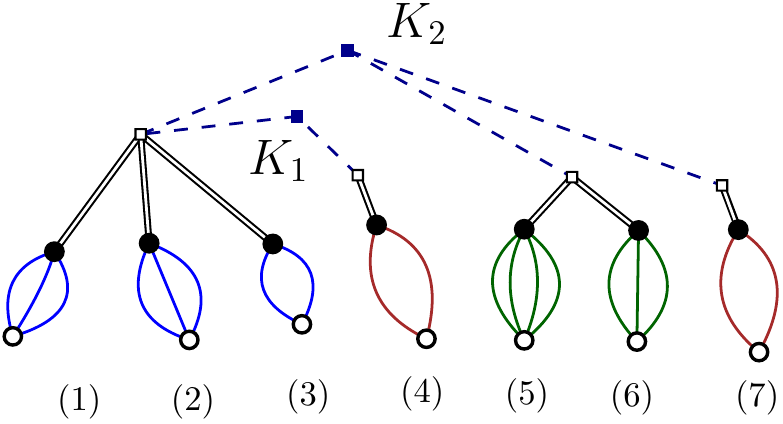}\hspace{2cm}\includegraphics[scale=0.7]{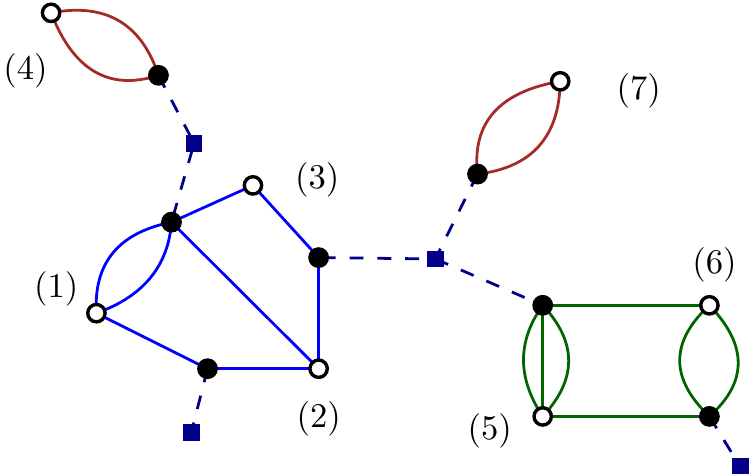}
\caption{Top: a bipartite map $(\gamma_\lambda, \tau^{-1})$ with factorized $\tau =\bigtimes_{i=1}^p \tau_i$, $\tau_i\in \mathrm{NC}(\lambda_i)$ (the black and white vertices and the edges linking them), and tree-like partition $\pi'$ with three blocks (the three red vertices). Bottom left: Focusing on the block $B$ of $\pi'$, a choice is made of $\bar\pi=\{1,2,3\}\{4\}\{5,6\}\{7\}$ (edges in the same blocks carry the same color) and $T$ with two blue vertices $K_1=\{1,2,3\}\{4\}$ and $K_2=\{1,2,3\}\{5,6\}\{7\}$. Bottom right:  choice of a collection of non-separable hypermaps $\{\nu_G\}$, one for each block of $\bar \pi$, and injective attribution of the edges of $T$ to the cycles of the $\{\nu_G\}$. The two cycles which have not been attributed any edge of the tree are attributed a blue leaf. The blue vertices of valency $q$ correspond to order-$q$ free cumulants. 
}
\label{fig:first-expl}
\end{figure}

To put it in words and pictures,  the cost of using factorized permutations $\tau = \tau_1 \times \ldots \times \tau_p$, where $\tau_i\in  \mathrm{NC}(\lambda_i)$ (top of Fig.~\ref{fig:first-expl}) as in \eqref{eq:split-perm-1}, instead of using generic permutations $\tau$ such that $g(\tau, \gamma_\lambda^{-1})=0$, as in \eqref{eq:assumption-facto} (left of Fig.~\ref{fig:Bij-arbres}), is to replace the weight  $f\bigl(\Pi( \tau_{\lvert_B})\bigr)$ for the blocks $B\in \pi'$ in \eqref{eq:assumption-facto} (indicated by the red square vertices in Fig.~\ref{fig:first-expl}) by a new weight \eqref{eq:split-perm-2}, which is the sum of the usual weight and a sum of corrections. 
The corrections for the block $B$ are labeled by $\bar \pi \in \mathcal{P}^\star(I_B)$ (thicker edges and colors on the bottom left of Fig.~\ref{fig:first-expl}) and $T\in \mathcal{G}_{\bar\pi}$, where $I_B$ is the set of white disc vertices (the cycles of $\gamma_\lambda$) involved in this block $B$.  The corrections for the partition $\bar \pi$ re-create the connections between the  white vertices grouped in $\bar \pi$, connections that we have  removed in factorizing $\tau$. These connections take the form of planar non-separable hypermaps (bottom right of Fig.~\ref{fig:first-expl}). For a fixed $\bar \pi$, the sum over trees  (dashed blue edges on the left of Fig.~\ref{fig:first-expl}) indicates  how the higher-order cumulants are arranged. The sum over $\mathcal L_G$ corresponds to choosing which black vertices of the non-separable hypermaps are grouped in order-$q$ free cumulants (one for each blue vertex with $q$ incident edges). On the bottom right of  Fig.~\ref{fig:first-expl} for instance, the contribution is 
$$
\kappa_{2, 2, 4}[b] \kappa_{2, 4}[b] \kappa_{3}[b] \kappa_{2}[b] . 
$$

As in Thm.~\ref{th:first-tree-formula}, we can provide a more detailed expression, useful in Sec.~\ref{sub:small-orders} for concrete combinatorial computations. Given a planar  bipartite map with edges labeled from 1 to $n$ and white vertices labeled with an integer in $G\subset \{1, \ldots, p\}$, we provide a labeling of the black vertices as follows. First consider  the white vertex $v_1$ with smallest label in $G$, and label from 1 to $\mathrm{deg}(1)$ the black vertices to which $v_1$ is connected, according to the smallest label of the edges that connect them. Then go to the white vertex $v_1$ with second smallest label in $G$, and label in the same way the black vertices to which $v_2$ is connected and that haven't been labeled yet, and so on, until all the black vertices are labeled.

Letting $I\subset \{1, \ldots, p\}$ and  $\mu_I=\{\mu_i\}_{i\in I}$, $\mu_i\ge 1$, as well as $\bar \pi \in \mathcal{P}^\star(I)$, $G\in \bar \pi$ and  $\mu_G=\{\mu_i\}_{i\in G}$, we denote by
$$
\mathrm{NS} (\mu_G):=\mathrm{NS} (\gamma_{\mu_G}).
$$
To a bipartite planar map with white vertices given by $\gamma_{\mu_G}$ (thus labeled by $G$) and $r$ black vertices given by $\nu_G$ (so that $\#(\nu_G)=r$) and labeled as above, we can therefore  associate its adjacency matrix 
$$
\mathcal{A}(\nu_G)=\{j_{i, s}\}_{\substack{{i\in G}\\{1\le s \le r}}},
$$ 
where $j_{i, s}$ is the number of edges connecting the white vertex $i$ and the black vertex $s$. We define the cardinal of the set of planar non-separable hypermaps with fixed permutation $\gamma_{\mu_G}$ for the white vertices and fixed adjacency matrix as:
$$
\mathcal{N}_{G, r}\left(\{j_{i,s}\}_{\substack{{i\in G}\\{1\le s \le r}}}\right) = \mathrm{Card}\left\{ \nu_G\in\mathrm{NS} (\mu_G),\, \mu_G=\Bigl(\sum_{s} j_{i,s}\Bigr)_{i\in G}, \, \mathcal{A}(\nu_G) =  \{j_{i,s}\}_{i,s}\right\}.
$$

With these notations and fixing $T\in \mathcal G _{\bar \pi}$, we can now follow steps analogous to the proof of Thm.~\ref{th:first-tree-formula} (see also Cor.~\ref{cor:bij-with-trees-form}) and rewrite \eqref{eq:bar-f-tree-1} as: 
\begin{align}
\label{eq:combinatorial-expr-correction-with-trees}
&\bar \kappa _{\bar \pi, T}  (\gamma_{\mu_I} )   = 
\Biggl(\sum_{G\in \bar \pi}\sum_{r_G\ge \mathrm{deg}_T(G)}
\ 
 \sum_{\substack{{\{j_{i,s}\ge 0\}_{\substack{{i\in G},\, {1\le s \le r_G}}}}\\[+0.2ex]{\mathrm{s.t.}  \sum_{1\le s \le r_G} j_{i,s} = \mu_i}}}
 \ 
  \sum_{\substack{{\{\theta({G,K})\}_{K\in \mathcal{I}_G(T)}
}\\[+0.2ex]{\mathrm{s.t.}\, 1\le \theta({G,K})\le r_G}\\[+0.3ex]{\mathrm{2-by-2\  different}}}}  \Biggr)
\\&\nonumber\hspace{3cm} \Biggl(\prod_{K \in \mathcal{I}(T)}  \kappa_{\left\{j_{G}^{\theta(G, K)}\right\}_{G\in K}}  \Biggr) \prod_{G\in \bar \pi}\Biggl(\mathcal{N}_{G, r_G}\left(\{j_{i,s}\}_{\substack{{i\in G}\\{1\le s \le r_G}}}\right) \prod_{\substack{{1\le s \le r_G}\\{s\notin \{\theta(G, K)\}_{K\in \mathcal{I}_G}} }}  \kappa_{j_G^s} \Biggr), 
\end{align}
where $j_G^s = \sum_{i\in G} j_{i, s}$. The sum over $\theta(G,K)$ corresponds to the sum over attributions $\mathcal L_G$. This formula will be useful for looking at concrete examples in Sec.~\ref{sub:small-orders}.

\paragraph{Generating functions.} Having now a sum over factorized $\tau$, we can apply Thm.~\ref{thm:functional-relations-for-any-H}. Fixing $I\subset \{1, \ldots, p\}$, we therefore introduce the generating function $H_{\lvert I \rvert}\bigl(\{Y_i\}_{i\in I}\bigr)$ of Thm.~\ref{thm:functional-relations-for-any-H}
as
\be
\label{eq:def-HI}
H_{\lvert I \rvert} = C_{\lvert I \rvert} + \sum_{\bar \pi \in \mathcal{P}^*(I)} \bar H_{\lvert I \rvert, \bar \pi}, \qquad \bar H_{\lvert I \rvert, \bar \pi} = \sum_{T \in \mathcal{G}_{\bar \pi}} 
\bar C_{\lvert I \rvert, \bar \pi, T},
\ee
where:
\be
\label{eq:bar-C-T}
\bar C_{\lvert I \rvert, \bar \pi, T} \bigl(\{Y_i\}_{i\in I}\bigr) = \sum_{\{\mu_i\}_{i\in I}} \  \prod_{i\in I} Y_i^{\mu_i} \Biggl[\hspace{-0.1cm} \sum_{\substack{{\nu =\bigtimes_{G\in \bar \pi}\nu_G}\\[+0.3ex]{\nu_G\in \mathrm{NS} (\mu_G)}\\[+0.3ex]{\#(\nu_G) \ge \mathrm{deg}_T(G) }}}
\sum_{\substack{{\{\mathcal L_G\}_{G}
}\\{\textrm{inj. att.}}}} \ 
\prod_{K \in \mathcal{I}(T)}  \kappa_{\{\lvert \mathcal L_G(K)\rvert\}_{G\in K}}  \ \prod_{G\in \bar \pi}\hspace{-0.1cm} \prod_{\substack{{V\in \Pi(\nu_G)}\\{V\notin \mathcal L_G( \mathcal{I}_G(T))} }}\hspace{-0.3cm}  \kappa_{\lvert V\rvert}  \Biggr] , 
\ee
In the case where $T$ is the tree with one blue vertex $\lvert \mathcal{I}(T)\rvert=1$, we let:
$\bar C_{\lvert I \rvert, \bar \pi} : = \bar C _{\lvert I \rvert, \substack{{\bar \pi, T},\, \mathrm{s.t.} \, {\lvert \mathcal{I}(T)\rvert=1}}} $:
\be
\label{eq:bar-C}
\bar C_{\lvert I \rvert, \bar \pi} \bigl(\{Y_i\}_{i\in I}\bigr) = \sum_{\{\mu_i\}_{i\in I}} \  \prod_{i\in I} Y_i^{\mu_i}\ \Biggl[\sum_{\bigl\{\nu_G\in \mathrm{NS} (\mu_G)\bigr\}_{G\in \bar \pi}} \ 
\sum_{\{F_G\in \Pi(\nu_G)\}_{G\in \bar \pi}} 
 \kappa_{\{\lvert F_G\rvert\}_{G\in \bar \pi}}  
 \prod_{\substack{{V \in \Pi(\bigtimes_G \nu_G)}\\{V \notin \{F_G\}} }} \kappa_{\lvert V \rvert} \Biggr] , 
\ee
which involves only free cumulants of order $\#(\bar \pi)$ and one. If $\bar\pi$ has a single block, we let $\tilde C_{\lvert I \rvert}  := \bar C_{\lvert I \rvert, \bar \pi} $, which is the generating function for non-separable hypermaps with prescribed vertex valencies, whose hyperedges are weighted by free cumulants $\kappa_q$: 
\be
\label{eq:bar-C-one-block}
\tilde C_{\lvert I \rvert} \bigl(\{Y_i\}_{i\in I}\bigr) := \bar C_{\lvert I \rvert, \bar \pi} \bigl(\{Y_i\}_{i\in I}\bigr) = \sum_{\{\mu_i\}_{i\in I}} \  \prod_{i\in I} Y_i^{\mu_i}\ \Biggl(\sum_{\nu\in \mathrm{NS} (\mu_I)} \ 
 \kappa\bigl(\Pi(\nu), \nu\bigr) \Biggr) , 
\ee
where we recall that $\kappa\bigl(\Pi(\nu), \nu\bigr)  = \prod_{G\in \Pi(\nu)} \kappa_{\lvert G \rvert}$.

\paragraph{Proof of the theorem.}We can rewrite \eqref{eq:split-perm-2} with a sum over partitions $\pi_B\ge \Pi(\tau_{\lvert_B})$ instead of $\bar \pi \in \mathcal{P}(I_B):$
\be
g (\pi', \tau)= \prod_{B\in \pi'} \sum_{\pi_B \ge \Pi(\tau_{\lvert_B}) }   \sum_{\substack{{\nu_B =\bigtimes_{G\in \pi_B}\nu^B_G}\\[+0.3ex]{\nu^B_G\in \mathrm{NS} (\tau_{\lvert_G})}}} \ 
\sum_{\substack{{\pi''_B  \ge \Pi(\nu_B)}\\[+0.3ex]{\pi''_B \vee\Pi(\tau_{\lvert_B})= 1_{\lvert B \rvert}}\\[+0.3ex] {L[\pi''_B,\, \Pi(\tau_{\lvert_B})\vee \Pi(\nu_B)\, ; \,\Pi(\nu_B)] =0
}}} f(\pi''_B, \nu_B). 
\ee

For each block $B$, the picture is that of the bottom of Fig.~\ref{fig:first-expl}: $B$ is a subset of the edges of the bipartite map $(\gamma_{\lambda}, \tau^{-1})$, $\pi_B$ provides groups of edges of the same color, to each group is associated a non-separable hypermap, and $\pi''_B$ is a tree -like partition of the black vertices  of the non-separable hypermaps  (cycles of $\nu_B$). 

The partitions $\{\pi_B\}$ form a partition $\pi$ in $\mathcal{P}(p)$ which satisfies $\Pi(\tau) \le \pi \le \pi'$. The permutations $\nu_B$ also form a permutation $\nu$ in $S_n$ which factorizes over the blocks of $\pi$. 
\be
g (\pi', \tau)= \sum_{\substack{{\pi \in \mathcal{P}(n)}\\{\pi'\ge \pi \ge \Pi(\tau)}} }\    \sum_{\substack{{\nu =\bigtimes_{G\in \pi}\nu_G}\\[+0.3ex]{\nu_G\in \mathrm{NS} (\tau_{\lvert_G})}}} \ 
\ 
\sum_{\substack{{\pi''\in \mathcal P(n)}\\{\pi' \ge \pi''  \ge \Pi(\nu)}\\[+0.3ex]{\forall\, B\in \pi',\ \pi''_{\lvert_B} \vee\Pi(\tau_{\lvert_B})= 1_{\lvert B \rvert}}\\[+0.3ex] {L[\pi'',\, \Pi(\tau)\vee \Pi(\nu)\, ; \,\Pi(\nu)] =0
}}}   f(\pi'', \nu). 
\ee
Note that a consequence of $\nu =\bigtimes_{G\in \pi}\nu_G$, $\nu_G\in \mathrm{NS} (\tau_{\lvert_G})$ is that $\pi=\Pi(\nu)\vee\Pi(\tau)$.

We now want to exchange the  summations over $\pi'$ and $\pi$ in \eqref{eq:split-perm-1}. To this aim, it must be noted that if $\pi' \ge \pi \ge \Pi(\tau)$ and $L[\pi',\, \Pi(\gamma_{\lambda})\, ; \,\Pi(\tau)] =0$, then $L[\pi,\, \Pi(\gamma_{\lambda})\, ; \,\Pi(\tau)] =0$, as from the definition in Lemma~\ref{lem:L}:
\[
 L[\pi,\, \Pi(\gamma_{\lambda})\, ; \,\Pi(\tau)] =   L[\pi',\, \Pi(\gamma_{\lambda})\, ; \,\Pi(\tau)] - L[\pi', \pi \vee \Pi(\gamma_{\lambda}) ; \pi ]
\]
where $L[\pi', \pi \vee \Pi(\gamma_{\lambda}) ; \pi ]=  \#(\pi) -\#(\pi')-  \#(\pi \vee \Pi(\gamma_{\lambda})) + \#(\pi' \vee \Pi(\gamma_{\lambda}))\ge 0$, so that   
if $L[\pi',\, \Pi(\gamma_{\lambda})\, ; \,\Pi(\tau)] =0$:
$$
0\le  L[\pi,\, \Pi(\gamma_{\lambda})\, ; \,\Pi(\tau)] = - L[\pi', \pi \vee \Pi(\gamma_{\lambda}) ; \pi ] \le 0. 
$$

\begin{figure}[!h]
\centering
\includegraphics[scale=0.7]{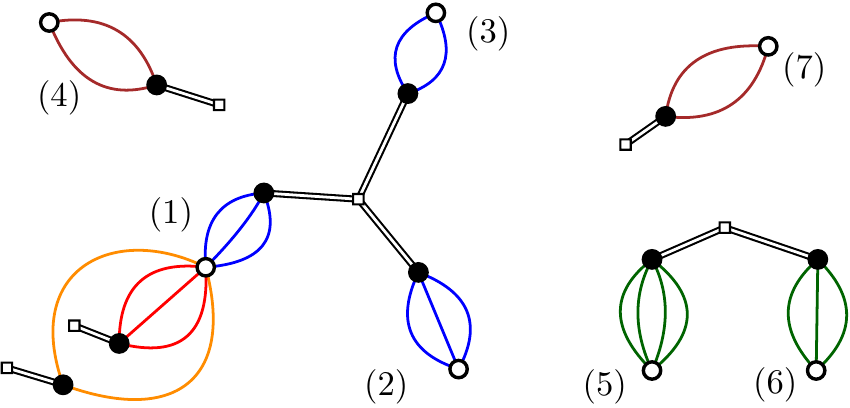}\hspace{1cm}\raisebox{6ex}{;}\hspace{1cm}\includegraphics[scale=0.65]{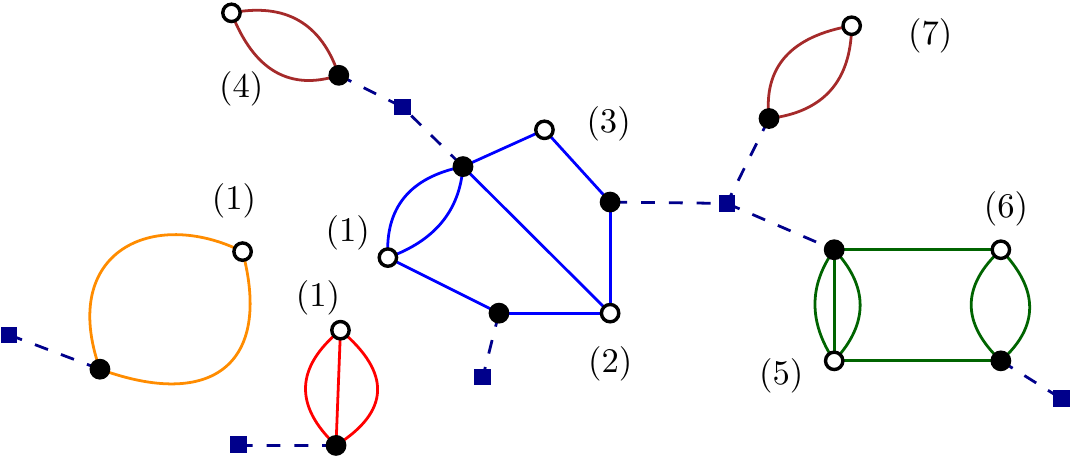}
\caption{Left: the factorized map $(\gamma_\lambda, \tau^{-1})$ and  the blocks of $\pi$ (represented by the white squares and the colors).  Right: the collection of non-separable hypermaps $(\nu, \tau)=\{({\gamma_\lambda}_{\lvert_G}, \nu_G)\}_{G\in \pi}$, together with a forest-like partition $\pi''$, whose blocks are indicated by the blue vertices. 
}
\label{fig:triplet}
\end{figure}

Therefore: 
\be
\beta_{\lambda_1, \ldots, \lambda_p} =
\sum_{\substack{{\tau =\bigtimes_{i} \tau_i}\\[+0.3ex]{\tau_i\in \mathrm{NC}(\lambda_i)}}}\ \sum_{\substack{{\pi \ge \Pi(\tau)}\\[+0.3ex] {L[\pi, \Pi(\gamma_{\lambda})\, ; \,\Pi(\tau)] =0
}}} 
\  \sum_{\substack{{\nu =\bigtimes_{G\in \pi}\nu_G}\\[+0.3ex]{\nu_G\in \mathrm{NS} (\tau_{\lvert_G})}}}
\ \sum_{\substack{{\pi' \ge \pi}\\[+0.3ex]{\pi'\vee\Pi(\gamma_{\lambda})= 1_n}\\[+0.3ex] {L[\pi',\, \Pi(\gamma_{\lambda})\, ; \,\Pi(\tau)] =0
}}}    
\ \,
\sum_{\substack{{\pi''\in \mathcal P(n)}\\{\pi' \ge \pi''  \ge \Pi(\nu)}\\[+0.3ex]{\forall  B\in \pi', \pi''_{\lvert_B} \vee\Pi(\tau_{\lvert_B})= 1_{\lvert B \rvert}}\\[+0.3ex] {L[\pi'',\, \pi\, ; \,\Pi(\nu)] =0
}}}   f(\pi'', \nu). 
\ee
The partition $\pi'$ has become superfluous as it is just $\pi'' \vee \Pi(\tau)$. It corresponds to the connected components of the graph obtained from the collection of non-separable hypermaps $(\tau, \nu)$ by drawing the blocks of the partition $\pi''$ as on the right of Fig.~\ref{fig:triplet}.

Since $\Pi(\gamma_\lambda)\ge \Pi(\tau)$, the condition $\pi'\vee\Pi(\gamma_{\lambda})= 1_n$ is translated as $\pi'' \vee \Pi(\gamma_\lambda)= 1_{n}$. Conditionally to $L[\pi, \Pi(\gamma_{\lambda})\, ; \,\Pi(\tau)] =0$, the following is true (where we recall that $\Pi(\tau)\vee \Pi(\nu)$):
\be
 \left\{
    \begin{array}{cc}
        L[\pi'' \vee \Pi(\tau),\, \Pi(\gamma_{\lambda})\, ; \,\Pi(\tau)] =0\\
	 L[\pi'',\, \pi\, ; \,\Pi(\nu)]  =0
    \end{array}
    \right\}\qquad \Leftrightarrow  \qquad L[\pi'',\, \Pi(\gamma_{\lambda})\vee \Pi(\nu)\, ; \,\Pi(\nu)] =0,
\ee
as: 
$$
L[\pi'' \vee \Pi(\tau),\, \Pi(\gamma_{\lambda})\, ; \,\Pi(\tau)]  + L[\pi'',\, \pi\, ; \,\Pi(\nu)] = L[\pi'',\, \Pi(\gamma_{\lambda})\vee \Pi(\nu)\, ; \,\Pi(\nu)] + L[\pi, \Pi(\gamma_{\lambda})\, ; \,\Pi(\tau)]. 
$$
Therefore:
\be
\label{eq:eq-in-proof-facto-1}
\beta_{\lambda_1, \ldots, \lambda_p} =
\sum_{\substack{{\tau =\bigtimes_{i} \tau_i}\\[+0.3ex]{\tau_i\in \mathrm{NC}(\lambda_i)}}}\ \sum_{\substack{{\pi \ge \Pi(\tau)}\\[+0.3ex] {L[\pi, \Pi(\gamma_{\lambda})\, ; \,\Pi(\tau)] =0
}}} 
\  \sum_{\substack{{\nu =\bigtimes_{G\in \pi}\nu_G}\\[+0.3ex]{\nu_G\in \mathrm{NS} (\tau_{\lvert_G})}}}
\ 
\sum_{\substack{{\pi''  \ge \Pi(\nu)}\\[+0.3ex]{\pi'' \vee \Pi(\gamma_\lambda)= 1_{n}}\\[+0.3ex] {L[\pi'',\, \Pi(\gamma_{\lambda})\vee \Pi(\nu)\, ; \,\Pi(\nu)] =0
}}}   f(\pi'', \nu). 
\ee

Consider a triplet $(\tau, \pi, \nu)$, where $\tau$ is  a factorized planar bipartite map $\tau =\bigtimes_{i} \tau_i$, $\tau_i \in  \mathrm{NC}(\lambda_i)$;  $\pi \ge \Pi(\tau)$ satisfies $L[\pi,\, \Pi(\gamma_{\lambda})\, ; \,\Pi(\tau)] =0$, that is, it is a forest-like partition of the black vertices of the factorized map $(\gamma_\lambda, \tau^{-1})$; and  $\nu =\bigtimes_{G\in \pi}\nu_G$, where the $\nu_G\in \mathrm{NS} (\tau_{\lvert_G})$ are  non-separable hypermaps, one for each block of $\pi$, so that in particular, $\pi=\Pi(\nu)\vee\Pi(\tau)$.

\begin{figure}[!h]
\centering
\includegraphics[scale=0.75]{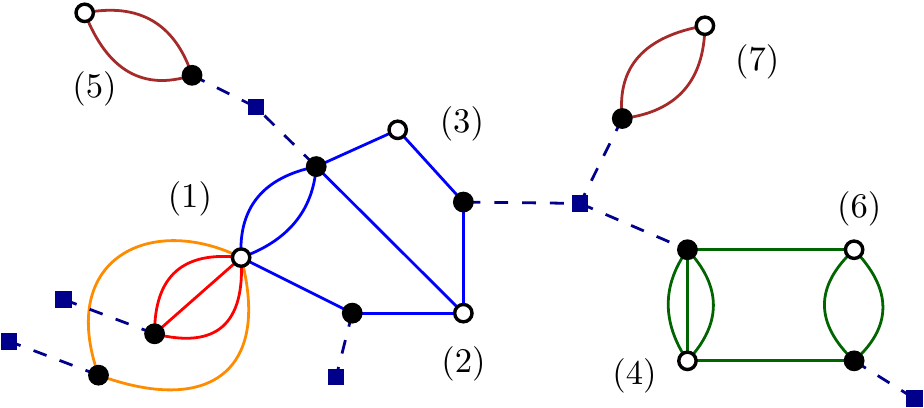}\hspace{1cm}
\caption{The generic planar map $(\gamma_\lambda, \nu)$, together with the forest-like partition $\pi''$. The colors correspond to its decomposition in non-separable hypermaps. 
}
\label{fig:gen-map}
\end{figure}

The permutation $\nu$ defines a non-necessarily-factorized planar bipartite map with $p$ white vertices $(\gamma_\lambda, \nu)$. The partition $\pi\in \mathcal{P}(n)$ can then be read on $(\gamma_\lambda, \nu)$: it is the partition of its edges given by the unique decomposition in non-separable hypermaps (see Sec.~\ref{sub:bip-map}).  This map can be seen as obtained gluing back the non-separable hypermaps on the right of Fig.~\ref{fig:triplet} to the white vertices of the factorized map on the left of Fig.~\ref{fig:triplet} (since each edge is an edge both around a white vertex of the factorized map and of a non-separable hypermap), obtaining the example of Fig.~\ref{fig:gen-map}. 

A consequence of $\pi \ge \Pi(\tau)$,  $L[\pi,\, \Pi(\gamma_{\lambda})\, ; \,\Pi(\tau)] =0$ is that the restriction of $(\gamma_\lambda, \tau^{-1})$ to a block $G\in \pi$ is a collection of disjoint planar bipartite maps with two vertices (bottom left of Fig.~\ref{fig:first-expl}), that is,  $\tau_G={\gamma_\lambda}_{\lvert_G}$. This means that the decomposition of $(\gamma_\lambda, \nu)$ in non-separable hypermaps  is $\{({\gamma_\lambda}_{\lvert_G}, \nu_G)\}_{G\in \pi} = (\tau, \nu)$. Consider a recursive splitting of edges of  $(\gamma_\lambda, \nu)$ leading to its  terminal form $(\tau, \nu)$. At every splitting, the two corners belong to the same face and a new face is formed after the splitting. The total number of faces created is  the difference between the number of connected components of $(\tau, \nu)$ and $(\gamma_\lambda, \nu)$, so that:
\be
\label{eq:faces-when-glued}
\#( \nu\gamma_\lambda) = \#( \nu\tau) + \#(\Pi(\nu)\vee \Pi(\gamma_\lambda)) - \#(\Pi(\nu)\vee \Pi(\tau)).
\ee
From this, we see that 
\begin{align*}
2g (\gamma_\lambda, \nu) &= 2\#(\Pi(\nu)\vee \Pi(\gamma_\lambda)) + n - \# (\nu) - \# (\gamma_\lambda) -  \#(\nu\gamma_\lambda) 
\\&= 2g(\nu, \tau) + L[\Pi(\nu)\vee\Pi(\tau),\, \Pi(\gamma_{\lambda})\, ; \,\Pi(\tau)]  
\end{align*}
so that $g (\gamma_\lambda, \nu)$ vanishes as long as  $L[\pi,\, \Pi(\gamma_{\lambda})\, ; \,\Pi(\tau)] =0$ and  $g(\nu, \tau)=0$.

Reciprocally, given $\nu\in S_n$ which defines a generic planar bipartite map $g(\gamma_\lambda, \nu)=0$, we decompose it into a triplet $(\tau, \pi, \nu)$ as follows: We consider its unique decomposition in non-separable hypermaps. It induces a partition $\pi\in \mathcal{P}(n)$ of the edges of the original map, which is such that or each block $G$  of this partition, the permutations ${\gamma_\lambda}_{\lvert_G}$ and $\nu_G$ respectively induced by ${\gamma_\lambda}$ and $\nu$ on $G$ satisfy  $\nu_G\in \mathrm{NS} ({\gamma_\lambda}_{\lvert_G})$. Letting $\gamma_i$ be the restriction of  $\gamma_\lambda$ to $\{\sum_{j<i} \lambda_j + 1, \ldots, \sum_{j<i} \lambda_j + \lambda_i \}$, we  now set  $\tau_i:= \bigtimes_{G\in \pi}{ \gamma_i }_{\lvert_G}$. From the construction, $\pi \ge \Pi(\tau)$. 

We have to justify that $L[\pi,\, \Pi(\gamma_{\lambda})\, ; \,\Pi(\tau)] =0$: it is equivalent to the following incidence graph $\mathcal G$ being a forest: consider the white vertices of the map (the cycles of $\gamma_\lambda$), add a new black vertex for each non-separable hypermap  component $G\in \pi$, and an edge between a black vertex and a white vertex whenever some of the edges around that white vertex belong to $G$ in  $(\gamma_\lambda, \nu^{-1})$. Consider a recursive splitting of edges of the map $(\gamma_\lambda, \nu)$ leading to its decomposition in non-separable hypermaps. Each such step corresponds to the removal of an edge of $\mathcal G$, and the number of connected components are raised both in the map and the incidence graph. The terminal map whose connected components are the non-separable hypermap components of $(\gamma_\lambda, \nu)$ corresponds to the terminal incidence graph which is just a collection of vertices with no edges. This shows that $\mathcal G$ is a forest.

Therefore, \eqref{eq:eq-in-proof-facto-1} is equivalent to
$$
\beta_{\lambda_1, \ldots, \lambda_p} =
\sum_{\substack{{\nu\in S_n,}\\[+0.3ex]{g(\gamma_\lambda, \nu)=0}}}  
\ 
\sum_{\substack{{\pi''  \ge \Pi(\nu)}\\[+0.3ex]{\pi'' \vee \Pi(\gamma_\lambda)= 1_{n}}\\[+0.3ex] {L[\pi'',\, \Pi(\gamma_{\lambda})\vee \Pi(\nu)\, ; \,\Pi(\nu)] =0
}}}   f(\pi'', \nu). 
$$
This indeed coincides with \eqref{eq:assumption-facto} by opting for the inverse of $\nu$. \qed

\subsection{Corrections to the vertex weights for small orders}
\label{sub:small-orders}

In this subsection, we take a closer look at the weight $g$  \eqref{eq:split-perm-2}  associated to the blocks of the tree-like partition $\pi'$ in Th.~\ref{thm:facto}.

\paragraph{First order.}Consider \eqref{eq:split-perm-2}. If $\tau_{\lvert_B}$ has one cycle of length $d$ only, then it is a cycle of one of the $\tau_i$, so that $I_B=\{i\}$, and $\bar \pi$ is the trivial partition of a single element, which is not an element of $\mathcal{P}^*(\{i\})$. In this case, there is no correction to $f$:
$$
g (1_{d}, \gamma_{d}) = f(1_{d}, \gamma_{d}). 
$$

\subsubsection{Second order}
If $\tau_{\lvert_B}$ has only two cycles respectively of lengths  $m, n$, \eqref{eq:split-perm-2}  is the trick of \cite{CMSS}. There are two contributions, and we can write \eqref{eq:split-perm-2} more simply as:
\be
g(1_{m+n}, \gamma_{m,n}) =  f(1_{m+n}, \gamma_{m,n}) +\sum_{\nu\in \mathrm{NS}(m,n)} f(\Pi(\nu), \nu).
\ee
Indeed,  the two cycles of  $\tau_{\lvert_B}$ must be respectively cycles of $\tau_a$ and $\tau_b$, $a < b$, and 
the sum over $\bar \pi$ in  \eqref{eq:split-perm-2} has one term: the one-block partition $\bar \pi = \{a,b\}$, which is given by \eqref{eq:bar-f-one-block}.
Setting   $f(1_{k_1 + \cdots + k_p}, \gamma_{k_1, \ldots, k_p})= \kappa_{k_1, \ldots, k_p} $, we have in terms of generating functions:
\be
\label{eq:H2}
H_2(Y_1, Y_2) = C_2 (Y_1, Y_2) + \tilde C_2 (Y_1, Y_2).
\ee 
Letting
\be
\label{eq:def-tilde-kappa}
\tilde \kappa_{m, n}  : =  \sum_{\nu\in \mathrm{NS} (m, n)} \kappa(\Pi(\nu), \nu) ,
\ee 
we see from \eqref{eq:combinatorial-expr-correction-with-trees} (see also the proof of Prop.~6.1 of \cite{CMSS}), that
\be
\label{eq:tilde-kappa}
\tilde \kappa_{m,n}  =  \sum_{r\ge 1} \sum_{\substack{{i_1, \ldots, i_r\ge 1}\\{i_1+\cdots + i_r = m}}}\,\,\sum_{\substack{{j_1, \ldots, j_r\ge 1}\\{j_1+\cdots + j_r = n}}}i_1 n \ \prod_{s=1}^r 
\kappa_{i_s + j_s} . 
\ee
Indeed, the elements of $\mathrm{NS} (m, n)$ consist of two white vertices  $a< b$ and $r$ black vertices, $r\ge 1$, such that each black vertex is connected to both white vertices. Labelling the black edges as  above \eqref{eq:combinatorial-expr-correction-with-trees}, we consider the simpler notations $i_s := j_{a, s}$ and $j_s =  j_{b, s}$. The tree $T$ being trivial, we only have the bottom right contribution of \eqref{eq:combinatorial-expr-correction-with-trees} with $j_G^s= i_s + j_s$. 
The maps that have the same permutation for the white vertices and the same adjacency matrix differ only by the position of the edge with smallest label among the $i_1$ edges connecting $a$ and $1$, and the $n=\sum_s j_s$ ways to rotate the cycle $b$, so that  $\mathcal{N}_{r}(\{i_s,j_{s}\}_{s})=i_1 n$.

This leads to the generating function \eqref{eq:bar-C-one-block}:
$$
\tilde C_{2} \bigl(Y_1, Y_2\bigr)  = \sum_{r\ge 1 } Y_2   \frac{\partial}{\partial Y_2}  \Biggl\{ \hat C(Y_1, Y_2)^{r-1}\,  Y_1 \frac{\partial}{\partial Y_1} \hat C(Y_1, Y_2) \Biggr\},
$$
with 
$\hat C(Y_1, Y_2) = \sum_{i, j\ge 1} \kappa_{i+j} Y_1^i Y_2^j$, which are computed as:
\be
\label{eq:C2Tilde-log}
\tilde C_{2} \bigl(Y_1, Y_2\bigr) = -Y_1 Y_2 \frac{\partial}{\partial Y_1} \frac{\partial}{\partial Y_2} \log\bigl(1 - \hat C(Y_1, Y_2) \bigr), 
\ee
and:
\be
\label{eq:hat-C}
\hat C(Y_1, Y_2) = 1 + \frac{Y_1}{Y_2 - Y_1} C_1(Y_2) + \frac{ Y_2} {Y_1 - Y_2} C_1 (Y_1).
\ee
For $Y_i = X_i C_1(Y_i) = X_i M_1(X_i) $, this  can be computed as:
\be
\label{eq:Tilde-C2-in-pat-ex}
\tilde C_{2} \bigl(Y_1, Y_2\bigr) = \frac{Y_1Y_2}{(Y_1 - Y_2)^2}- \frac{Y_1 Y_2}{(X_1 - X_2)^2} \frac {\mathrm{d} X_1} {\mathrm d Y_1} \frac {\mathrm{d} X_2} {\mathrm d Y_2} , 
\ee 
and we recover the second-order relation \eqref{eq:second-order-functional-rel}.

\subsubsection{Third order}
\begin{figure}[!h]
\centering
\includegraphics[scale=0.7]{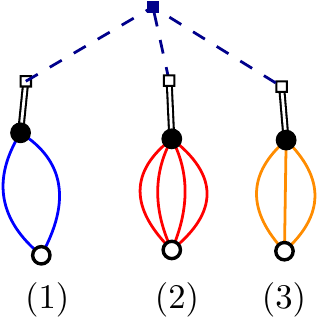}\hspace{2cm}\includegraphics[scale=0.7]{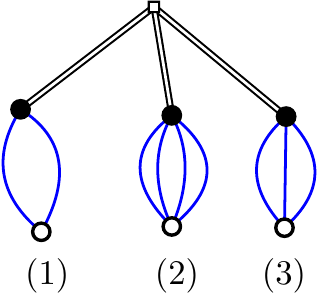}\hspace{2cm}\includegraphics[scale=0.7]{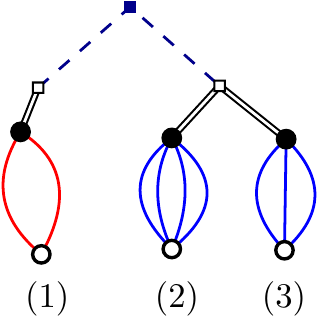}
\caption{Contributions at third order. From left to right:  $C_3$, $\tilde C_3$, $\bar C_{3, \{1\}\{2,3\}}$. 
}
\label{fig:Third-order}
\end{figure}

Assuming now that $\tau_{\lvert_B}$ has only three disjoint cycles $\eta_a$, $\eta_b$, and $\eta_c$, respectively cycles of $\tau_a$, $\tau_b$, and $\tau_c$  for $a<b<c$, and respectively of lengths $l, m, n$,  there are now three kinds of contributions to \eqref{eq:split-perm-2}, illustrated in Fig.~\ref{fig:Third-order}: that given by $f$, that given by $\tilde f$ \eqref{eq:bar-f-one-block} for $\bar \pi$ the one-block partition and $T$ the trivial tree, and the terms of the form \eqref{eq:bar-f-one-tree} for $\bar \pi$ having two blocks and $T$ being the only element of $\mathcal{G}_2$. If $\bar \pi = \{a\}\{b, c\}$:
\be
\label{eq:first-expr-pibardeuxterms}
\bar f _{\{a\}\{b,c\}}  (\tau_{\lvert_B}) = \sum_{\substack{{\nu =\eta_i \times \nu_2}\\[+0.4ex]{\nu_2 \in \mathrm{NS} (\eta_j \times \eta_k )}}} \ 
\sum_{F_2\in \Pi(\nu_2)}  \  f\Bigl(1_{l + \lvert F_2 \rvert}, \eta_i \times ({\nu_2}_{\lvert _{F_2}})\Bigr) \prod_{\substack{{V \in \Pi(\nu_2)}\\{V \notin F_2} }}  f(1_{\lvert V \rvert}, {\nu_2}_{\lvert _{V }}). 
\ee
Setting   $f(1_{k_1 + \cdots + k_p}, \gamma_{k_1, \ldots, k_p})= \kappa_{k_1, \ldots, k_p} $, we have in terms of generating functions:
\be
\label{eq:H3}
H_3 = C_3  + \tilde C_3 + \bar C_{3, \{1\}\{2,3\}}  + \bar C_{3, \{1,2\}\{3\}}  + \bar C_{3, \{2\}\{1,3\}}  .
\ee 

\begin{proposition} 
\label{prop:gen-fun-third-order}
For the third order case, the corrections to $C_3$ in $H_3$ are the three terms given by \eqref{eq:bar-C} for $\#(\bar \pi)=2$, of the form
\be
\label{eq:gen-func-pibardeuxterms}
\bar C_{3, \{1\}\{2,3\}} \bigl(Y_1, Y_2, Y_3\bigr) =  Y_2 Y_3  \frac{\partial}{\partial Y_2} \frac{\partial}{\partial Y_3} \left[ \frac{Y_3 C_2(Y_1, Y_2) - Y_2 C_2(Y_1, Y_3)}{  Y_2 C_1(Y_3) - Y_3 C_1(Y_2)}\right], 
\ee
and similarly for $\bar \pi = \{2\}\{1,3\}$ and $\bar \pi=\{3\}\{1,2\}$, and the remaining term  $\tilde C_3$ \eqref{eq:bar-C-one-block} given by
 \begin{align}
\nonumber
\tilde C_3(Y_1, Y_2, Y_3) &= Y_1 Y_2 Y_3 \frac{\partial}{\partial Y_1}\frac{\partial}{\partial Y_2} \frac{\partial}{\partial Y_3}  \Biggl[ \frac{(1 + \hat C(Y_1, Y_2, Y_3))^2 + \hat C(Y_1, Y_2)\hat C(Y_1, Y_3)\hat C(Y_2, Y_3) -1}{(1- \hat C(Y_1, Y_2))(1- \hat C(Y_1, Y_3))(1- \hat C(Y_2, Y_3))}\\&\hspace{5cm}+ \sum_{\substack{{1\leftrightarrow 2}\\{1\leftrightarrow 3}}}   \frac{Y_1^2\frac{\partial}{\partial Y_1}\left(\frac 1 {Y_1}\hat C(Y_1, Y_2, Y_3) \right) } {(1 - \hat C(Y_1, Y_2))(1 - \hat C(Y_1, Y_3))}  \Biggr], 
\label{eq:tilde-c3-gen}
\end{align}
where $\hat C$ was defined in \eqref{eq:hat-C}, and:
\be
\hat C(Y_1, Y_2, Y_3) =   \frac{Y_3} {Y_2 - Y_3} \hat C(Y_1, Y_2)  +   \frac{Y_2} {Y_3 - Y_2} \hat C (Y_1, Y_3).
\ee
With this and from Thm.~\ref{thm:functional-relations-for-any-H} and Thm.~\ref{thm:facto}, we have all the ingredients of the third-order moment-cumulant relation. 
\end{proposition}

\proof From \eqref{eq:combinatorial-expr-correction-with-trees},
the expression for $\bar \kappa_{\{1\}\{2,3\}}  (\gamma_{l,m,n}) $ is very similar to that of $\tilde \kappa_{m,n}  $ in \eqref{eq:tilde-kappa}, with the difference that a cycle $\theta \in \{1, \ldots, r\}$ of $\nu_2\in \mathrm{NS} (m,n)$ is chosen, for which the term  $\kappa_{i_{\theta} + j_{\theta}}$  is replaced by the second order free cumulant $\kappa_{l, i_{\theta} + j_{\theta}}$:
\be
\label{eq:gen-func-pibar-unbloc}
\bar \kappa _{\{1\}\{2,3\}}  (\gamma_{l,m,n}) =  \sum_{r\ge 1} \sum_{\substack{{i_1, \ldots, i_r\ge 1}\\{i_1+\cdots + i_r = m}}}\,\,\sum_{\substack{{j_1, \ldots, j_r\ge 1}\\{j_1+\cdots + j_r = n}}}\sum_{\theta=1}^r i_1 n\,  
\kappa_{i_1 + j_1} \ldots \kappa_{l, i_{\theta} + j_{\theta}}  \ldots \kappa_{i_r + j_r} . 
\ee
The generating function for these numbers is found to be
\begin{align}
\bar C_{3, \{1\}\{2,3\}} \bigl(Y_1, Y_2, Y_3\bigr) &= \sum_{l, m, n \ge 1}\bar \kappa _{\{1\}\{2,3\}}  (\gamma_{l,m,n}) Y_1^l Y_2^m Y_3^n\\&= \sum_{r\ge 1 } Y_2 Y_3  \frac{\partial}{\partial Y_2} \frac{\partial}{\partial Y_3} \Biggl\{ \hat C(Y_2, Y_3)^{r-1}  \sum_{l, i, j\ge 1} \kappa_{l, i+j}  Y_1^l Y_2^i Y_3^j  \Biggr\} \\
&= Y_2 Y_3  \frac{\partial}{\partial Y_2} \frac{\partial}{\partial Y_3} \Biggl\{ \frac{1}{1- \hat C(Y_2, Y_3)} \hat C_{\{1\}\{2,3\}}(Y_1, Y_2, Y_3) \Biggr\},
\end{align}
where $ \hat C_{\{1\}\{2,3\}}(Y_1, Y_2, Y_3)=\sum_{l, i, j\ge 1} \kappa_{l, i+j}  Y_1^l Y_2^i Y_3^j$ is computed as:
\be
\label{eq:hat-C-12}
\bigl( Y_2  - Y_3 \bigr) \hat C_{\{1\}\{2,3\}}(Y_1, Y_2, Y_3)  = Y_3 C_2(Y_1, Y_2) - Y_2 C_2(Y_1, Y_3), 
\ee
which concludes the proof of \eqref{eq:gen-func-pibardeuxterms}. 

\begin{figure}[!h]
\centering
\includegraphics[scale=0.6]{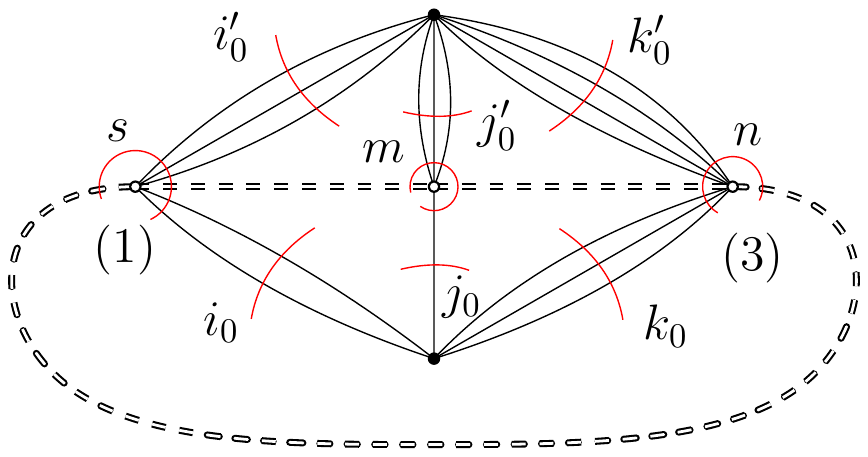}\hspace{1cm}\includegraphics[scale=0.6]{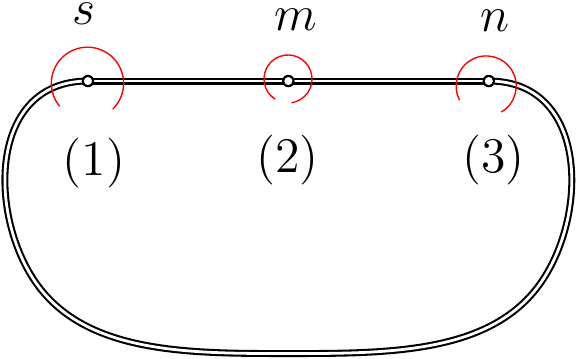}\hspace{1cm}\includegraphics[scale=0.6]{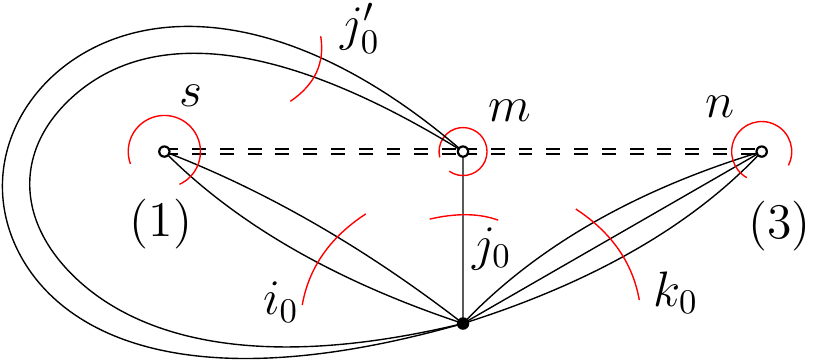}
\caption{Non-separable hypermaps with three white vertices. The thick dotted edges represent non-separable hypermaps with two white vertices. In the case where these thick edges are dotted, these hypermaps may be absent. The black vertices (and incident edges) on the left diagram may be absent or not, but there must be at least one.
}
\label{fig:proof-C3tilde}
\end{figure}

We now turn to $\tilde C_3$ \eqref{eq:bar-C-one-block}, the generating function of non-separable hypermaps with white vertices given by $\gamma_{l, m, n}$, whose hyper-edges are weighted by first-order free cumulants: 
\be
\tilde C_{3} \bigl(Y_1, Y_2, Y_3\bigr) = \sum_{s, m, n}  Y_1^s Y_2^m Y_3^n \sum_{\nu\in \mathrm{NS} (s, m, n)} \ 
 \kappa(\Pi(\nu), \nu). 
 \ee
 There are three kinds of contributions, illustrated in Fig.~\ref{fig:proof-C3tilde}. 

The case where both back vertices are present contributes as:
\begin{align} 
s m n \sum_{\substack{
{i_0, j_0, k_0\ge 1}\\{i_0', j_0', k_0'\ge 1}
}} 
\ 
\sum_{r_{12}, r_{13}, r_{23}\ge 0} 
\ 
\sum_{\substack{{\{i_a\}_{1\le a \le r_{12}}, \{i_a'\}_{1\le a \le r_{13}}}\\{\sum i_a + \sum i_a'+ i_0 + i_0' = s}}}
\sum_{\substack{{\{j_a\}_{1\le a \le r_{12}}, \{j_a'\}_{1\le a \le r_{23}}}\\{\sum j_a + \sum j_a'+ j_0 + j_0' = m}}}
\sum_{\substack{{\{k_a\}_{1\le a \le r_{23}}, \{k_a'\}_{1\le a \le r_{13}}}\\{\sum k_a + \sum k_a'+ k_0 + k_0' = n}}}\nonumber\\
\prod_{a=1}^{r_{12}} \kappa_{i_a+j_a} \prod_{a=1}^{r_{13}} \kappa_{i_a'+k_a'}  \prod_{a=1}^{r_{23}} \kappa_{j_a'+k_a}  \kappa_{i_0 + j_0 + k_0}  \kappa_{i'_0 + j'_0 + k'_0} , 
\end{align}
where the factor $s m n$ comes from rotating the cycles of  $\gamma_{s, m, n}$ for each white vertex (choosing which edge is its smallest element). This leads to the generating function:
$$
Y_1 Y_2 Y_3 \frac{\partial}{\partial Y_1}\frac{\partial}{\partial Y_2} \frac{\partial}{\partial Y_3}  \Biggl[ \frac{\hat C(Y_1, Y_2, Y_3)^2 }{(1- \hat C(Y_1, Y_2))(1- \hat C(Y_1, Y_3))(1- \hat C(Y_2, Y_3))}  \Biggr], 
$$
where $ \hat C(Y_1, Y_2, Y_3)=\sum_{i, j, k\ge 1} \kappa_{i+j+k}  Y_1^l Y_2^i Y_3^j$ is computed as:
\be
\label{eq:hat-C3}
\bigl( Y_2  - Y_3 \bigr) \hat C (Y_1, Y_2, Y_3)  = Y_3 \hat C(Y_1, Y_2) - Y_2 \hat C(Y_1, Y_3). 
\ee 
The contributions when a single black vertex is present is given by 
$$
Y_1 Y_2 Y_3 \frac{\partial}{\partial Y_1}\frac{\partial}{\partial Y_2} \frac{\partial}{\partial Y_3}  \Biggl[ \frac{2 \hat C(Y_1, Y_2, Y_3) }{(1- \hat C(Y_1, Y_2))(1- \hat C(Y_1, Y_3))(1- \hat C(Y_2, Y_3))}  \Biggr], 
$$
and the contribution to the diagram in the middle of Fig.~\ref{fig:proof-C3tilde} by
$$
Y_1 Y_2 Y_3 \frac{\partial}{\partial Y_1}\frac{\partial}{\partial Y_2} \frac{\partial}{\partial Y_3}  \Biggl[ \frac{C(Y_1, Y_2)}{(1- \hat C(Y_1, Y_2))} \frac{C(Y_1, Y_3)}{(1- \hat C(Y_1, Y_3))} \frac{C(Y_2, Y_3)}{(1- \hat C(Y_2, Y_3))}  \Biggr]. 
$$
The right diagram contributes as
\begin{align} 
&s m n \sum_{\substack{
{i_0, j_0, i_0'\ge 1}
}} 
\ 
\sum_{r_{12},  r_{23}\ge 0} 
\ 
\sum_{\substack{{\{i_a\}_{1\le a \le r_{12}}}\\{\sum i_a + i_0 = s}}}
\sum_{\substack{{\{j_a\}_{1\le a \le r_{12}}, \{j_a'\}_{1\le a \le r_{23}}}\\{\sum j_a + \sum j_a'+ j_0 + j_0' = m}}}
\sum_{\substack{{\{k_a\}_{1\le a \le r_{23}}}\\{\sum k_a + k_0  = n}}}\nonumber\\&
\hspace{7cm}\prod_{a=1}^{r_{12}} \kappa_{i_a+j_a} \prod_{a=1}^{r_{23}} \kappa_{j_a'+k_a}  \kappa_{i_0 + j_0 + j_0' +  k_0}  , 
\end{align}
 leading to the generating function
 $$
Y_1 Y_2 Y_3 \frac{\partial}{\partial Y_1}\frac{\partial}{\partial Y_2} \frac{\partial}{\partial Y_3}  \Biggl[ \frac{\sum_{i, j,j',  k\ge 1} \kappa_{i + j + j' +  k} Y_1^i Y_2^{j+j'} Y_3^k }{(1- \hat C(Y_1, Y_2))(1- \hat C(Y_2, Y_3))}  \Biggr], 
 $$
 and we compute:
 $$
\sum_{i, j,j',  k\ge 1}\hspace{-0.2cm} \kappa_{i + j + j' +  k} Y_1^i Y_2^{j+j'} Y_3^k  =\hspace{-0.3cm} \sum_{i, j, k\ge 1} (p-1)  \kappa_{i + j  +  k} Y_1^i Y_2^{j} Y_3^k = Y_2\frac{\partial}{\partial Y_2} \hat C(Y_1, Y_2, Y_3) - \hat C(Y_1, Y_2, Y_3).
 $$
 This concludes the proof of \eqref{eq:tilde-c3-gen}.
  \qed

 \paragraph{Remark.} One may alternatively compute $\tilde C_3$ as
 \begin{align}
\label{eq:tilde-c3-alt}
\tilde C_3(Y_1, Y_2, Y_3) &= Y_1 Y_3 \frac{\partial}{\partial Y_1} \frac{\partial}{\partial Y_3} \left(\frac {Y_3 \tilde C_2(Y_1, Y_2) - Y_1 \tilde C_2 (Y_2, Y_3)}{Y_1 C_1(Y_3) - Y_3 C_1(Y_1)}\right) \\&+ Y_1 Y_2 Y_3 \frac{\partial}{\partial Y_1}\frac{\partial}{\partial Y_2} \frac{\partial}{\partial Y_3} \left[ \frac{Y_2^2\frac{\partial}{\partial Y_2}\left(\frac 1 {Y_2}\hat C(Y_1, Y_2, Y_3) \right)  - \hat C(Y_1, Y_2)\hat C(Y_2, Y_3)} {(1 - \hat C(Y_1, Y_2))(1 - \hat C(Y_2, Y_3))}   \right]\nonumber, 
\end{align}
 In a similar way as above, recalling the definition of $\tilde \kappa_{i,j}$ in \eqref{eq:def-tilde-kappa}, one finds that the following series 
 \be
 \label{eq:firs-contrib-tilde-c3}
  \sum_{r\ge 1} \sum_{\substack{{i_1, \ldots, i_r\ge 1}\\{i_1+\cdots + i_r = s}}}\,\,\sum_{\substack{{j_1, \ldots, j_r\ge 1}\\{j_1+\cdots + j_r = n}}}\sum_{\theta=1}^r i_1 n\,  
\kappa_{i_1 + j_1} \ldots \tilde \kappa_{m, i_{\theta} + j_{\theta}}  \ldots \kappa_{i_r + j_r} ,
 \ee
generates exactly once the major part of the elements of $\mathrm{NS} (s, m, n)$, with the exception of the terms for the diagram on the right of Fig.~\ref{fig:proof-C3tilde}, which we have to add, and  which we have computed above. It also generates elements as in the middle of Fig.~\ref{fig:proof-C3tilde} but without the thick line linking 1 and 3, which  are not in $\mathrm{NS} (s, m, n)$, and which we have to substract.
The series \eqref{eq:firs-contrib-tilde-c3} is similar to what we have detailed above and provides the contribution in the first line of \eqref{eq:tilde-c3-gen}. The terms to be substracted are generated by:
\be
Y_1Y_2Y_3 \frac{\partial}{\partial Y_1}\frac{\partial}{\partial Y_2} \frac{\partial}{\partial Y_3}  \left(\frac{\hat C(Y_1, Y_2)} {(1 - \hat C(Y_1, Y_2))}\frac{\hat C(Y_2, Y_3)}{(1 - \hat C(Y_2, Y_3))}  \right). 
\ee

\subsubsection{Fourth order}

\begin{figure}[!h]
\centering
\includegraphics[scale=0.7]{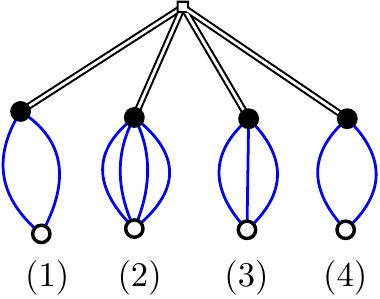}\hspace{0.6cm}\includegraphics[scale=0.7]{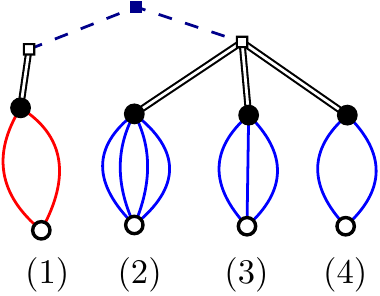}\hspace{0.6cm}\includegraphics[scale=0.7]{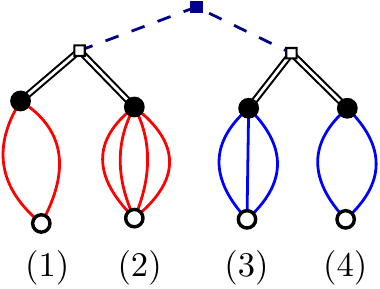}\hspace{0.6cm}\includegraphics[scale=0.7]{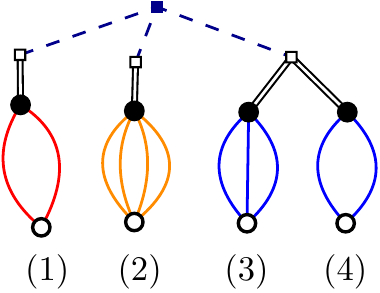}\hspace{0.6cm}\includegraphics[scale=0.7]{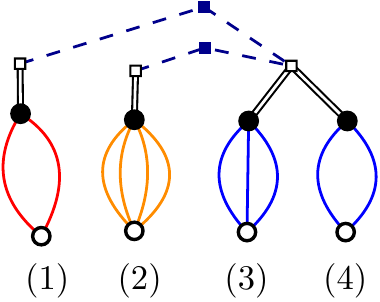}
\caption{Corrections at fourth order. From left to right:  $\tilde C_4$, $\bar C_{4, \{1\}\{2,3, 4\}}$,  $\bar C_{4, \{1,2\}\{3, 4\}}$, $\bar C_{4, \{1\}\{2\}\{3, 4\}}$,  and $\bar C_{4, \{1\}\{2\}\{3, 4\}, T}$ for $T$ with two blue vertices. 
}
\label{fig:Fourth-order}
\end{figure}

 Consider $\tau_{\lvert_B}$ with four disjoint cycles $\eta_a$, $\eta_b$, $\eta_c$ and $\eta_d$, respectively cycles of $\tau_a$, $\tau_b$, $\tau_c$ $\tau_d$  for $a<b<c<d$, and respectively of lengths $l, m, n, q$,  there are now six kinds of contributions to \eqref{eq:split-perm-2}, illustrated in Fig.~\ref{fig:Fourth-order}:
that given by $f$; that given by $\tilde f$ \eqref{eq:bar-f-one-block} for $\bar \pi$ the one-block partition and $T$ the trivial tree; two kinds of  terms of the form \eqref{eq:bar-f-one-tree} for $\bar \pi$ having two blocks and $T$ being the only element of $\mathcal{G}_2$:  one for $\bar \pi$ having blocks with 1 and 3 elements, the other with 2 and 2 element;   and finally  two kinds of terms for $\bar \pi$ with three blocks: one of the form \eqref{eq:bar-f-one-tree} for $T$ with a single blue vertex, and the one where $T$ has two blue vertices. 
\be
\label{eq:H4}
H_4 = C_4  + \tilde C_4 + \sum_{\mathrm{perm}} \bar C_{4, \{1\}\{2,3, 4\}}   +   \sum_{\mathrm{perm}} \bar C_{4, \{1,2\}\{3, 4\}} +   \sum_{\mathrm{perm}} \bar C_{4, \{1\}\{2\}\{3, 4\}}  +  \sum_{\mathrm{perm}} \bar C_{4, \{1\}\{2\}\{3, 4\}, T}, 
\ee 
where the sums are over appropriate permutations of the indices. 

Here we only compute $C_{4, \{1,2\}\{3,4\}}$ and  $\bar C_{4, \{1\}\{2\}\{3, 4\}, T}$. $\bar C_{4, \{1\}\{2\}\{3, 4\}}$ will be computed in \eqref{eq:bar-C-ab} as the particular case of $\bar C_{p, \bar \pi}$ for $\bar \pi$ having only one non-trivial block with two elements, see also Prop.~\ref{prop:order_p-1_at-orderp}. $ \bar C_{4, \{1\}\{2,3, 4\}}$ can be computed from the expression of $\tilde C_3$ above, but we  will instead compute it in  \eqref{eq:hat-C3-13}  after deriving a simpler expression for $\tilde C_3$ in Prop.~\ref{prop:simpl-third-order}, to avoid a three-lines-equation. Computing $\tilde C_4$ as we did for $\tilde C_3$ is painful and  a more efficient approach should be prefered. 
 \begin{proposition} 
\label{prop:Bar-C-4-22}
The following corrections can be easily computed:
\be
\label{eq:expr-C-4-13-tree-prop}
\bar C_{4, \{1\},\{2\}\{3,4\}, T}(Y_1, Y_2, Y_3, Y_4) = Y_3 Y_4 \frac{\partial}{\partial Y_3}\frac{\partial}{\partial Y_4} \left[\frac{\hat C_{\{1\}\{3,4\}}(Y_1,  Y_3, Y_4) \hat C_{\{2\}\{3,4\}}(Y_2,  Y_3, Y_4)}{(1 - \hat C(Y_3, Y_4))^2}\right], 
\ee
where $T$ has two blue vertices and $\hat C_{\{a\}\{b,c\}}$ is given in \eqref{eq:hat-C-12}, and:
\be
\label{eq:eq:expr-C-4-22-prop}
\bar C_{4, \{1,2\}\{3,4\}}(Y_1, Y_2, Y_3, Y_4) = \left(\prod_{i=1}^4 Y_i \frac{\partial}{\partial Y_i} \right)\left[\frac{\hat C_{ \{1,2\}\{3,4\}}(Y_1, Y_2, Y_3, Y_4)}{(1 - \hat C(Y_1, Y_2))(1 - \hat C(Y_3, Y_4)}\right], 
\ee
where
\begin{align}
\label{eq:hat-C-22}
\hat C_{\{1,2\}\{3,4\}}(Y_1, Y_2, Y_3, Y_4) &= \frac{Y_2}{Y_1 - Y_2}\frac{Y_4}{Y_3 - Y_4}C_2(Y_1, Y_3)+\frac{Y_2}{Y_1 - Y_2}\frac{Y_3}{Y_4 - Y_3}C_2(Y_1, Y_4)\\&\hspace{2cm}\nonumber +\frac{Y_1}{Y_2 - Y_1}\frac{Y_4}{Y_3 - Y_4}C_2(Y_2, Y_3)+\frac{Y_1}{Y_2 - Y_1}\frac{Y_3}{Y_4 - Y_3}C_2(Y_2, Y_4).
\end{align}
\end{proposition}
\proof From  \eqref{eq:combinatorial-expr-correction-with-trees}, the expression \eqref{eq:bar-f-gen} for $\bar \kappa _{\{1\}\{2\}\{3,4\}, T}(\gamma_{m_1,n_1,m_2,n_2})$ is very similar to  \eqref{eq:gen-func-pibar-unbloc}, with the difference that two edges of the tree must be attributed to black vertices of the non-separable hypermap  defined by $\nu \in \mathrm{NS} (m_2,n_2)$:
\begin{align*}
\bar \kappa _{\{1,2\}\{3,4\}}  (\gamma_{m_1,n_1,m_2,n_2}) &=   \sum_{r\ge 2} \sum_{\substack{{i_1, \ldots, i_r\ge 1}\\{i_1+\cdots + i_r = m_2}}}\,\,\sum_{\substack{{j_1, \ldots, j_r\ge 1}\\{j_1+\cdots + j_r = n_2}}}\sum_{\theta_1 \neq \theta_2=1}^r i_1 n\,  
\kappa_{i_1 + j_1} \\&\hspace{4cm}\ldots  \kappa_{m_1, i_{\theta_1} + j_{\theta_1}}  \ldots  \kappa_{n_1, i_{\theta_2} + j_{\theta_2}}  \ldots\kappa_{i_r + j_r} ,
\end{align*}
The generating function is found to be:
$$
\bar C_{4, \{1\},\{2\}\{3,4\}, T} = Y_3 Y_4 \frac{\partial}{\partial Y_3}\frac{\partial}{\partial Y_4} \left[\hat C_{\{1\}\{3,4\}}(Y_1,  Y_3, Y_4) \hat C_{\{2\}\{3,4\}}(Y_2,  Y_3, Y_4) \sum_{r\ge 2} (r-1) \hat C(Y_3, Y_4)^{r-2}\right], 
$$
which is summed to \eqref{eq:expr-C-4-13-tree-prop}.

The expression \eqref{eq:bar-f-gen} for $\bar \kappa _{\{1,2\}\{3,4\}}(\gamma_{m_1,n_1,m_2,n_2})$ is as \eqref{eq:gen-func-pibar-unbloc}, but now two summations over $\nu_1 \in \mathrm{NS} (m_1,n_1)$  and $\nu_2 \in \mathrm{NS} (m_2,n_2)$, and two cycles are chosen, a cycle $\theta_1$ of $\nu_1$ and one  $\theta_2$ of $\nu_2$,  for which the product of first order free cumulants   is replaced by a second order free cumulant. This leads to:
\begin{align}
\bar \kappa _{\{1,2\}\{3,4\}}  (\gamma_{m_1,n_1,m_2,n_2}) &=  \prod_{s=1, 2}\Biggl(\sum_{r_s\ge 1} \sum_{\substack{{i^s_1, \ldots, i^s_{r_s}\ge 1}\\{i^s_1+\cdots + i^s_{r_s} = m_s}}}\,\,\sum_{\substack{{j^s_1, \ldots, j^s_{r_s}\ge 1}\\{j^s_1+\cdots + j^s_r = n_s}}}\sum_{\theta_s=1}^{r_s}\Biggr) \\&\hspace{3cm \nonumber } \times  i^1_1  i^2_1 n_1 n_2 \prod_{s=1,2}\prod_{\substack{{1\le a_s \le r_s}\\{a_s \neq \theta_s}}} \kappa_{i^s_a + j_a^s}  \times \kappa_{i^1_{\theta_1} + j^1_{\theta_1}, i^2_{\theta_2} + j^2_{\theta_2}}. 
\end{align}
The generating function for these numbers is therefore found to be \eqref{eq:eq:expr-C-4-22-prop}, where 
\be
\label{eq:C-hat-22}
 \hat C_{\{1,2 \}\{3,4\}}(Y_1, Y_2, Y_3, Y_4)  =  \sum_{i,j,k,l\ge 1} \kappa_{i+j, k+l}  Y_1^i Y_2^j Y_3^k Y_4^l.  
\ee
As for \eqref{eq:hat-C-12},  this is computed as:
\be
\bigl( Y_1  - Y_2 \bigr) \hat C_{\{1,2 \}\{3,4\}}(Y_1, Y_2, Y_3, Y_4) = Y_2 \hat C_{\{1\}\{3,4\}}(Y_1, Y_3, Y_4)  -   Y_1 \hat C_{\{2\}\{3,4\}}(Y_2, Y_3, Y_4), 
\ee
which leads to \eqref{eq:hat-C-22}. \qed 

\newpage
\section{Efficient treatment of the corrections, simplifications}
\label{sec:corrections}

\subsection{Corrections for several-blocks partitions}

For $I\subset \{1, \ldots, p\}$, the corrections to $C_{\lvert I \rvert}$ in $H_{\lvert I \rvert}$ \eqref{eq:def-HI} are of the form $\bar C_{\lvert I \rvert, \bar \pi, T}$ \eqref{eq:bar-C-T} for $ \bar \pi \in \mathcal{P}^\star(I)$ and $T\in \mathcal{G}_{\bar \pi}$. The correction $\bar C_{\lvert I \rvert, \bar \pi, T}$ simplifies to $\bar C_{\lvert I \rvert, \bar \pi}$ \eqref{eq:bar-C} if the tree $T$ has a single blue vertex,  and to 
$\tilde C_{\lvert I \rvert}$ if $\bar \pi$ is the one-block partition $\bar\pi= \{I\}$. We have already computed three of the several-block corrections ($\bar C_{3,\{1\}\{2,3\}}$, $\bar C_{4,\{1,2 \}\{3,4\}}$ and $\bar C_{4, \{1\},\{2\}\{3,4\}, T}$), and  at the level of the coefficients $\bar \kappa_{\bar \pi, T}$ \eqref{eq:combinatorial-expr-correction-with-trees}, the picture is a bit more clear, and it would be good to make these computations more efficient (and less repetitive) by \emph{working directly with the generating functions}. 

The picture of \eqref{eq:combinatorial-expr-correction-with-trees} is the following: to each pair $(G,K)$ where $G,K$ are respectively a white and a blue vertex  of the tree $T$ corresponds an edge of $T$. For each such edge, a choice is made of  one of the first order free cumulants from the expansions of the $\tilde C_{\lvert G \rvert}$ (this choice corresponds to $\theta({G,K})$ in \eqref{eq:combinatorial-expr-correction-with-trees}). For each $K$, the cumulants for the edges $(G,K)$ incident to $K$ are replaced by a higher order cumulant. In this subsection, we explain how these computations can be performed at the level of generating functions instead of coefficients,  thus showing how all the corrections $\bar C_{\lvert I \rvert, \bar \pi, T}$ can be computed efficiently by acting directly on the generating function $\tilde C_{\lvert I \rvert}$ of planar non-separable hyper-maps with weighted hyperedges (the one-block correction).

\subsubsection{Trees with a single blue vertex}

We consider $\mathbb{K}$, the set of formal power series in the \emph{first order} free cumulants  $1\cup \{\prod_{a}\kappa_{i_a}\}_{\{i_a> 0\}}$, 
$$
c_0 + \sum_{q} \sum_{i_1, \ldots, i_q \ge 1} c_{i_1, \ldots, i_q} \, \kappa_{i_1} \cdots \kappa_{i_q} 
$$
where the sum over $q$ has a finite support, whose  coefficients $c_0$ and $c_{i_1, \ldots, i_q}$ are themselves some formal power series in a set of variables $\{Y_i\}$ with complex coefficients (that may include higher order free cumulants). For a typical example in ${\mathbb{K}}_2 $: $\sum_{i_1, i_2 \ge 0} \Bigl(\sum_{j_1, j_2 \ge1} j_1 \kappa_{j_1, j_2 + i_2} Y_1^{i_i + j_1} Y_2^{j_2} \Bigr)\kappa_{i_1}\kappa_{i_2},$ where $\kappa_0=1$. We define an operator by its action on the $\{\kappa_i\}_{i\ge 0}$:
\be
\mathcal{D} : (\kappa_{i_1}, \ldots, \kappa_{i_d}) \mapsto \left\{\begin{array}{rcl}\textrm{if }\forall a,\, i_a>0&: &\kappa_{i_1, \ldots, i_d} 
\\
\textrm{otherwise }&: &0
    \end{array}\right. .
\ee
The $\kappa$ being symmetric in their arguments, $\mathcal{D}$ is symmetric too. We extend the definition of $\mathcal{D}$ to $\mathbb{K}^d$ by requiring the following:
\begin{itemize}
\item $\mathcal{D}$ is multilinear: let $\{V^1_{i_1}, \ldots , V^d_{i_d}\}_{\{i_a\}}$ and $\{P^1_{i_1}, \ldots , P^d_{i_d}\}_{\{i_a\}}$ such that for any $a, i_a$,  $V^a_{i_a}\in\mathbb K$, and  $P^a_{i_a}$ is some formal series in the $\{Y_i\}$ that does not involve first order free cumulants:
\be
\mathcal{D}\Bigl(\sum_{i_1\ge 0} P^1_{i_1}V^1_{i_1}, \ldots, \sum_{i_d\ge 0} P^d_{i_d}V^d_{i_d} \Bigr) = \sum_{i_1\ge 0} P^1_{i_1}\cdots \sum_{i_d\ge 0} P^d_{i_d} \, \mathcal{D}(V^1_{i_1}, \ldots, V^d_{i_d}).
\ee
\item  $\mathcal{D}$ satisfies a Leibniz rule on each one of its arguments, that is, considering $F_1, F_2 \in \bar{\mathbb F}$: 
$$
\mathcal{D}(F_1 F_2, \ldots ) = F_1 \mathcal{D}(F_2, \ldots ) +F_2 \mathcal{D}(F_1, \ldots ).
$$
\end{itemize}
In other words, $\mathcal{D}$ is a multilinear partial differential form (since the higher order free cumulants are seen as coefficients) acting on the first order free cumulants as variables and considering any polynomial in the $Ys$ as constant (including higher order free cumulants).  From the rules above, we find that if $F_1, \ldots, F_d \in \bar{\mathbb F}$:
\be
\frac{\partial}{\partial Y_i} \mathcal{D}(F_1 , \ldots , F_d) =\sum_{a=1}^d  \mathcal{D}\Bigl(F_1 , \ldots , \frac{\partial F_a}{\partial Y_i} , \ldots,  F_d\Bigr).
\ee
If different arguments of $\mathcal{D}$ involve different variables $Y_i$, ${\partial}/{\partial Y_i}$ commutes with $\mathcal{D}$.

From the Leibniz rule,   $\mathcal{D}$ satisfies the chain rule on each argument, that is, if for instance $F_1(X)=\sum_k c_k X^k$ and $H_1 \in \mathbb K$, then $\mathcal{D}\bigl(F_1(H_1), \ldots \bigr)  = F_1'(H_1) \mathcal{D}\, \bigl(H_1, \ldots \bigr)$,  
and if $F_1(X_1, \ldots, X_q)=\sum_k c_{k_1, \ldots, k_q} X_1^{k_1} \cdots X_q^{k_q}$, then: 
\be
\label{eq:chain-rule-multivar}
\mathcal{D}\bigl(F_1(H_1, \ldots, H_q), \ldots \bigr) = \sum_{\theta=1}^q \frac{\partial F_1}{\partial X_\theta} (H_1, \ldots, H_q)\  \mathcal{D}\bigl(H_\theta, \ldots \bigr). 
\ee

In practice, we will use $\mathcal{D}$ through its action on  $C_1$ as follows:
\be
\mathcal{D}\bigl(C_1(Y_1), \ldots, C_1(Y_d)\bigr) = C_{d}(Y_1, \ldots, Y_d), 
\ee
and applied to functions of the $C_1(Y_i)$.

\begin{lemma}
\label{lem:bar-from-tilde}
For $I\subset\{1, \ldots, p\}$ and $ \bar \pi \in \mathcal{P}^\star(I)$, the several-blocks corrections $\bar C_{\lvert I \rvert, \bar \pi}$ can be derived from the one-block corrections $\tilde C_{\lvert B \rvert}$ for the $B\in \bar \pi$ (with $\tilde C_{1}= C_1$) as:
\be
\bar C_{\lvert I \rvert, \bar \pi}\bigl(Y_I\bigr) = \mathcal{D}\Bigl(\{\tilde C_{\lvert B \rvert}(Y_B)\}_{B\in \bar \pi}\Bigr).
\ee 
\end{lemma}
\proof This is seen directly on \eqref{eq:bar-C}. With the same notations:
\be
 \kappa_{\{\lvert F_G\rvert\}_{G\in \bar \pi}}[b]  =
 \mathcal{D}\Bigl(\{\kappa_{\lvert F_G\rvert}  \}_{G\in \bar \pi}\Bigr),
\ee
so that since the Leibniz rule holds for each argument:
\be
\sum_{\{F_G\in \Pi(\nu_G)\}_{G\in \bar \pi}}  \kappa_{\{\lvert F_G\rvert\}_{G\in \bar \pi}}  
 \prod_{\substack{{V \in \Pi(\bigtimes_G \nu_G)}\\{V \notin \{F_G\}} }} \kappa_{\lvert V \rvert}  = \mathcal{D}\Bigl(\Bigl\{
 \kappa\bigl(\Pi(\nu_G), \nu_G\bigr) 
 \Bigr\}_{G\in \bar \pi}\Bigr).
\ee
By multilinearity:
\be
\bar C_{\lvert I \rvert, \bar \pi} \bigl(\{Y_i\}_{i\in I}\bigr) =  \mathcal{D}\Biggl(\Biggl\{\sum_{\{\mu_i\}_{i\in G}} \  \prod_{i\in G} Y_i^{\mu_i}\ \Bigl[  \sum_{\nu_G\in \mathrm{NS} (\mu_G)}\kappa\bigl(\Pi(\nu_G), \nu_G\bigr) \Bigr]\Biggr\}_{G\in \bar \pi}\Biggr), 
\ee
where we recognize \eqref{eq:bar-C-one-block}. \qed

\ 

This provides a more efficient way of computing the corrections $\bar C_{\lvert I \rvert, \bar \pi}$ from known $\tilde C_p$. For instance, let $\bar C_{p, \{a,b\}}$ be $\bar C_{p, \bar \pi},$ for $\bar\pi$ the partition of $\{1, \ldots, p\}$ whose blocks all contain a single element, apart from the block $\{a,b\}$. $\bar C_{p, \{a,b\}}$ can be computed from the coefficients as $\bar C_{3, \{1\}\{2,3\}}$, but applying Prop.~\ref{lem:bar-from-tilde}, we find that 
\be
\label{eq:term-p-1-order-p-D}
\bar C_{p, \{a,b\}} \bigl(Y_1, \ldots, Y_p\bigr) = \mathcal{D}\Bigl(\tilde C_{2}(Y_a, Y_b), \{C_1(Y_i)\}_{i \neq a,b}\Bigr).
\ee
Using the expression of $\tilde C_2$ in \eqref{eq:C2Tilde-log}: 
\begin{align}
\label{eq:bar-C-ab}
\bar C_{p, \{a,b\}} \bigl(Y_1, \ldots, Y_p\bigr) = Y_a Y_b  \frac{\partial}{\partial Y_a} \frac{\partial}{\partial Y_b} \left[  \frac 1 {1-\hat C(Y_a, Y_b)} \mathcal{D}\Bigl(\hat C(Y_a, Y_b), \{C_1(Y_i)\}_{i \neq a,b}\Bigr) \right],
\end{align}
and from the expression of $\hat C$ \eqref{eq:hat-C}:
\be
\label{eq:bar-C-ab2}
\mathcal{D}\Bigl(\hat C(Y_a, Y_b), \{C_1(Y_i)\}_{i \neq a,b} \Bigr)= \frac{Y_a}{Y_b - Y_a} C_{p-1}(\{Y_i\}_{i\neq a}) + \frac{ Y_b} {Y_a - Y_b}  C_{p-1}(\{Y_i\}_{i\neq b}).
\ee
In the same way, we can re-derive $\bar C_{4, \{1,2\}\{3,4\}}$  \eqref{eq:eq:expr-C-4-22-prop} to verify that everything adds up: 
\be
\bar C_{4, \{1,2\}\{3,4\}}(Y_1, \ldots, Y_4)=  \mathcal{D}\Bigl(\tilde C_{2}(Y_1, Y_2),\tilde C_{2}(Y_3, Y_4)\Bigr).
\ee
Using again the expression of $\tilde C_2$ in \eqref{eq:C2Tilde-log}:
\be
\bar C_{4, \{1,2\}\{3,4\}}(Y_1, \ldots, Y_4)=  \left(\prod_{i=1}^4 Y_i \frac{\partial}{\partial Y_i} \right)\left[\frac{ \mathcal{D}\Bigl(\hat C(Y_1, Y_2),\hat C(Y_3, Y_4)\Bigr)}{(1 - \hat C(Y_1, Y_2))(1 - \hat C(Y_3, Y_4)}\right],
\ee
where indeed,  $ \mathcal{D}_{2}\bigl(\hat C(Y_1, Y_2),\hat C(Y_3, Y_4)\bigr)= \hat C_{\{1,2 \}\{3,4\}}$ defined in \eqref{eq:C-hat-22}. 

\ 

More generally, we define for $p\ge 2$:
\be
\hat C (Y_1, \ldots, Y_p) = \sum_{i_1, \ldots, i_p} \kappa_{i_1 + \cdots + i_p}  \, Y_1^{i_1} \cdots Y_p^{i_p}
\ee
and the higher order generalizations for $\pi \in \mathcal{P}(p)$: 
\be
\hat C_{\pi} (Y_1, \ldots, Y_p) = \sum_{i_1, \ldots, i_p} \kappa_{\{\sum_{a\in G}i_a\}_{G\in \pi}}  \,Y_1^{i_1} \cdots Y_p^{i_p}.
\ee
\begin{lemma} 
\label{lem:The-hats}
These generating functions are computed as
\be
\label{eq:hat-C-gen}
\hat C (Y_1, \ldots, Y_p) = (-1)^p +  \sum_{i=1}^p C_1(Y_i)\prod_{j\neq i} \frac{Y_j}{Y_i - Y_j},
\ee
and if $\pi \in \mathcal{P}(p)$ with $\#(\pi)=d>1$:
\be
\label{eq:hat-C-pi-gen}
\hat C_{\pi} (Y_1, \ldots, Y_p) =  \sum_{\{i_G \in G\}_{_{G\in \pi}}} C_d(\{Y_{i_G}\}_{G\in \pi}) \prod_{G\in \pi}\prod_{j_G\neq i_G} \frac{Y_{j_G}}{Y_{i_G} - Y_{j_G}}. 
\ee
\end{lemma}
\proof As already done several times:
\be
(Y_1 - Y_2) \hat C (Y_1, \ldots, Y_p) = Y_2 \hat C (Y_1, Y_3, \ldots, Y_p) - Y_1 \hat C (Y_2,  \ldots, Y_p), 
\ee
and \eqref{eq:hat-C-gen} follows by induction noticing that for $i\neq 1, 2$, 
$$
\frac{Y_2}{Y_1 -Y_2}\frac{Y_1}{Y_i -Y_1} + \frac{Y_1}{Y_2 -Y_1}\frac{Y_2}{Y_i -Y_2} = \frac{Y_1}{Y_i -Y_1}\frac{Y_2}{Y_i -Y_2}. 
$$
Equation \eqref{eq:hat-C-pi-gen} follows directly from \eqref{eq:hat-C-gen} since from the definition:
$$
\hat C_{\pi} (Y_1, \ldots, Y_p) = \mathcal{D}\bigl(\{\hat C(Y_B)\}_{B\in \pi}\bigr).
$$
\qed

In the same way, we can easily compute $\bar C_{4,\{1,2,3\}\{4\}}$, $\bar C_{5,\{1,2,3\}\{4,5\}}$ and $\bar C_{6,\{1,2,3\}\{4,5,6\}}$ from $\tilde C_2$ and $\tilde C_3$, whose expressions we have derived explicitly in Prop.~\ref{prop:gen-fun-third-order}. In this form, the explicit expressions get two lengthy so we don't write them here  but $\bar C_{4,\{1,2,3\}\{4\}}$ will be computed this way in  \eqref{eq:hat-C3-13}  after a simpler expression for $\tilde C_3$ is obtained in Prop.~\ref{prop:simpl-third-order}.

\subsubsection{Trees with several blue vertices}

We now turn to the general case of $\bar C_{\lvert I \rvert, \bar \pi, T}$ for $T$ with more than one blue vertex. 

\begin{proposition}
\label{prop:bar-from-tilde-gen}
For $I\subset\{1, \ldots, p\}$,  $ \bar \pi \in \mathcal{P}^\star(I)$, and $T\in \mathcal G_{\bar \pi}$, the several-blocks corrections $\bar C_{\lvert I \rvert, \bar \pi, T}$ can be derived from the one-block corrections $\tilde C_{\lvert B \rvert}$ for the $B\in \bar \pi$ (with $\tilde C_{1}= C_1$) as:
\be
\bar C_{\lvert I \rvert, \bar \pi, T}\bigl(Y_I\bigr) = \sum_{\{j_G^K\}_{\substack{{K \in \mathcal{I}(T)}\\{G\in K}}}}\, \prod_{K \in \mathcal{I}(T)} \kappa_{\{j_G^K\}_{G\in K}}\prod_{G\in \bar \pi} \biggl(\prod_{K \in \mathcal{I}_G(T)}\frac{\partial}{\partial \kappa_{j_G^K}}\biggr)  \Bigl(\tilde C_{\lvert G \rvert}(Y_G)\Bigr).
\ee 
\end{proposition}
This formula can also be understood and used using $\mathcal D$ as for Lemma~\ref{lem:bar-from-tilde}, as discussed after the proof. 

\proof 
Our starting point is \eqref{eq:bar-f-tree-1}, where we recall that $\tau_G=\bigtimes_{i\in G} \bigl({\tau_i}_{\lvert_{B}}\bigr)$ (which is also  $\tau_G={\gamma_\lambda}_{\lvert_G}$). Grouping the terms in this formula according to the sizes of sizes of the blocks attributed to the edges, we obtain the following:
\be 
\label{eq: proof-D-tree-1}
\bar \kappa _{\bar \pi, T}  (\tau_{\lvert_B})  = \sum_{{\substack{ {\{j_G^K\}_{G\in \bar \pi,\,  K \in \mathcal{I}_G(T)}}\\[+0.2ex]{\forall G,\, \sum_{K \in \mathcal{I}_G} j_G^K \le \lvert G \rvert} }}}  \ \prod_{K \in \mathcal{I}(T)}  \kappa_{\{j_G^K\}_{G\in K}}  \prod_{G\in \bar \pi}  \Biggl(\sum_{\substack{{\nu_G\in \mathrm{NS} (\tau_G)}\\[+0.4ex]{\#(\nu_G) \ge \mathrm{deg}_T(G) }}} \ \sum_{\substack{{\mathcal L_G }\\{\textrm{inj. att.}}\\{\forall K,\, \lvert \mathcal L_G(K)\rvert =j_G^K}}} \  \ \prod_{\substack{{V\in \Pi(\nu_G)}\\{V\notin \mathcal L_G( \mathcal{I}_G(T))} }}  \kappa_{\lvert V\rvert} \Biggr).  
\ee
Note that the sum over attributions can vanish if $\Pi(\nu_G)$ does not have a block of size $j_G^K$ for each $K\in \mathcal I_G(T)$. Recalling that \eqref{eq:bar-f-one-block}:
\be
\tilde \kappa (\tau_G)  = \sum_{\nu\in \mathrm{NS}\left(\tau_G\right)} \ 
\prod_{V\in \Pi(\tau_G)} \kappa_{\lvert V \rvert }, 
\ee
we compute for $q\ge 1$:
\be
\frac{\partial}{\partial \kappa_{q}} \tilde \kappa (\tau_G)  = \sum_{\nu \in \mathrm{NS} (\tau_G) } \ \sum_{\substack{{\tau \in \Pi(\nu)}\\{\lvert \theta \rvert =q}}} \ \prod_{\substack{{V \in \Pi(\nu)}\\{V\neq \tau}}} \ \kappa_{\lvert V \rvert },
\ee
with the convention that the sum over $\tau$ vanishes if $\nu$ does not have a cycle of length $q$. Differentiating twice:
\be
\frac{\partial}{\partial \kappa_{q_1}} \frac{\partial}{\partial \kappa_{q_2}} \tilde \kappa (\tau_G)  = \sum_{\nu \in \mathrm{NS} (\tau_G) } \ \sum_{\substack{{\tau_1, \tau_2 \in \Pi(\nu)}\\{\lvert \theta_i \rvert =q_i}}} \ \prod_{\substack{{V \in \Pi(\nu)}\\{V\neq \tau_1, \tau_2}}} \ \kappa_{\lvert V \rvert },
\ee
and so on, so that we can replace the term in parenthesis in \eqref{eq: proof-D-tree-1} to obtain:
\be 
\label{eq: proof-D-tree-2}
\bar \kappa _{\bar \pi, T}  (\tau_{\lvert_B})  = \sum_{{\substack{ {\{j_G^K\}_{G\in \bar \pi,\,  K \in \mathcal{I}_G(T)}}\\[+0.2ex]{\forall G,\, \sum_{K \in \mathcal{I}_G} j_G^K \le \lvert G \rvert} }}}  \ \prod_{K \in \mathcal{I}(T)}  \kappa_{\{j_G^K\}_{G\in K}}  \prod_{G\in \bar \pi}  \,  \biggl(\prod_{K \in \mathcal{I}_G(T)}\frac{\partial}{\partial \kappa_{j_G^K}}\biggr)\,  \tilde \kappa (\tau_G) .  
\ee
The expression for $\bar \kappa _{\bar \pi, T}  (\{\mu_i\}_{i\in I}) $ is obtained replacing in this expression $\lvert G \rvert$ by $\sum_{i\in G} \mu_i$, and $\tau_G$ by $\mu_G=\{\mu_i\}_{i\in G}$. Going then to the generating function  \eqref{eq:bar-C-T}:
$$
\bar C_{\lvert I \rvert, \bar \pi, T}\bigl(Y_I\bigr)= \sum_{\{\mu_i\}_{i\in I}} \bar \kappa _{\bar \pi, T}  (\{\mu_i\}_{i\in I}) \, \prod_{i\in I} Y_i^{\mu_i}, 
$$
and exchanging the summations over the $\mu_i$ and the $j_G^K$, we find: 
\be
\bar C_{\lvert I \rvert, \bar \pi, T}\bigl(Y_I\bigr)= \sum_{\{j_G^K\}_{\substack{{K \in \mathcal{I}(T)}\\{G\in K}}}}\, \prod_{K \in \mathcal{I}(T)} \kappa_{\{j_G^K\}_{G\in K}}\prod_{G\in \bar \pi} \Biggl(\sum_{\substack{{\{\mu_i\}_{i\in G}}\\{\sum_{i\in G} \mu_i \ge \sum_{K\in \mathcal I_G} j_G^K}}}\prod_{i\in G} Y_i^{\mu_i}\biggl(\prod_{K \in \mathcal{I}_G(T)}\frac{\partial}{\partial \kappa_{j_G^K}}\biggr)  \tilde\kappa (\mu_G)\Biggr), 
\ee
and factorizing the product of derivatives for each $G$ leads to the result. \qed

\

For $I\subset\{1, \ldots, p\}$,  $ \bar \pi \in \mathcal{P}^\star(I)$, and $T\in \mathcal G_{\bar \pi}$, and considering $K_1\in \mathcal I(T)$, we see recalling that $\kappa_{i_1, \ldots, i_q}=\mathcal{D}(\kappa_{i_1}, \ldots, \kappa_{i_q})$ applying  \eqref{eq:chain-rule-multivar} for each argument that:
\be
\label{eq:first-layer-K1}
 \sum_{\{j_G^{K_1}\}_{G\in K_1}}\, \kappa_{\{j_G^{K_1}\}_{G\in K_1}}\prod_{G\in K_1} \biggl(\frac{\partial}{\partial \kappa_{j_G^{K_1}}}\biggr)  \Bigl(\tilde C_{\lvert G \rvert}(Y_G)\Bigr) = \mathcal{D} \Bigl(\bigl\{\tilde C_{\lvert G \rvert}(Y_G)\bigr\}_{G\in K_1}\Bigr) .
\ee 
We now want to act with the derivatives  corresponding to another blue vertex $K_2\in \mathcal{I}(T)$ to obtain
\be
 \sum_{\{j_G^{K_2}\}_{G\in K_2}}\, \kappa_{\{j_G^{K_2}\}_{G\in K_2}} \sum_{\{j_G^{K_1}\}_{G\in K_1}}\, \kappa_{\{j_G^{K_1}\}_{G\in K_1}}\prod_{G\in K_1 \cup K_2} \biggl(\frac{\partial}{\partial \kappa_{j_G^{K_1}}}\biggr)^{\epsilon_1(G)}\biggl(\frac{\partial}{\partial \kappa_{j_G^{K_2}}}\biggr)^{\epsilon_2(G)}  \Bigl(\tilde C_{\lvert G \rvert}(Y_G)\Bigr),
\ee 
where $\epsilon_i(G)$ is $1$ if $G\in K_i$ and 0 otherwise,
and we would like to express this  as a differential operator acting on what has already been computed. We introduce two kinds  of operators on the tensor product $\bigotimes_{G\in \bar \pi} \mathbb{K}$. We define a $\mathcal{D}^\otimes_{K}:\mathbb K^{\otimes \#(\bar \pi)} \rightarrow \mathbb K^{\otimes \#(\bar \pi)} $ for  any $K \in \mathcal I(T)$:
\be
\mathcal{D}^\otimes_{K}  =   \sum_{\{j_G\}_{G\in K}}  \kappa_{\{j_G\}_{G\in K}} \biggl(\un_{\hat K} \otimes \bigotimes_{G\in K} \frac{\partial}{\partial_{\kappa_{j_G}} }\biggr), 
\ee
 where $\hat K$ denotes the set of blocks of $\bar \pi$ that are not in $K$. More explicitly:
\be
\mathcal{D}^\otimes_{K} \biggl(\, \bigotimes_{G\in \bar \pi} f_G(Y_G)\biggr) =  \sum_{\{j_G\}_{G\in K}}  \kappa_{\{j_G\}_{G\in K}} \biggl(\, \bigotimes_{G\in \hat K} f_G(Y_G)\biggr) \otimes \biggl(\,  \bigotimes_{G\in K} \frac{\partial}{\partial_{\kappa_{j_G}} }  f_G(Y_G) \biggr). 
\ee
We also let $\mathcal P: \mathbb K^{\otimes \#(\bar \pi)} \rightarrow \mathbb K$ be the multilinear map which gives the product of the elements:
\be
\mathcal P\biggl(\, \bigotimes_{G\in \bar \pi} f_G(Y_G)\biggr) = \prod_{G\in \bar \pi} f_G(Y_G).
\ee
We clearly have from \eqref{eq:first-layer-K1} for $K\in \mathcal{I}(T)$: 
\be
\mathcal{D}(\{f_G(Y_G)\}_{G\in K})= \bigl(\mathcal P \circ \mathcal{D}^\otimes_{K}\bigr) \biggl(\, \bigotimes_{G\in K} f_G(Y_G)\biggr). 
\ee
More generally, the following follows directly from Prop.~\ref{prop:bar-from-tilde-gen}.
\begin{corollary}
\label{cor:tensor-product}
For $I\subset\{1, \ldots, p\}$,  $ \bar \pi \in \mathcal{P}^\star(I)$, and $T\in \mathcal G_{\bar \pi}$, the several-blocks corrections $\bar C_{\lvert I \rvert, \bar \pi, T}$ are computed as:
\be
\bar C_{\lvert I \rvert, \bar \pi, T}\bigl(Y_I\bigr) = \mathcal P \circ \Bigl(\prod_{K\in \mathcal{I(T)}} \mathcal{D}^\otimes_{K}\Bigr) \biggl(\, \bigotimes_{G\in K} \tilde C_{\lvert G \rvert}(Y_G) \biggr).
\ee 
Note that the $\mathcal{D}^\otimes_{K}$ commute.
\end{corollary}

\emph{In practice, we don't need to carry out the sums over coefficients in the computations. } Indeed, assume that we have an expression of $\tilde C_q$ in terms of generating functions such as the $\hat C$, or other generating functions linear in $\{\kappa_i\}_{i\ge 1}$, involving explicit coefficients as $Y_i$, $\partial / \partial Y_i$. Then the coefficients involving functions of the $Y_i$, $\partial / \partial Y_i$, etc, can be factorized as for $\mathcal{D}$, letting the derivatives $\partial / \partial \kappa_q$ act on the $\hat C$ or similar generating functions, resulting in a term that does not depend in the $\{\kappa_i\}_{i\ge 1}$ and is therefore also treated as a coefficient. We have for instance for $J\subset \{1, \ldots, p\}$:
\be
\frac{\partial \hat C(Y_{J}) }{\partial \kappa_q} = \sum_{i\in J}\, Y_i ^q\,  \prod_{j\neq i }\frac{Y_j}{Y_i - Y_j}.
\ee
The sums over $j_G^K$in Prop.~\ref{prop:bar-from-tilde-gen} or Cor.~\ref{cor:tensor-product} can therefore be carried out explicitly outside of the tensor product. We have for instance, considering a $J_G\subset G$ for each $G\in K$ and $J=\cup_{G\in K} J_G$ (Lemma~\ref{lem:The-hats}):
\be
\label{eq:sum-hat}
\sum_{\{j_G^K\}_{G\in K}} \kappa_{\{j_G^K\}_{G\in K}} \prod_{G\in K}\frac{\partial \hat C(Y_{J_G}) }{\partial \kappa_{j_G^K}} = \hat C_{\{J_G\}_{G\in K}}(Y_J) = \mathcal{D}\bigl(\{\hat C(Y_{J_G})\}_{G\in K}\bigr).
\ee
It should be thought of as a chain rule seing the $\tilde C_{\lvert G \rvert}(Y_G)$ as composed with the $\hat C$:
\be
\mathcal{D}^\otimes_{K}  =   \sum_{\{J_G\subset G\}_{G\in K}}  \mathcal{D}\bigl(\{\hat C(J_G)\}_{G\in K}\bigr) \biggl(\un_{\hat K} \otimes \bigotimes_{G\in K} \frac{\partial}{\partial{\hat C(J_G)} }\biggr), 
\ee
\emph{where the full operator commutes with the $\partial / \partial Y_i$ but not }$ \mathcal{D}\bigl(\{\hat C(J_G)\}_{G\in K}\bigr)$.

\begin{figure}[!h]
\centering
\includegraphics[scale=0.8]{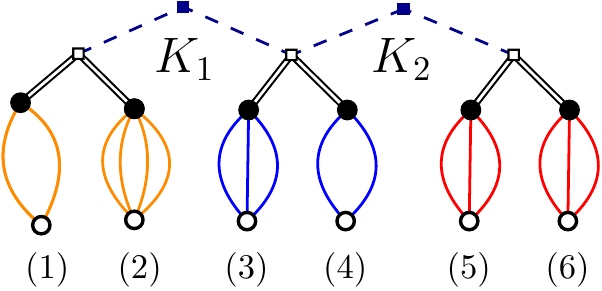}
\caption{The treated example $\bar C_{6, \{1,2\}\{3,4\}\{5,6\}, T}$. 
}
\label{fig:Ex-order-6}
\end{figure}

\

 This is better seen on an example: let us compute $\bar C_{6, \{1,2\}\{3,4\}\{5,6\}, T}$, for $T$ as shown in Fig.~\ref{fig:Ex-order-6}. From Cor.~\ref{cor:tensor-product}, we have:
\be
\bar C_{6, \{1,2\}\{3,4\}\{5,6\}, T}(Y_1, \ldots, Y_6) = \mathcal P \circ \Bigl(\mathcal{D}^\otimes_{K_2}\, \mathcal{D}^\otimes_{K_1}\Bigr) \biggl( \tilde C_{2}(Y_1, Y_2) \otimes  \tilde C_{2}(Y_3, Y_4) \otimes  \tilde C_{2}(Y_5, Y_6) \biggr). 
\ee 
We first let $\mathcal{D}^\otimes_{K_1}$ act on the argument, replacing $\tilde C_2$ by the expression \eqref{eq:C2Tilde-log}, giving:
\begin{align}
&\mathcal{D}^\otimes_{K_1} \Bigl( \tilde C_{2}\otimes  \tilde C_{2}\otimes  \tilde C_{2} \Bigr) = \prod_{i=1}^4 Y_i \frac{\partial}{\partial Y_i} \Biggl(\mathcal{D}^\otimes_{K_1} \biggl( \log(1-\hat C(Y_1, Y_2)) \otimes  \log(1-\hat C(Y_3, Y_4))\otimes  \tilde C_{2}(Y_5, Y_6) \biggr)
\\&\hspace{2.5cm}\nonumber\prod_{i=1}^4 Y_i \frac{\partial}{\partial Y_i} \Biggl(  \hat C_{\{1,2\}\{3,4\}}(Y_1, \ldots Y_4) \biggl( \frac{1}{1 - \hat C(Y_1, Y_2)}\otimes \frac{1}{1 - \hat C(Y_3, Y_4)} \otimes \tilde C_2(Y_5, Y_6) \biggr)\Biggr).\end{align}
Letting now $\mathcal{D}^\otimes_{K_2}$ act on the result:
\begin{align}
&\mathcal{D}^\otimes_{K_2} \biggl(\mathcal{D}^\otimes_{K_1} \Bigl( \tilde C_{2}\otimes  \tilde C_{2}\otimes  \tilde C_{2} \Bigr) \biggr)  = \prod_{i=1}^6 Y_i \frac{\partial}{\partial Y_i} \Biggl(\hat C_{\{1,2\}\{3,4\}}(Y_1, \ldots Y_4) \hat C_{\{3,4\}\{5,6\}}(Y_3, \ldots Y_6)
\\&\hspace{6cm}\nonumber \biggl( \frac{1}{1 - \hat C(Y_1, Y_2)}\otimes \frac{1}{(1 - \hat C(Y_3, Y_4))^2} \otimes \frac{1}{1 - \hat C(Y_5, Y_6)} \biggr)\Biggr).\end{align}
Applying $\mathcal P$, we get the final result:
\be
\bar C_{6, \{1,2\}\{3,4\}\{5,6\}, T} (Y_1, \ldots, Y_6)= \prod_{i=1}^6 Y_i \frac{\partial}{\partial Y_i} \Biggl(\frac{\hat C_{\{1,2\}\{3,4\}}(Y_1, \ldots Y_4) \hat C_{\{3,4\}\{5,6\}}(Y_3, \ldots Y_6)}{(1 - \hat C(Y_1, Y_2))(1 - \hat C(Y_3, Y_4))^2(1 - \hat C(Y_5, Y_6))}\Biggr).
\ee

\subsection{Simplifications}
\label{sub:simplifications}
The relations of Borot et al.~\cite{Analytic-higher} reformulated in Thm.~\ref{thm:theorem-functional-reformulated} are much simpler however, as some drastic simplifications occur. From Thm.~\ref{thm:functional-relations-for-any-H} and Thm.~\ref{thm:facto}, we see that the combinatorial proof of Thm.~\ref{thm:analytic-higher} equivalently formulated as Thm.~\ref{thm:theorem-functional-reformulated} boils down to showing the following remarkable simplifications:
\be
\label{eq:simplifications}
\sum_{T\in \mathcal{G}_p} \biggl( \prod_{i=1}^p 
  X_i^{\mathrm{d}_T(i) -1} \frac {\mathrm d^{\mathrm{d}_T(i) -1}} {\mathrm d X_i^{\mathrm{d}_T(i) -1}}\biggr)\Biggl\{\prod_{i=1}^p
   {\frac {\mathrm{d} Y_i} {\mathrm d X_i}}\frac 1 {C(Y_i)^{\mathrm{d}_T(i)}} \Biggl(\prod_{I\in \mathcal{I}(T)} H_{\lvert I \rvert}\bigl(Y_I\bigr) - \prod'_{I\in \mathcal{I}(T)} \ C_{\lvert I \rvert }\bigl(Y_I\bigr)\Biggr)\Biggr\} = 0,   
\ee
where $\mathrm{d}_T(i) = \mathrm{deg}_T(i)$ and $Y_I=\{Y_i\}_{i\in I}$, $H_{\lvert I \rvert}$ is given by \eqref{eq:def-HI}
\be
H_{\lvert I \rvert} = C_{\lvert I \rvert} + \sum_{\bar \pi \in \mathcal{P}^*(I)} \bar H_{\lvert I \rvert, \bar \pi}, \qquad \bar H_{\lvert I \rvert, \bar \pi} = \sum_{T \in \mathcal{G}_{\bar \pi}} 
\bar C_{\lvert I \rvert, \bar \pi, T},
\ee
where $ \bar C_{\lvert I \rvert, \bar \pi, T}$ is defined in \eqref{eq:bar-C-T}, and $\prod'$ means that any occurence of $C_2(Y_1, Y_2)$ should be replaced by
\be
\label{eq:correction-c2} 
C_2 (Y_1, Y_2) \quad \rightarrow \quad C_2 (Y_1, Y_2) + \frac {Y_1Y_2}{(Y_1 - Y_2)^2}.
\ee
If $T$ has $I\in \mathcal{I}(T)$ with two elements (black vertices of valency two), the term $\prod_{I} H_{\lvert I \rvert} - \prod_{I} \ C_{\lvert I \rvert }$ is a sum of terms of the form $\prod_{I} O_{\lvert I \rvert} $, where $O_{\lvert I \rvert}$ is either $C_{\lvert I \rvert }$, or $\bar C_{\lvert I \rvert, \bar \pi}$ for some $\bar \pi \in \mathcal{P}^*(I)$, with the condition that not all  $O_{\lvert I \rvert}$ are  $C_{\lvert I \rvert }$. The corrections \eqref{eq:correction-c2} must also be taken into account. 

\

Because there are no relations among the $C_p(Y_I)$, it should be true that the cancellations in \eqref{eq:simplifications} occur ``order-by-order'' in the $C_k(Y_I)$, that is, all the terms involving for instance $C_3(Y_1, Y_2, Y_3), C_2(Y_3, Y_4),  C_2(Y_4, Y_5)$ should cancel \emph{among each other}. Since there are no relations either among the $C_p$ and their derivatives for $p>1$, we expect the the cancellations to occur for each coefficient of an expression of the form $\frac{\partial}{\partial Y_1}C_3(Y_1, Y_2, Y_3) \frac{\partial}{\partial Y_1} C_2(Y_3, Y_4)  C_2(Y_4, Y_5)$. 
  In particular, \emph{we expect that the terms that only involve $C_1$ should cancel among each other, implying an identity satisfied by $\tilde C_p$ for $p>2$}: 
\begin{align}
\label{eq:identity}
&\sum_{T\in \mathcal{G}_p} \biggl( \prod_{i=1}^p 
  X_i^{\mathrm{d}_T(i) -1} \frac {\mathrm d^{\mathrm{d}_T(i) -1}} {\mathrm d X_i^{\mathrm{d}_T(i) -1}}\biggr)\Biggl\{\prod_{i=1}^p
   {\frac {\mathrm{d} Y_i} {\mathrm d X_i}}\frac 1 {C(Y_i)^{\mathrm{d}_T(i)}} \Biggl(\prod_{I\in \mathcal{I}(T)} \tilde C_{\lvert I \rvert}\bigl(Y_I\bigr)\\&\hspace{7cm\nonumber } - \prod_{I\in \mathcal{I}(T)} \prod_{i_1, i_2 \in \{1, \ldots, p\}} \delta_{I, \{i_1, i_2\}} \frac{Y_{i_1}Y_{ i_2}}{(Y_{i_1} - Y_{i_2})^2})\Biggr)\Biggr\} \ \raisebox{1ex}{$\substack{{?}\\{=}}$} \ 0,   
\end{align}
The last term is non-zero only if all the black vertices of the tree $T$ are of valency two.
This identity would correspond to a recurrence satisfied by the $\tilde C_p$, as the term for the tree of $\mathcal{G}_p$  that has a single black vertex and all white vertices of valency one is just 
$$
\prod_{i=1}^p {\frac {\mathrm{d} Y_i} {\mathrm d X_i}}\frac 1 {C(Y_i)} \tilde C_p(Y_1, \ldots, Y_p),
$$
and all the other terms involve $\tilde C_k$ for $k<p$. The relation can be extended to $\tilde C_2$ but it takes the different form  \eqref{eq:Tilde-C2-in-pat-ex}: 
\be
\label{eq:eq:order1-of-order2}
 \tilde C_2(Y_1, Y_2) - \frac{Y_1 Y_2}{(Y_1 - Y_2)^2} = -\frac{Y_1Y_2}{(X_1 - X_2)^2} \frac {\mathrm{d} X_1} {\mathrm d Y_1} \frac {\mathrm{d} X_2} {\mathrm d Y_2},
\ee 
instead of zero.

We then expect the contributions to \eqref{eq:simplifications} involving higher order free cumulants to \emph{derive from \eqref{eq:identity} by application of some differential operators on the identities satisfies by the $\tilde C_k$ for }$k<p$. For instance, the other contributions to  \eqref{eq:simplifications}  for the tree of $\mathcal{G}_p$  that has a single black vertex and all white vertices of valency one are of the form 
$$
\prod_{i=1}^p {\frac {\mathrm{d} Y_i} {\mathrm d X_i}}\frac 1 {C(Y_i)} \bar C_{\lvert I \rvert, \bar \pi, T}  =  \prod_{i=1}^p {\frac {\mathrm{d} Y_i} {\mathrm d X_i}}\frac 1 {C(Y_i)} \mathcal P \circ \Bigl(\prod_{K\in \mathcal{I(T)}} \mathcal{D}^\otimes_{K}\Bigr) \biggl(\, \bigotimes_{G\in K} \tilde C_{\lvert G \rvert}(Y_G) \biggr) ,
$$
for  $I\subset\{1, \ldots, p\}$,  $ \bar \pi \in \mathcal{P}^\star(I)$, and $T\in \mathcal G_{\bar \pi}$, where we have applied Cor.~\ref{cor:tensor-product}. 

\

Below, we detail this for the orders $p=3$ and $p=4$.

\subsubsection{Third order}
At third order $p=3$ for instance, there are two kinds of trees in $\mathcal{G}_3$, one with one black vertex and three with two black vertices, such that one white vertex has valency two and the two others valency one. The term in \eqref{eq:simplifications} for the first kind of tree has a weight $H_3 - C_3$. From \eqref{eq:H3} (see also Fig.~\ref{fig:Third-order}), the contribution to  \eqref{eq:simplifications}  from that tree is:
\be
\label{eq:contrib-order3-simplif}
\prod_{i=1}^3   {\frac {\mathrm{d} Y_i} {\mathrm d X_i}}\frac 1 {C(Y_i)}\left( \tilde C_3 + \bar C_{3, \{1\}\{2,3\}}  + \bar C_{3, \{1,2\}\{3\}}  + \bar C_{3, \{2\}\{1,3\}}\right)(Y_1, Y_2, Y_3) ,
\ee
whose expressions have been computed in Prop.~\ref{prop:gen-fun-third-order}. The term $\tilde C_3$ has only contributions of order 1, whereas the other three terms have contributions in $C_2$. For a tree with two black vertices and such that the white vertex labeled 1 has valency two, we get a weight 
$$
H_2(Y_1, Y_2)H_2(Y_1, Y_3) - \Bigl(C_2 (Y_1, Y_2) + \frac {Y_1Y_2}{(Y_1 - Y_2)^2}\Bigr)\Bigl(C_2 (Y_1, Y_3) + \frac {Y_1Y_3}{(Y_1 - Y_3)^2}\Bigr), 
$$
where $H_2 = C_2 + \tilde C_2$ \eqref{eq:H2}, contributing to \eqref{eq:simplifications} as:
\be
X_1 \frac{\partial}{\partial X_1}\Biggl(\frac{C_2(Y_1, Y_2) \mathring{C}_2(Y_1, Y_3) + \mathring C_2(Y_1, Y_2) {C}_2(Y_1, Y_3) + \tilde C_2 (Y_1, Y_2)\tilde C_2 (Y_1, Y_3) - \frac{Y_1Y_2}{(Y_1 - Y_2)^2}\frac{Y_1Y_3}{(Y_1 - Y_3)^2}}{\prod_{i=1}^3 {\frac {\mathrm{d} X_i} {\mathrm d Y_i}}  C_1(Y_1)^2C_1(Y_2)C_1(Y_3)}\Biggr),
\ee
where:
\be
\label{eq:C2Circ}
\mathring C_2(Y_1, Y_2) = \tilde C_2(Y_1, Y_2) - \frac{Y_1 Y_2}{(Y_1 - Y_2)^2} = -\frac{Y_1Y_2}{(X_1 - X_2)^2} \frac {\mathrm{d} X_1} {\mathrm d Y_1} \frac {\mathrm{d} X_2} {\mathrm d Y_2}.
\ee 
The contributions for the two other similar trees are obtained exchanging 1 and 2, and 1 and 3. From \eqref{eq:simplifications}, the summation of these terms should vanish, and in fact the terms cancel ``order by order'', that is the contributions involving $C_2$ cancel among each other. We simplify by the prefactor in \eqref{eq:contrib-order3-simplif}.

\begin{proposition} 
\label{prop:simpl-third-order}
The following simplifications occur at third order: 
\be
\label{eq:eq:order1-of-order3}
\tilde C_3(Y_1, Y_2, Y_3) + \sum_{\substack{{1\leftrightarrow 2}\\{1\leftrightarrow 3}}} Y_1 \frac{\partial}{\partial Y_1}\left({\frac {\mathrm{d} Y_1} {\mathrm d X_1}} \frac{\tilde C_2 (Y_1, Y_2)\tilde C_2 (Y_1, Y_3) - \frac{Y_1Y_2}{(Y_1 - Y_2)^2}\frac{Y_1Y_3}{(Y_1 - Y_3)^2}}{C_1(Y_1)^2}\right) = 0,
\ee
which involves only first order free cumulants and corresponds to \eqref{eq:identity}, and
\begin{align}
\label{eq:eq:order2-of-order3}
&\bigl(\bar C_{3, \{1\}\{2,3\}}  + \bar C_{3, \{1,2\}\{3\}}  + \bar C_{3, \{1,3\}\{2\}}  \bigr)\bigl(Y_1, Y_2, Y_3\bigr) 
\\&\hspace{3cm}+ \sum_{\substack{{1\leftrightarrow 2}\\{1\leftrightarrow 3}}} Y_1 \frac{\partial}{\partial Y_1}\left({\frac {\mathrm{d} Y_1} {\mathrm d X_1}} \frac{ C_2(Y_1, Y_2) \mathring{C}_2(Y_1, Y_3) + \mathring C_2(Y_1, Y_2) {C}_2(Y_1, Y_3) }{C_1(Y_1)^2}\right) = 0,\nonumber
\end{align}
which involves second order free cumulants and where $\mathring C_2$ is given in \eqref{eq:C2Circ}. 
This concludes the combinatorial proof of Thm.~\ref{thm:theorem-functional-reformulated} for $p=3$.
\end{proposition}
 Indeed, the second relation \eqref{eq:eq:order2-of-order3} is easily verified from \eqref{eq:gen-func-pibardeuxterms} and is a particular case of Prop.~\ref{prop:order_p-1_at-orderp}. It is indeed true that the cancellations occur for each coefficient $C_2(Y_a, Y_b)$ and $\frac{\partial}{\partial Y_a}C_2(Y_a, Y_b)$ as can be seen from the proof of Prop.~\ref{prop:order_p-1_at-orderp}.
The first relation \eqref{eq:eq:order1-of-order3} is much more involved, but we have verified it from \eqref{eq:tilde-c3-gen} using Mathematica.  
 
 \ 
 
 By replacing the expressions of $\tilde C_2$ and recognizing the results of some derivatives, we find an expression of $\tilde C_3$ which is well adapted for applying differential operators such as $\mathcal D$ or $\mathcal P \circ \Bigl(\prod_{K\in \mathcal{I(T)}} \mathcal{D}^\otimes_{K}\Bigr)$:
\begin{align}
\nonumber
&\tilde C_3(Y_1, Y_2, Y_3) = \sum_{\substack{{1\leftrightarrow 2}\\{1\leftrightarrow 3}}} Y_1 \frac{\partial}{\partial Y_1}Y_2 \frac{\partial}{\partial Y_2}Y_3 \frac{\partial}{\partial Y_3}\Biggl(  \frac{Y_3}{Y_1 - Y_3} \frac{Y_2}{Y_1-Y_2}  \Biggl[  \frac 1{1-\hat C(Y_1, Y_2)}+  \frac 1{1-\hat C(Y_1, Y_3)} \\&\hspace{7cm}-  \frac{C_1(Y_1)^2 \frac {\mathrm{d} X_1} {\mathrm d Y_1}} {(1-\hat C(Y_1, Y_2))(1-\hat C(Y_1, Y_3))}\Biggr]\Biggr) .
\end{align}

From this expression, we can indeed easily compute for instance $\bar C_{3,\{1,2,3\}\{4\}} (Y_1, \ldots, Y_4)  = \mathcal{D}(\tilde C_3(Y_1, Y_2, Y_3), C_1(Y_4))$:
\begin{align}
\nonumber
&\bar C_{3,\{1,2,3\}\{4\}}(Y_1, \ldots, Y_4) = \sum_{\substack{{1\leftrightarrow 2}\\{1\leftrightarrow 3}}} Y_1 \frac{\partial}{\partial Y_1}Y_2 \frac{\partial}{\partial Y_2}Y_3 \frac{\partial}{\partial Y_3}\Biggl(  \frac{Y_3}{Y_1 - Y_3} \frac{Y_2}{Y_1-Y_2}  
\Biggl[  \frac{  Y_1 \frac {\partial} {\partial  Y_1}C_2(Y_1,Y_4)- C_2(Y_1, Y_4)} {(1-\hat C(Y_1, Y_2))(1-\hat C(Y_1, Y_3))}\\&
+ \frac {\hat C_{\{1,2\}\{4\}}(Y_1, Y_2, Y_4)}{(1-\hat C(Y_1, Y_2))^2}\Biggl(1 - \frac{C_1(Y_1)^2\frac {\mathrm{d} X_1} {\mathrm d Y_1}}{1- \hat C (Y_1, Y_3)}  \Biggr)+    \frac {\hat C_{\{1,3\}\{4\}}(Y_1, Y_3, Y_4)}{(1-\hat C(Y_1, Y_3))^2}\Biggl(1 - \frac{C_1(Y_1)^2\frac {\mathrm{d} X_1} {\mathrm d Y_1}}{1- \hat C (Y_1, Y_2)} \Biggr)\Biggr]\Biggr) .
\label{eq:hat-C3-13}
\end{align}
 
\subsubsection{Contribution for $C_{p-1}$ at order $p$}

 It is easy to verify that at any order $p$, the terms involving the highest possible cumulants - the cumulants of order $p-1$ - cancel among themselves in \eqref{eq:simplifications}, that is: 
\begin{align}
\label{eq:order_p-1_at-orderp}
\prod_{i=1}^p {\frac {\mathrm{d} Y_i} {\mathrm d X_i}} \frac1{C_1(Y_i)}  \sum_{1\le a<b\le p} \bar C_{p, \{a,b\}} \bigl(Y_1, \ldots, Y_p\bigr)+ \sum_{a\neq b} X_a \frac{\partial}{\partial X_a}\left(\prod_{i=1}^p {\frac {\mathrm{d} Y_i} {\mathrm d X_i}} \frac{C_{p-1}(\{Y_i\}_{i\neq b}) \mathring{C}_2(Y_a, Y_b)}{C_1(Y_a)^2\prod_{i\neq a} C_1(Y_i)}\right) = 0,\nonumber
\end{align}
 which is in fact nothing more than \eqref{eq:eq:order2-of-order3}. It follows from the following stronger fact.
\begin{proposition} 
\label{prop:order_p-1_at-orderp}
For any $1\le a<b\le p$, it is true that:
\be
\bar C_{p, \{a,b\}} \bigl(Y_1, \ldots, Y_p\bigr)+ Y_a \frac{\partial}{\partial Y_a}\left( {\frac {\mathrm{d} Y_a} {\mathrm d X_a}} \frac{C_{p-1}(\{Y_i\}_{i\neq b}) \mathring{C}_2(Y_a, Y_b)}{C_1(Y_a)^2}\right)\ +  \ a\leftrightarrow b = 0.
\ee
where we we recall that $\bar C_{p, \{a,b\}}$ is $\bar C_{p, \bar \pi},$ for $\bar\pi$ the partition of $\{1, \ldots, p\}$ whose blocks all contain a single element, apart from the block $\{a,b\}$. 
\end{proposition}

\proof The  generating function $\bar C_{p, \{a,b\}}$ has been computed in \eqref{eq:bar-C-ab} and \eqref{eq:bar-C-ab2}  and found to be:
\begin{align}
\label{eq:proof-order-p-1-p-1}
\bar C_{p, \{a,b\}} \bigl(Y_1, \ldots, Y_p\bigr) = Y_a Y_b  \frac{\partial}{\partial Y_a} \frac{\partial}{\partial Y_b} \left[ \frac{Y_b C_{p-1}(\{Y_i\}_{i\neq b}) - Y_a C_{p-1}(\{Y_i\}_{i\neq a})}{  Y_a C_1(Y_b) - Y_b C_1(Y_a)} \right].
\end{align}
Noticing that
\be
\label{eq:c2dd}
C_1(Y) - Y C_1'(Y)=  C_1(Y)^2  \frac {\mathrm{d} X} {\mathrm d Y},
\ee
we can carry out the derivatives, obtaining:
\begin{align}
\bar C_{p, \{a,b\}} \bigl(Y_1, \ldots, Y_p\bigr) = Y_a Y_b  \frac{\partial}{\partial Y_a}  \left[  \frac{Y_a C_1(Y_b)^2}{ ( Y_a C_1(Y_b) - Y_b C_1(Y_a))^2} \frac {\mathrm{d} X_b} {\mathrm d Y_b} C_{p-1}(\{Y_i\}_{i\neq b})\right] \ + \ a \leftrightarrow b. 
\end{align}
On the other hand, 
$$
Y_a \frac{\partial}{\partial Y_a}\Biggl({\frac {\mathrm{d} Y_a} {\mathrm d X_a}} \frac{C_{p-1}(\{Y_i\}_{i\neq b}) \mathring{C}_2(Y_a, Y_b)}{C_1(Y_a)^2}\Biggr) =  -  Y_a \frac{\partial}{\partial Y_a}\Biggl[  \frac{Y_aY_b C_1(Y_b)^2 \frac {\mathrm{d} X_b} {\mathrm d Y_b}}{ ( Y_a C_1(Y_b) - Y_b C_1(Y_a))^2}  C_{p-1}(\{Y_i\}_{i\neq b}) \Biggr],
$$
which concludes the proof\qed

\ 

While \eqref{eq:proof-order-p-1-p-1} was obtained in  \eqref{eq:bar-C-ab} as \eqref{eq:term-p-1-order-p-D}, that is, applying $\mathcal{D}$ to $\tilde C_{2}$ expressed in its logarithm form \eqref{eq:C2Tilde-log},  the relation of Prop.~\ref{prop:order_p-1_at-orderp} can be seen as applying $\mathcal{D}$ to $\tilde C_{2}$  in its formulation \eqref{eq:Tilde-C2-in-pat-ex} but without the $X's$, that is:
\be
\tilde C_2(Y_1, Y_2) = \frac{Y_1Y_2}{(Y_1 - Y_2)^2} \left(  1 - \frac{(C_1(Y_1) - Y_1 C_1'(Y_1))(C_1(Y_2) - Y_2 C_1'(Y_2))}{(1-\hat C(Y_1, Y_2))^2}\right) . 
\ee

\subsubsection{Fourth order}

We now look more closely at \eqref{eq:simplifications} at  order $p=4$. There are four kinds of trees in $\mathcal{G}_4$. The weight for the tree with one black vertex and four white vertices of valency one is $H_4 - C_4$, where $H_4$ is given by \eqref{eq:H4}. Recalling that $\mathring{C}_2$ is given by \eqref{eq:C2Circ}, the contribution to \eqref{eq:simplifications} is:
\be
\label{eq:contrib-order4-simplif}
\prod_{i=1}^4   {\frac {\mathrm{d} Y_i} {\mathrm d X_i}}\frac 1 {C(Y_i)}\Bigl( \tilde C_4 + \sum_{\mathrm{perm}} \bar C_{4, \{1\}\{2,3, 4\}}   +   \sum_{\mathrm{perm}} \bar C_{4, \{1,2\}\{3, 4\}} +   \sum_{\mathrm{perm}} \bar C_{4, \{1\}\{2\}\{3, 4\}}  +  \sum_{\mathrm{perm}} \bar C_{4, \{1\}\{2\}\{3, 4\}, T}\Bigr).
\ee
The second kind of trees has one black vertex of valency 3, one black vertex $\{a,b\}$ of valency 2, one white vertex of valency two, and the other white vertices are of valency one. The weight for such trees is of the form $H_3 H_2 - C_3 (C_2 + Y_a Y_b /(Y_a - Y_b)^2)$. Carrying out the simplifications, the weight is of the form $$ C_3 \mathring C_2 + \bar C_{3, \{1\}\{2,3\}} C_2 +  \tilde C_3 C_2 + \bar C_{3, \{1\}\{2,3\}} \tilde C_2+ \tilde C_3 \tilde C_2, $$
where all appropriate choices of labels and arguments must be taken into account. From left to right, the first term involves some $C_3$, the second some $C_2 C_2$, the third and fourth some $C_2$, and the last only some $C_1$.

The third kind of trees has three black vertices $\{a, b\}$, $\{b,c\}$, $\{c,d\}$ of valency 2, two white vertices of valency two, and two white vertices of valency one. The weight for the black vertices of such trees is of the form $H_2 H_2 H_2 - (C_2 + Y_a Y_b /(Y_a - Y_b)^2)(C_2 + Y_b Y_c /(Y_b - Y_c)^2)(C_2 + Y_c Y_d /(Y_c - Y_d)^2) $.  This can be recast as
$$
\sum_{\mathrm{perm}} C_2 C_2 \mathring C_2 + \sum_{\mathrm{perm}} C_2 \Bigl(\tilde C_2 \tilde C_2 - \frac{Y_a Y_b }{(Y_a - Y_b)^2}\frac{ Y_bY_c }{(Y_b - Y_c)^2}\Bigr) + \tilde C_2 \tilde C_2 \tilde C_2 -  \frac{Y_a Y_b }{(Y_a - Y_b)^2}\frac{ Y_bY_c }{(Y_b - Y_c)^2} \frac{ Y_cY_d }{(Y_c - Y_d)^2} . 
$$
The leftmost terms involve some $C_2 C_2$, the second terms some $C_2$, and the two remaining terms only some $C_1$. The last kind of trees has three black vertices  $\{a, b\}$, $\{a,c\}$, $\{a,d\}$ of valency two, but one white vertex of valency 3. The weights are just as above, but with adapted indices and arguments. 

\

As mentioned above, the simplification \eqref{eq:simplifications} should occur ``order-by-order'', that is here, we expect the terms corresponding to each item in the following list to cancel among each other (and coefficientwise for the various partial derivatives): 
\begin{enumerate}
\item The terms involving $C_3$ (this we have already shown in Prop.~\ref{prop:order_p-1_at-orderp}),
\item The terms of the form $C_2(Y_1, Y_2),C_2(Y_1, Y_3)$ with one repeated argument, 
\item The terms of the form  $C_2(Y_1, Y_2),C_2(Y_3, Y_4)$, 
\item The terms of the form  $C_2(Y_1, Y_2)$, 
\item The terms involving only the $C_1$.  
\end{enumerate}

In addition, we expect the simplifications to derive from \eqref{eq:eq:order1-of-order2} and \eqref{eq:eq:order1-of-order3}, which are the simplifications for the terms involving only the $C_1$ for the lower orders $\tilde C_2$ and $\tilde C_3$.

\

We detail below the terms in $C_2 C_2$. They are obtained from the following corrections: 
\be
\label{eq:single-out-C2C2}
 \bar C_{4, \{1\}\{2\}\{3, 4\}, T},\  \ \bar C_{3, \{1\}\{2,3\}} C_2,\  \textrm{ and } \  C_2 C_2 \mathring C_2 .\ee In addition to $C_2C_2$, these only involve  $\tilde C_2$ (and $C_1$), so we expect these simplifications to derive from applying differential operators such as $\mathcal{D}$ and $\mathcal{P}\circ \mathcal{D^\otimes} \mathcal{D^\otimes}$ to the identity \eqref{eq:eq:order1-of-order2}. So to summarize, we expect the simplifications corresponding to the terms \eqref{eq:single-out-C2C2} to occur among each other:\\[+1ex]
 $-\ $among terms deriving from for instance $\tilde C_2(Y_1, Y_4)$,\\[+1ex]
 $-\ $among terms involving of the form $C_2(Y_1, Y_2),C_2(Y_1, Y_3)$ with one repeated variable or of the form $C_2(Y_1, Y_2),C_2(Y_3, Y_4)$ with disjoint variables, coefficientwise in the various partial derivatives.

\begin{figure}[!h]
\centering
\includegraphics[scale=0.7]{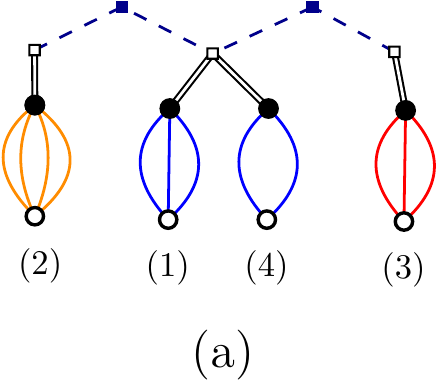}\hspace{0.7cm}\includegraphics[scale=0.7]{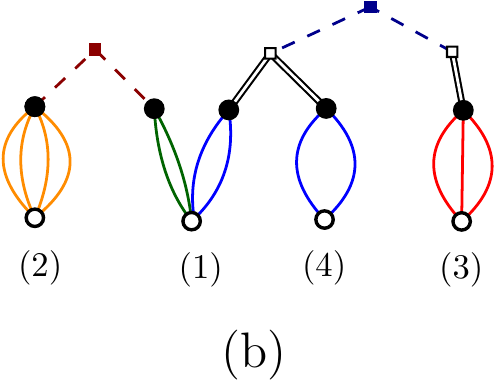}\hspace{0.7cm}\includegraphics[scale=0.7]{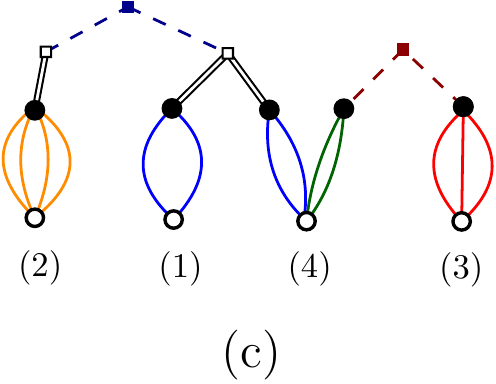}\hspace{0.7cm}\includegraphics[scale=0.7]{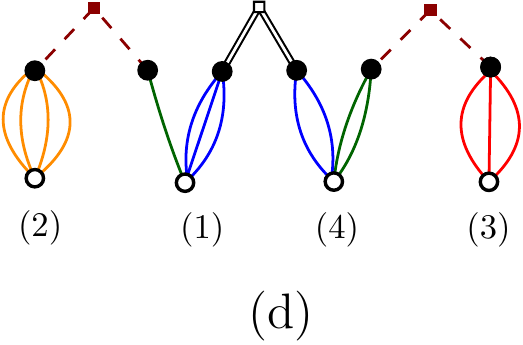}
\caption{Diagrammatic representation of the contributions in $C_2(Y_1, Y_2)$ and $C_2(Y_3, Y_4)$  to the simplifications at order 4 that derive from the expression of $\tilde C_2(Y_1, Y_4)$. 
The red vertices and edges represent direct insertions of $C_2$, while the thick white edges and the blue vertices and edges represent corrections as in Fig.~\ref{fig:Third-order} and Fig.~\ref{fig:Fourth-order}. 
}
\label{fig:Fourth-order-simpl}
\end{figure}

Focusing for instance on $\tilde C_2(Y_1, Y_4)$, $C_2(Y_1, Y_2)$ and $C_2(Y_3, Y_4)$, we get contributions from  the diagrams in Fig.~\ref{fig:Fourth-order-simpl}. 
We verify that 
these contributions to \eqref{eq:simplifications} simplify among each other, that is:
\begin{align}
\nonumber
&- Y_1 Y_4  \frac{\partial}{\partial Y_1}  \frac{\partial}{\partial Y_4} \Bigl(\frac{Y_1 Y_4}{(Y_1 - Y_4)^2}   \frac{C_2(Y_1, Y_2)C_2(Y_3, Y_4)}{(1 - \hat C (Y_1, Y_4))^2} \Bigr) \\&\nonumber
+ Y_1 \frac{\partial}{\partial Y_1}\biggl[\frac {\mathrm{d} Y_1} {\mathrm d X_1}\frac{C_2(Y_1, Y_2)}{C(Y_1)^2}Y_1 Y_4  \frac{\partial}{\partial Y_1} \frac{\partial}{\partial Y_4}\biggl(\frac{Y_1}{Y_4 - Y_1}\frac{C_2(Y_3, Y_4)}{1-\hat C(Y_1, Y_2)}\biggr)\biggr] \ +\  \left \{\begin{array}{ll}1\leftrightarrow 4\\ 2 \leftrightarrow 3\end{array}\right. 
\\&
 + Y_1 Y_4 \frac{\partial}{\partial Y_1}\frac{\partial}{\partial Y_4}\biggl[\frac {\mathrm{d} Y_1} {\mathrm d X_1}\frac {\mathrm{d} Y_4} {\mathrm d X_4}\frac{C_2(Y_1, Y_2)C_2(Y_3, Y_4)\mathring C_{2}(Y_1, Y_4)}{C(Y_1)^2C(Y_4)^2}\biggr] = 0, 
 \end{align}
where line by line from the top, the terms come from $\bar C_{4, \{2\}\{3\}\{1, 4\}, T}$ (Fig.~\ref{fig:Fourth-order-simpl}, (a)), $\bar C_{3, \{1,4\}\{3\}}$  (Fig.~\ref{fig:Fourth-order-simpl}, (b)) and $\bar C_{3, \{1,4\}\{2\}}$ for the bracket  (Fig.~\ref{fig:Fourth-order-simpl}, (c)), and $C_2 C_2 \mathring C_2$  (Fig.~\ref{fig:Fourth-order-simpl}, (d)). 

\begin{figure}[!h]
\centering
\includegraphics[scale=0.7]{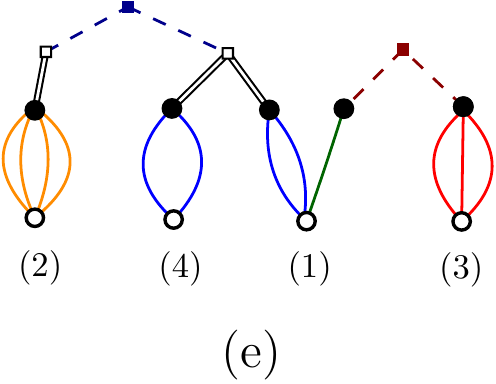}\hspace{0.5cm}\includegraphics[scale=0.7]{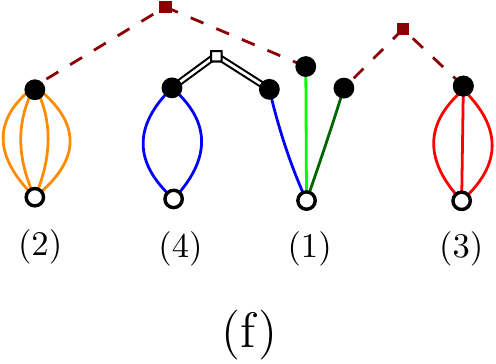}\hspace{0.5cm}\includegraphics[scale=0.7]{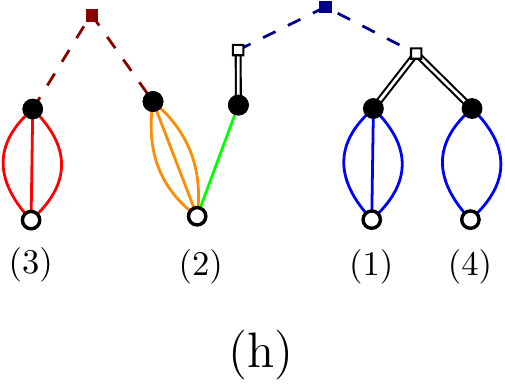}\hspace{0.5cm}\includegraphics[scale=0.7]{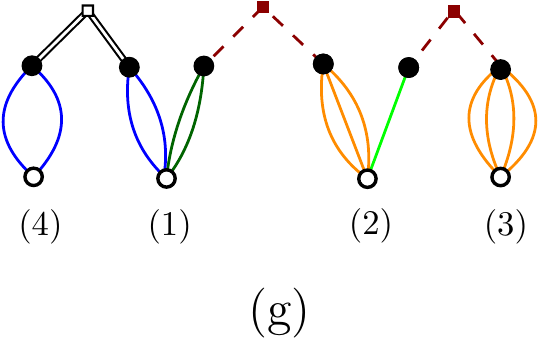}
\caption{Remaining contributions  in $C_2, C_2$  to the simplifications at order 4 that derive from the expression of $\tilde C_2(Y_1, Y_4)$. 
}
\label{fig:Fourth-order-simpl2}
\end{figure}

In the same way, the terms deriving from $\tilde C_2(Y_1, Y_4)$, and involving $C_2(Y_1, Y_2)$ and $C_2(Y_1, Y_3)$ cancel as:
\begin{align}
\nonumber
&Y_1 Y_4  \frac{\partial}{\partial Y_1}  \frac{\partial}{\partial Y_4} \Bigl(\frac{Y_4^2}{(Y_1 - Y_4)^2}   \frac{C_2(Y_1, Y_2)C_2(Y_1, Y_3)}{(1 - \hat C (Y_1, Y_4))^2} \Bigr) 
\\&+ Y_1 \frac{\partial}{\partial Y_1}\biggl[\frac {\mathrm{d} Y_1} {\mathrm d X_1}\frac{C_2(Y_1, Y_3)}{C(Y_1)^2}Y_1 Y_4  \frac{\partial}{\partial Y_1} \frac{\partial}{\partial Y_4}\biggl(\frac{Y_4}{Y_1 - Y_4}\frac{C_2(Y_1, Y_2)}{1-\hat C(Y_1, Y_4)}\biggr)\biggr] \quad + 2\leftrightarrow 3
\\&
+ X_1^2 \frac{\partial^2}{\partial X_1^2} \Bigl( \frac {\mathrm{d} Y_1} {\mathrm d X_1}\frac{C_2(Y_1, Y_2)C_2(Y_1, Y_3) \mathring C_2(Y_1, Y_4)}{C_1(Y_1)^3}\Bigr)  = 0, 
 \end{align}
 where line by line from the top, the terms come from $\bar C_{4, \{2\}\{3\}\{1, 4\}, T}$ (Fig.~\ref{fig:Fourth-order-simpl}, (a)), $\bar C_{3, \{1,4\}\{3\}}$  (Fig.~\ref{fig:Fourth-order-simpl}, (b)) and $\bar C_{3, \{1,4\}\{2\}}$ for the bracket  (Fig.~\ref{fig:Fourth-order-simpl2}, (e)), and $C_2 C_2 \mathring C_2$  (Fig.~\ref{fig:Fourth-order-simpl2}, (f)).
 
Finally, the remaining contributions in $C_2 C_2$  involve other $C_2 C_2$ for instance $C_2(Y_1, Y_2) C_2(Y_2, Y_3)$ and therefore also cancel among each other: 
\begin{align}
&+ Y_2 \frac{\partial}{\partial Y_2}\biggl[\frac {\mathrm{d} Y_2} {\mathrm d X_2}\frac{C_2(Y_2, Y_3)}{C(Y_2)^2}Y_1 Y_4  \frac{\partial}{\partial Y_1} \frac{\partial}{\partial Y_4}\biggl(\frac{Y_4}{Y_1 - Y_4}\frac{C_2(Y_1, Y_2)}{1-\hat C(Y_1, Y_4)}\biggr)\biggr]
 \\&\hspace{3cm} Y_1 Y_2 \frac{\partial}{\partial Y_1}\frac{\partial}{\partial Y_2}\biggl[\frac {\mathrm{d} Y_1} {\mathrm d X_1}\frac {\mathrm{d} Y_2} {\mathrm d X_2}\frac{C_2(Y_1, Y_2)C_2(Y_2, Y_3)\mathring C_{2}(Y_2, Y_4)}{C(Y_1)^2C(Y_2)^2}\biggr]  =0,
\end{align}
where line by line from the top, the terms come from  $\bar C_{3, \{1,4\}\{2\}}$ (Fig.~\ref{fig:Fourth-order-simpl2}, (g)), and $C_2 C_2 \mathring C_2$  (Fig.~\ref{fig:Fourth-order-simpl2}, (h)).
There are terms of this kind involving $C_2(Y_1, Y_2) C_2(Y_1, Y_3)$ but they derive from $\tilde C_2(Y_3, Y_4)$ or  $\tilde C_2(Y_3, Y_4)$ and cancel among each other.

In each case above,  the simplifications presented are satisfied coefficientwise for the various partial derivatives $ \frac{\partial}{\partial Y_1}C_2(Y_1, Y_2)C_2(Y_1, Y_3)$, $ \frac{\partial}{\partial Y_2}C_2(Y_1, Y_2)C_2(Y_1, Y_3)$, etc.

This includes all the contributions in $C_2, C_2$ by permutations of the indices. We formalize this in a proposition:
\begin{proposition}
All the terms contributing to \eqref{eq:simplifications} that involve $C_2, C_2$ cancel among each other.
\end{proposition}

We expect the same to occur among terms involving only one $C_2$, implying that \eqref{eq:identity} would be satisfied for $\tilde C_4$, and so on. This gives ground to the structure proposed at the beginning of Sec.~\ref{sub:simplifications}, but a more efficient approach is needed to on one hand verify  \eqref{eq:identity}, and on the other, verify that it implies all the other simplifications in \eqref{eq:simplifications} via the application of differential operators such as $\mathcal P \circ \bigl(\prod_{K\in \mathcal{I(T)}} \mathcal{D}^\otimes_{K}\bigr)$. 
This is left for a future study.

\end{document}